\documentclass[11pt]{article}
\usepackage{amsmath,amssymb}
\usepackage{graphicx}
\usepackage{epstopdf}
\newtheorem{theorem}{Theorem}[section]
\newtheorem{proposition}{Proposition}[section]
\newtheorem{corollary}{Corollary}[section]

\newcommand{\CQFD}{\nolinebreak\hfill\rule{2mm}{2mm}\medbreak\par}  
\newtheorem{Rem}{Remark}[section]  
  
\newtheorem{Example}{Example}[section]

\numberwithin{equation}{section}

\def\bbe{\mathbb E}

\def\bbp{\mathbb P}
\def\bbr{\mathbb R}

\def\udel{\Delta}

\def\udel{\Delta}
\def\del{\delta}
\def\call{{\cal L}}
\def\lam{\lambda}

\def\bbr{{\mathbb R}}
\def\bbe{{\mathbb E}}

\def\diag{\mathop{\rm diag}}
\def\udiag{\mathop{\rm Diag}}
\DeclareMathOperator{\textvar}{Var}

\begin{document}

\title{On the Limiting Shape of Young Diagrams Associated With Markov Random Words}
\author{Christian Houdr\'e \thanks{Georgia Institute of Technology,
School of Mathematics, Atlanta, Georgia, 30332-0160,
houdre@math.gatech.edu.  This work was partly supported by the grants \#246283 and 
\#524678 from 
the Simons Foundation and a Simons Foundation Fellowship grant \#267336.}
\and Trevis J.~Litherland\thanks{LexisNexis, Inc., Alpharetta, GA 30005; trevis.litherland@lexisnexisrisk.com} }

\maketitle

\vspace{0.5cm}

\begin{abstract}
\noindent Let $(X_n)_{n \ge 0}$ 
be an irreducible, aperiodic, homogeneous
Markov chain, with state space 
a totally ordered finite alphabet of size $m$.
Using combinatorial constructions and 
weak invariance principles, we obtain
the limiting shape of the associated RSK Young diagrams 
as a multidimensional Brownian functional.  Since the 
length of the top row of
the Young diagrams is also
the length of the longest weakly increasing
subsequences of $(X_k)_{1\le k \le n}$, 
the corresponding limiting law follows.
We relate our results to a conjecture of
Kuperberg by providing, under a cyclic condition,
a spectral characterization of the Markov
transition matrix precisely characterizing when
the limiting shape is the spectrum of the
$m \times m$ traceless GUE.  For each $m \ge 4$, 
this characterization identifies a proper, non-trivial 
class of cyclic transition matrices
producing such a limiting shape.  
However, for $m=3$, all cyclic Markov chains 
have such a limiting shape,
a fact previously only known for $m=2$. 
For $m$ arbitrary, we also study 
reversible Markov chains and obtain a characterization 
of symmetric Markov chains for which the limiting shape 
is the spectrum of the traceless GUE.  
To finish, we explore, in this general setting,
connections between various limiting laws and
spectra of Gaussian random matrices, focusing in particular
on the relationship between the terminal points of the 
Brownian motions, the diagonal terms of the random matrix, and
the scaling of its off-diagonal terms, a scaling we conjecture
to be a function of the spectrum of the covariance matrix 
governing the Brownian motion.
\end{abstract}





\section{Introduction}

The identification of the, properly centered and scaled, limiting law of $LI_n$, the length of the 
longest weakly increasing subsequences of a random word of length $n$,
whose letters are iid taking their values in a totally ordered $m$-letter alphabet,
was first obtained by Kerov~\cite{Ke}.   For example, 
in case the letters are uniformly distributed, the limiting law is the maximal eigenvalue of the 
$m\times m$ traceless Gaussian unitary ensemble (GUE).  
Moreover, \cite[Chap.~3, Sec.~3.4, Theorem 2]{Ke} showed that the whole normalized limiting shape of 
the pair of Robinson-Schensted-Knuth (RSK) Young diagrams (having necessarily $m$-row and of the same shape) 
associated with the random word, is the spectrum 
of the $m\times m$ traceless GUE.    
Since the length of the top row of the diagrams is the length of the longest weakly increasing 
subsequences of the random word, the maximal eigenvalue result is recovered.  (The asymptotic 
length result was rediscovered by Tracy and Widom~\cite{TW} and the asymptotic shape one 
by Johansson~\cite{Jo}, using orthogonal polynomials methods,  
positively answering a conjecture of Tracy and Widom.)  
In the non-uniform case, the corresponding
limiting distribution is described as 
that of the largest eigenvalue of one
of the diagonal blocks 
(corresponding to the highest probability)
in a direct sum of certain independent GUE matrices.  
The number and respective dimensions of these matrices
are determined by the multiplicities 
of the probabilities of choosing the
letters, and the direct sum is again 
subject to an overall zero-trace type of condition, and extensions to non-uniform 
letters were also obtained    
by Its, Tracy and Widom~\cite{ITW1, ITW2}.

The results of \cite{TW} were in part motivated by the corresponding results on 
the limiting law of the length of the longest strictly increasing subsequences 
of a uniform random permutation of $\{1, \dots, n\}$, which was first determined by 
Baik, Deift, and Johansson \cite{BDJ1} to be the Tracy-Widom distribution, 
describing the limiting behavior of
the largest eigenvalue of the GUE, 
as the size of the GUE goes to infinity.  
Further, Okounkov \cite{Ok}, and
Borodin, Okounkov, and Olshanski \cite{BOO},
as well as Johansson \cite{Jo}, also answered a conjecture
of Baik, Deift, and Johansson \cite{BDJ2,BDJ2Add} regarding the limiting
shape of the Young diagrams associated with a uniform random
permutation of $\{1,2,\dots,n\}$.
In particular, as $n$ grows without bound,
the lengths $\lambda_1, \lambda_2, \dots ,\lambda_k$
of the first $k$ rows of the Young diagrams,
appropriately centered and scaled, have the same limiting law 
as the $k$ largest eigenvalues of a $n \times n$ 
element of the GUE,
a result first proved, for $k=2$, in \cite{BDJ2,BDJ2Add}.

Outside of the iid framework, 
Kuperberg \cite{Ku} conjectured that
if the word is generated by an irreducible,
doubly-stochastic, cyclic Markov chain,
then the limiting distribution of the shape is still
that of the joint distribution of the eigenvalues
of a traceless $m \times m$ element of the GUE.  
For $m=2$, this was shown to be true 
by Chistyakov and G\"otze \cite{ChG} (see also \cite{HL2}),
who, in view of further simulations,
expressed doubts concerning the validity for $m \ge 4$.
For $m=3$, we will show that the conjecture
holds as well. 
However, for $m\ge4$, this is no longer the case.
Indeed, some, but not all, cyclic Markov chains 
lead to a limiting law as in the iid uniform case.

The precise class of homogeneous Markov chains 
with which Kuperberg's conjecture is concerned is more specific than
the ones we shall study.  The irreducibility of
the chain is a basic property we certainly must demand:
each letter has to occur at some point following the
occurrence of any given letter.
Moreover, the doubly-stochastic hypothesis
ensures that we have a uniform stationary distribution.
However, the cyclic (also called {\it circulant})
criterion, {\it i.e.}, that the Markov transition
matrix $P$ has entries satisfying $p_{i,j} = p_{i+1,j+1}$,
for $1 \le i,j \le m$ (where $m+1 = 1$), 
is more restrictive: cyclicity implies but
is not equivalent to $P$ being doubly stochastic.  Starting
from a free-probability perspective, Kuperberg was led
to introduce this cyclicity restriction via simulations \cite{Ku}
which appear to show that at least some irreducible, 
doubly-stochastic, {\it non-}cyclic Markov chains do
{\it not} produce such limiting behavior.  
Our approach to studying this problems is in part motivated by 
the identity in law, due to Baryshnikov~\cite{Ba} and Gravner, 
Tracy and Widom \cite{GTW}, between 
the maximal eigenvalue of an $m\times m$ element of the GUE and a 
certain maximal functional of 
$m$-dimensional standard 
Brownian motion originating, with Glynn and Whitt~\cite{GW}, in queueing theory.

The connections between queueing theory and problems on strings in a Markovian setting are
multifold but different from our own setting.  In particular, the work of Malyshev \cite{Ma6}
examines the evolution of strings driven by $r$-tuple transitions governed by natural
constraints that provide the basis for the existence of stationary measures and ergocidity.  Further 
classification of such strings into ergodic, null-recurrent, and transient cases
via the maximal eigenvalue of certain matrices is given in \cite{GaMaMePe5}.  Null recurrence is  
studied in yet greater detail in \cite{GaIaMa4}, where laws of large numbers and a central limit theorem for the number
of occurrences of a given state (element of the alphabet) are proved.  Other perspectives
using ideas from random walks in random environments are given in \cite{Ma3}.  Finally, generalizations
to random grammars on graphs, for which certain thermodynamic limits are shown to exist, can be found in \cite{Ma2}, with
additional work on random grammars on strings, using cluster expansion methods, being explored in \cite{Ma1}.
Our work differs from these efforts in that we focus on the more straightforward asymptotics of the Markovian evolution of strings by leveraging an elementary construction that ultimately 
yields CLT-type results.

The paper is organized in the following manner.
In Section $2$, we present a combinatorial
formulation of the $LI_n$
problem and so obtain a (discrete) functional of
combinatorial quantities which describes
the shape of the entire Young diagrams with $n$ cells,
along with a concise expression for the
associated asymptotic covariance structure.  
In Section $3$, we apply
Markovian Invariance Principles
to express the limiting shape of the Young 
diagrams as a functional of $m$-dimensional, $\it dependent$ 
Brownian motion 
for all irreducible, aperiodic, homogeneous
Markov chains (without the cyclic or even
the doubly-stochastic constraint).  
Using this functional we are then able to
answer Kuperberg's conjecture.
In Section $4$, we investigate, in
further detail, various symmetries exhibited
by the Brownian functional and in particular obtain,
for $m$ arbitrary, a necessary and sufficient condition
for a cyclic Markov chain to have
the same limiting law as in the iid uniform case.
Still for $m$ arbitrary, a related characterization 
is also obtained for reversible Markov chains 
with symmetric transition matrices, {\it i.e.,} in the 
doubly stochastic case.
In Section $5$, we further explore connections between the
various Brownian functionals obtained as limiting laws
and eigenvalues of random matrices.  In particular, we relate
the behavior of the terminal points of the Brownian motions to
the diagonal elements of certain Gaussian random matrices, conjecturing
that the off-diagonal terms scale according to some yet to be determined 
function of the spectrum of the covariance matrix governing the
Brownian motions. Finally, in Section $6$, we conclude with a brief discussion
of natural extensions and complements to some of the ideas
and results presented in the paper.

\

\noindent
{\bf Acknowledgments:} We gratefully thank the referee 
whose detailed reading 
and numerous comments greatly helped to improve this paper.

\section{Combinatorics}

The combinatorial development for the $m$-letter
alphabet, as obtained in \cite{HL},
resulted in $m-1$ partial sums $S_n^r$, $1 \le r \le m-1$.  
Using an even more straightforward
development which involves $m$ quantities instead,
we can obtain more symmetric expressions for $LI_n$.
This is done next, and will prove useful when studying
the shape of the whole RSK Young diagrams.\\

Let $a^r_k$ be the number of occurrences
of $\alpha_r$ among $(X_i)_{1\le i \le k}$.  
Each weakly increasing subsequence of $(X_i)_{1\le i \le k}$ consists simply of 
consecutive identical values, 
with these values forming a weakly increasing subsequence of $\alpha_r$.
Moreover, the number of occurrences of $\alpha_r\in \{\alpha_1,
\dots,\alpha_m\}$ among
$(X_i)_{k+1 \le i \le \ell}$,
where $1 \le k < \ell \le n$, is simply
$a^r_{\ell}-a^r_k$. The length of the longest weakly increasing subsequences 
of $X_1,X_2,\dots, X_n$ is thus given by

\begin{equation}\label{item1}
LI_n=\max_{\stackrel{\scriptstyle 0\le k_1\le\cdots}{\le k_{m-1}\le k_m=n}}
     [(a^1_{k_1}-a^1_0)+(a^2_{k_2}-a^2_{k_1})+\cdots +
     (a^m_n-a^m_{k_{m-1}})],
\end{equation}

\noindent {\it i.e.},

\begin{equation}\label{item2}
LI_n=\max_{\stackrel{\scriptstyle 0\le k_1\le\cdots}{\le k_{m-1}\le k_m=n}}
     [(a^1_{k_1}-a^2_{k_1})+(a^2_{k_2}-a^3_{k_2})+\cdots +
     (a^{m-1}_{k_{m-1}}-a^m_{k_{m-1}})+a^m_n],
\end{equation}

\noindent where $a^r_0=0$, $r= 1, \dots, m$.  

Moving beyond the purely combinatorial setting, 
assume that $(X_k)_{k \ge 0}$ is an infinite
sequence generated by an irreducible homogeneous Markov chain
having a stationary distribution $(\pi_1,\pi_2,\dots,\pi_m)$.
(For no $k \ge 0$ is the law of $X_k$ necessarily
assumed to be the stationary distribution.)
For each $1 \le r \le m$, 
set $T^r_k = a^r_k - \pi_r k$,
for $k \ge 1$, and $T^r_0 = 0$.
Beginning again with \eqref{item1},
we find that

\begin{align}\label{item6a}
LI_n &=\max_{\stackrel{\scriptstyle 0\le k_1\le\cdots}{\le k_{m-1}\le k_m=n}}
      \Bigl[(a^1_{k_1}-a^1_0)+(a^2_{k_2}-a^2_{k_1})+\cdots +
      (a^m_n-a^m_{k_{m-1}})\Bigr]\nonumber\\
     &= \max_{\stackrel{\scriptstyle 0\le k_1\le\cdots}{\le k_{m-1}\le k_m=n}}
      \Bigl[( (T^1_{k_1}+ \pi_1 k_1) - (T^1_{k_0}+ \pi_1 k_0) )+ ( (T^2_{k_2}+ \pi_2 k_2) - (T^2_{k_1}+ \pi_2 k_1) )\nonumber\\
     &\qquad +\cdots + ( (T^m_{k_m}+ \pi_m k_m) - (T^m_{k_{m-1}}+ \pi_m k_{m-1}) )\Bigr]\nonumber\\
     &= \max_{\stackrel{\scriptstyle 0\le k_1\le\cdots}{\le k_{m-1}\le k_m=n}}
      \Bigl[( T^1_{k_1} - T^1_{k_0} )+ ( T^2_{k_2} - T^2_{k_1} ) +\cdots + ( T^m_{k_m} - T^m_{k_{m-1}} ) \nonumber\\
     & \qquad + \pi_1( k_1 - k_0 ) + \pi_2( k_2 - k_1 ) + \cdots  + \pi_m( k_m - k_{m-1} )    \Bigr],
\end{align}
where $k_0=0$.  

\noindent Setting $\pi_{max} = \max\{\pi_1, \pi_2, \dots, \pi_m\}$, 
\eqref{item6a} becomes

\begin{align}\label{item6b}
LI_n - \pi_{max}n 
&= \max_{\stackrel{\scriptstyle 0=k_0 \le k_1\le\cdots}{\le k_{m-1}\le k_m=n}}
      \sum_{r=1}^{m} \bigl[ (T^r_{k_r} - T^r_{k_{r-1}})
     + (\pi_r -  \pi_{max}) ( k_r - k_{r-1} ) \bigr].
\end{align}

\noindent For a uniform stationary distribution, 
$\pi_{max} = \pi_r = 1/m$, for all $r$, and \eqref{item6b} 
simplifies to

\begin{equation}\label{item6c}
LI_n - \frac{n}{m} = \max_{\stackrel{\scriptstyle 0=k_0\le k_1\le\cdots}{\le k_{m-1}\le k_m=n}}
      \sum_{r=1}^{m}  (T^r_{k_r} - T^r_{k_{r-1}}).
\end{equation}

To introduce a random walk formalism into the picture, 
we next set, for $i = 1, \dots ,n$ and $r = 1,2, \dots ,m$,  

\begin{equation}\label{item6d}
W^r_i=	\begin{cases} 	1, &\text{if $X_i=\alpha_r,$}\\
			0, & \text{otherwise.}
	\end{cases}
\end{equation}

\noindent Clearly, $a^r_k = \sum^k_{i=1}W^r_i$,
and so $T^r_k = \sum^k_{i=1}(W^r_i - \pi_r)$,
for $1 \le r \le m$.

To understand the limiting law of \eqref{item6b} or
\eqref{item6c}, we must have a more precise description
of the underlying Markovian structure.  To that end, let 
$p_{r,s} = \bbp(X_{k+1} = \alpha_s|X_k = \alpha_r)$,
and let $P = (p_{r,s})_{1\leq r, s \leq m}$
be the associated Markov transition matrix.  In this setting,

$$ (\pi^{n+1}_1, \pi^{n+1}_2,\dots,\pi^{n+1}_m) = (\pi^{n}_1, \pi^{n}_2,\dots,\pi^{n}_m)P.$$

\noindent Moreover, as usual, let $p^{(j)}_{r,s}$, $j\ge 1$, 
denote the $j$-step transition
probability from $\alpha_r$ to $\alpha_s$; its associated transition
matrix is simply $P^j$.

Assume now that the law of $X_0$ is the stationary distribution.
Thus, by construction, $\bbe T^r_k = 0$ for all $1\le r \le m$ and $1 \le k \le n$, and
our primary task is to describe the covariance structure of
these random variables $T_k^r$.  First, for $k \ge 1$, 

\begin{align}\label{item6e}
\mbox{Var}T^r_k 
= \sum^k_{i=1}\mbox{Var}W^r_i + \sum^{k-1}_{i_1=1} \sum^{k}_{i_2=i_1+1} \mbox{Cov}(W^r_{i_1},W^r_{i_2}) 
+ \sum^k_{i_1=2} \sum^{i_1 - 1}_{i_2=1} \mbox{Cov}(W^r_{i_1},W^r_{i_2}).
\end{align}

\noindent Next, by stationarity, and since $W^r_i$ is, simply, a Bernoulli random variable
with parameter $\pi_r$, \eqref{item6e} becomes

{\allowdisplaybreaks
\begin{align}\label{item6f}
\mbox{Var}T^r_k 
&= \sum^k_{i=1}\mbox{Var}W^r_i + \sum^{k-1}_{i_1=1} \sum^{k}_{i_2=i_1+1} \mbox{Cov}(W^r_{0},W^r_{i_2-i_1}) 
+ \sum^k_{i_1=2} \sum^{i_1 - 1}_{i_2=1} \mbox{Cov}(W^r_{0},W^r_{i_1-i_2})\nonumber\\
&= k\pi_r(1-\pi_r) + \sum^{k-1}_{i_1=1} \sum^{k}_{i_2=i_1+1} (\pi_r p^{(i_2-i_1)}_{r,r} - \pi_r^2) 
+ \sum^k_{i_1=2} \sum^{i_1 - 1}_{i_2=1} (\pi_r p^{(i_1-i_2)}_{r,r} - \pi_r^2)\nonumber\\
&= k\pi_r -k^2\pi_r^2 +  \pi_r\sum^{k-1}_{i_1=1} \sum^{k}_{i_2=i_1+1}  e_r P^{i_2-i_1} e_r^T 
+ \pi_r\sum^k_{i_1=2} \sum^{i_1 - 1}_{i_2=1}  e_r P^{i_1-i_2} e_r^T,
\end{align}
}

\noindent where $e_r = (0,0,\dots, 0,1,0, \dots 0)$ 
is the $r^{th}$ standard basis vector of $\bbr^m$.  
Setting

\begin{equation}\label{item6fa}
Q_k = \sum^{k-1}_{i_1=1} \sum^{k}_{i_2=i_1+1}  P^{i_2-i_1} = \sum^{k-1}_{i=1}  (k-i) P^{i},
\end{equation}

\noindent we can rewrite \eqref{item6f} in the simple form

\begin{align}\label{item6g}
\mbox{Var } T^r_k &= k\pi_r - k^2\pi_r^2 + 2\pi_r e_r Q_k e_r^T .
\end{align}

Our description of the covariance structure 
can now be completed using the above results.  
For $r_1 \ne r_2$ and $k \ge 1$,

{\allowdisplaybreaks
\begin{align}\label{item6h}
\mbox{Cov}(T^{r_1}_k, T^{r_2}_k) 
&= \sum^k_{i=1}\mbox{Cov}(W^{r_1}_i,W^{r_2}_i) +\sum^{k-1}_{i_1=1} \sum^{k}_{i_2=i_1+1} \mbox{Cov}(W^{r_1}_{i_1},W^{r_2}_{i_2}) 
+ \sum^k_{i_1=2} \sum^{i_1 - 1}_{i_2=1} \mbox{Cov}(W^{r_1}_{i_1},W^{r_2}_{i_2})\nonumber\\
&= \sum^k_{i=1}\mbox{Cov}(W^{r_1}_i,W^{r_2}_i)  + \sum^{k-1}_{i_1=1} \sum^{k}_{i_2=i_1+1} \mbox{Cov}(W^{r_1}_{0},W^{r_2}_{i_2-i_1}) 
+ \sum^k_{i_1=2} \sum^{i_1 - 1}_{i_2=1} \mbox{Cov}(W^{r_2}_{0},W^{r_1}_{i_1-i_2})\nonumber\\
&= -k\pi_{r_1}\pi_{r_2} + \sum^{k-1}_{i_1=1} \sum^{k}_{i_2=i_1+1} (\pi_{r_1} p^{(i_2-i_1)}_{r_1,r_2} - \pi_{r_1}\pi_{r_2}) 
+ \sum^k_{i_1=2} \sum^{i_1 - 1}_{i_2=1} (\pi_{r_2} p^{(i_1-i_2)}_{r_2,r_1} - \pi_{r_1}\pi_{r_2})\nonumber\\
&=  -k^2\pi_{r_1}\pi_{r_2} +  \pi_{r_1}\sum^{k-1}_{i_1=1} \sum^{k}_{i_2=i_1+1}  e_{r_1} P^{i_2-i_1} e_{r_2}^T 
+ \pi_{r_2}\sum^k_{i_1=2} \sum^{i_1 - 1}_{i_2=1}  e_{r_2} P^{i_1-i_2} e_{r_1}^T\nonumber\\
&= -k^2\pi_{r_1}\pi_{r_2} + \pi_{r_1} e_{r_1} Q_k e_{r_2}^T + \pi_{r_2} e_{r_2} Q_k e_{r_1}^T.
\end{align}
}

\begin{Rem} Both \eqref{item6g} and \eqref{item6h} 
appear to be asymptotically quadratic in $k$.  However, 
since $Q_k  = \sum^{k-1}_{i=i}  (k-i) P^{i}$,
cancellations will show that
when the Markov chain is irreducible and aperiodic, 
the order of the variance is, in fact, linear in $k$. 
\end{Rem}

In order to further analyze the asymptotics of $Q_k$, 
we first examine the Jordan decomposition of $P$ for a very
general class of transition matrices.

\begin{proposition}\label{prop1}
Let $P$ be the $m \times m$ transition matrix of
an irreducible, aperiodic, homogeneous Markov chain.
Further, let $P = S^{-1}\Lambda S$ be the Jordan decomposition of $P$, 
where the rows
of $S$ consist of the generalized {\it left}-eigenvectors of $P$
with, moreover, the first row of $S$ being the
stationary distribution $(\pi_1,\pi_2,\dots,\pi_m)$,
and where $\Lambda = \text{diag}(1,\tilde{\Lambda})$, 
with $\tilde{\Lambda}$ being the $(m-1) \times (m-1)$
(lower) Jordan-block matrix.
Then the first column of $S^{-1}$
is $(1,1,\dots,1)^T$.
\end{proposition}

\noindent \begin{Proof}
Since $P = S^{-1}\Lambda S$, then 
$PS^{-1} = S^{-1}\Lambda$.  Denoting the first column of
$S^{-1}$ by $c_1$, we have
$Pc_1 = c_1$.  But since each row of $P$ sums up to $1$,
we see that $c_1 = (1,1,\dots,1)^T$
satisfies $Pc_1 = c_1$.  Moreover, $c_1$ must be
unique, up to normalization, since the irreducibility
of $P$ implies that the eigenvalue 
$\lambda_1 = 1$ has multiplicity $1$.
Finally, since the inner product of
the first row of $S$ and the first column of $S^{-1}$
is $1$, the correct normalization is indeed
$(1,1,\dots,1)^T$.\CQFD
\end{Proof}

Returning to $Q_k$, as given in \eqref{item6fa},
and using Proposition~\ref{prop1}, we then
obtain:

\begin{theorem}\label{thm2}
Let $(X_n)_{n \ge 0}$ be a sequence generated
by an $m$-letter, aperiodic, irreducible, homogeneous
Markov chain with state space
${\cal A}_m = \{\alpha_1 < \cdots < \alpha_m\}$,
transition matrix $P$,
and stationary distribution $(\pi_1,\pi_2,\dots,\pi_m)$.  
Let also the law of $X_0$ be the
stationary distribution.
Moreover,
for $1 \le r \le m$,
let $T^r_k = a^r_k - \pi_r k$,
for $k \ge 1$, and $T^r_0 = 0$, where 
$a^r_k$ is the number of occurrences
of $\alpha_r$ among $(X_i)_{1\le i \le k}$.
Then, for $1 \le r \le m$,

\begin{equation}\label{item6k}
\lim_{k \rightarrow \infty} \frac{\text{Var } T^r_k}{k} 
= \pi_r\left(1 + 2e_r S^{-1} D S e_r^T  \right),
\end{equation}

\noindent and for $r_1 \ne r_2$,

\begin{equation}\label{item6l}
\lim_{k \rightarrow \infty} \frac{\mbox{Cov}(T^{r_1}_k, T^{r_2}_k)}{k} = \pi_{r_1}e_{r_1} S^{-1} D S e_{r_2}^T  
+ \pi_{r_2}e_{r_2} S^{-1} D S e_{r_1}^T ,
\end{equation}

\noindent where $P = S^{-1}\Lambda S = S^{-1} \text{diag}(1,\tilde{\Lambda}) S$ 
is the Jordan decomposition
of $P$ in Proposition~\ref{prop1}, and 
$D = \text{diag}(-1/2,\tilde{\Lambda}(I - \tilde{\Lambda})^{-1} )$.
That is, the asymptotic covariance matrix of
$(T^1_k,T^2_k,\dots,T^m_k)$ is given by

\begin{equation}\label{item6la}
 \Sigma = \Pi + \Pi(S^{-1} D S) + (S^{-1} D S)^T \Pi,
\end{equation}

\noindent where $ \Pi = \text{diag}(\pi_1,\pi_2,\dots,\pi_m)$.
\end{theorem}

\noindent \begin{Proof} Beginning with \eqref{item6fa},
we use the Jordan decomposition of $P$ 
in Proposition~\ref{prop1} to find that

{\allowdisplaybreaks
\begin{align}\label{itm6m}
Q_k = \sum^{k-1}_{i=1} (k-i) (S^{-1}\Lambda S)^{i} = S^{-1} \left( \sum^{k-1}_{i=1} (k-i) \Lambda ^{i} \right) S 
= S^{-1} \text{diag}\left( h(1, k), \sum^{k-1}_{i=1} (k-i) \tilde{\Lambda}^{i} \right) S,
\end{align}
}

\noindent where 
$h(\lambda, k) := \sum_{i=1}^{k-1} (k-i)\lambda^i$.
Now $h(1, k) = k(k-1)/2$ is quadratic in $k$, while for $\lambda \ne 1$,

$$h(\lambda, k) = k \frac{\lambda}{(1-\lambda)} + \frac{\lambda(\lambda^{k} - 1)}{(1-\lambda)^2},$$

\noindent so that $h(\lambda, k)$ is asymptotically linear in $k$ for $|\lambda| < 1$. 
By the aperiodicity of $P$, all 
the diagonal elements of the lower-triangular matrix 
$\tilde{\Lambda}$ are less than unity in modulus.
We can thus write the matrix sum in \eqref{itm6m} as

$$h(\tilde{\Lambda}, k) := k \tilde{\Lambda}(I-\tilde{\Lambda})^{-1} + o(k).$$

\noindent The matrix $Q_k$ can now be expressed as the (asymptotic) sum
of terms which are, respectively,
quadratic and linear in $k$.  Recalling, moreover,
that the first row of $S$ contains the
stationary distribution, and that the first
column of $S^{-1}$ is $(1,1,\dots,1)^T$, we have

{\allowdisplaybreaks
\begin{align}\label{itm6n}
Q_k &= S^{-1} \text{diag}( h(1, k), h(\tilde{\Lambda}, k) ) S,\nonumber\\
    &= \frac{k^2}{2} S^{-1} \text{diag}(  1, 0, \dots, 0) S  
    + k S^{-1} \text{diag}\left( -\frac12, \tilde{\Lambda}(I-\tilde{\Lambda})^{-1}\right) S + o(k)\nonumber\\
    &= \frac{k^2}{2} 
       \begin{pmatrix} 
          \pi_1 & \pi_2 & \cdots & \pi_m \\
	  \pi_1 & \pi_2 & \cdots & \pi_m \\
	  \vdots & \vdots & \cdots & \vdots\\
	  \pi_1 & \pi_2 & \cdots & \pi_m 
       \end{pmatrix} + k S^{-1} D S + o(k).
\end{align}
}

Starting with the variance in \eqref{item6g},
we next find that, for each $1 \le r \le m$,

\begin{align}\label{item6o}
\mbox{Var }T^r_k &= k\pi_r - k^2\pi_r^2 + 2\pi_r e_r Q_k e_r^T\nonumber\\
&= k\pi_r - k^2\pi_r^2 + 2\pi_r\left( \frac{k^2}{2} \pi_r  + k e_r S^{-1} D S e_r^T\right) + o(k)\nonumber\\
&= k\pi_r\left(1 + 2e_r S^{-1} D S e_r^T  \right) + o(k),
\end{align}

\noindent from which the asymptotic result \eqref{item6k}
follows immediately.

An identical development shows that, for $r_1 \ne r_2$,
\eqref{item6h} simplifies to

\begin{align}\label{item6q}
\mbox{Cov}(T^{r_1}_k, T^{r_2}_k) &= -k^2\pi_{r_1}\pi_{r_2} + \pi_{r_1} e_{r_1} Q_k e_{r_2}^T + \pi_{r_2} e_{r_2} Q_k e_{r_1}^T\nonumber\\
&= -k^2\pi_{r_1}\pi_{r_2} + \pi_{r_1}\left(\frac{k^2}{2} \pi_{r_2}  + k e_{r_1} S^{-1} D S e_{r_2}^T\right) \nonumber\\
&\qquad +  \pi_{r_2}\left(\frac{k^2}{2} \pi_{r_1}  + k e_{r_2} S^{-1} D S e_{r_1}^T\right) + o(k)\nonumber\\
&= k\left(\pi_{r_1}e_{r_1} S^{-1} D S e_{r_2}^T  + \pi_{r_2}e_{r_2} S^{-1} D S e_{r_1}^T\right) + o(k),\nonumber\\
\end{align}

\noindent from which the asymptotic result \eqref{item6l}
follows, and so does \eqref{item6la}.\CQFD
\end{Proof}

\begin{Rem}\label{Remiidnonunif}  
To see that \eqref{item6k} and \eqref{item6l} both
recover the covariance results for the iid case investigated by the
authors in \cite{HL}, let $P$ be the transition matrix whose
rows each consist of the stationary distribution
$(\pi_1,\pi_2,\dots,\pi_m)$.
In this case $\lambda_2 = \cdots = \lambda_m = 0$,
and so $D = \text{diag}( -1/2, 0,\dots,0)$. Hence,

\begin{align*}
e_{r_1} S^{-1} D S e_{r_2}^T
&= \left(1,*,\dots,*\right)D\left(\pi_{r_2},*,\dots,*\right)^T = -\frac{\pi_{r_2}}{2},
\end{align*}

\noindent for all $r_1$ and $r_2$, and so, for each $r$,

$$\lim_{k \rightarrow \infty} \frac{\mbox{Var }T^r_k}{k} = \pi_r\left(1 + 2\left(-\frac{\pi_r}{2}\right) \right) = \pi_r(1-\pi_r),$$

\noindent while, for $r_1 \ne r_2$,

$$\lim_{k \rightarrow \infty} \frac{\mbox{Cov}(T^{r_1}_k, T^{r_2}_k)}{k} 
= \pi_{r_1}\left( -\frac{\pi_{r_2}}{2} \right) + \pi_{r_2}\left( -\frac{\pi_{r_1}}{2} \right)     = - \pi_{r_1}\pi_{r_2}.$$

\noindent  In the uniform iid case, 
$\pi_r = 1/m$, for all $1 \le r \le m$.  
Hence, for $r_1 \ne r_2$, the 
asymptotic correlation between $T^{r_1}_k$ and $T^{r_2}_k$
is given by $(-1/m^2)/((1/m)(1-1/m)) = -1/(m-1)$, so that the covariance matrix
is indeed the permutation-symmetric one obtained in \cite{HL}.  
There is, moreover, another Brownian functional
representation for the iid uniform case in \cite{HL} in which the Brownian
motions have a tridiagonal covariance matrix.
\end{Rem}

\section{The Limiting Shape of the Young Diagrams}
Thus far, our results have centered on $LI_n$ alone, 
essentially ignoring the larger
question of the structure of the entire Young diagrams.
The present section extends
the combinatorial development of the previous section to answer the
question of the limiting shape of the Young diagrams.
Our first result in this direction is a purely combinatorial
expression generalizing \eqref{item1}.

\begin{theorem}\label{thm3}
Let $R^1_n, R^2_n, \dots, R^r_n$ be the lengths of the first $1 \le r \le m$ rows
of the RSK Young diagrams generated by the sequence $(X_k)_{1\le k \le n}$
whose elements belong to a totally ordered alphabet 
${\cal A}_m = \{\alpha_1 < \cdots < \alpha_m\}$.
Then, for each $1 \le r \le m$, the sum of the lengths of the first
$r$ rows of the Young diagrams is given by

\begin{equation}\label{item7e}
\sum_{j=1}^r R^j_n =
  \max_{k_{j,\ell} \in J_{r,m}} \sum_{j=1}^r \sum_{\ell=j}^{m-r+j} \left( a^{\ell}_{k_{j,\ell}} - a^{\ell}_{k_{j,\ell-1}} \right),
\end{equation}

\noindent where  
$J_{r,m} = \{(k_{j,\ell}, 1 \le j \le r, 0 \le \ell \le m): 
k_{j,\ell} = 0,  0 \le \ell \le j-1; 
k_{j,\ell} = n, m-r+j \le \ell \le m; 
k_{j,\ell-1} \le k_{j,\ell}, 1 \le \ell \le m; 
k_{j+1,\ell+1} \le k_{j,\ell}, 1 \le j \le r-1, 0 \le \ell \le m-1
\}$,
and where $a^{\ell}_k$ is the number of occurrences of
$\alpha_{\ell}$ among $(X_i)_{1\le i \le k}$.
\end{theorem}

\noindent \begin{Proof}
Recall that the sum of the lengths of the first $r$ rows of
the RSK Young diagrams generated by a sequence $(X_k)_{1\le k \le n}$,
whose letters arise from an $m$-letter alphabet,
has an interpretation in terms of the length of certain weakly increasing sequences.
Indeed, the sum $R^1_n + R^2_n + \cdots + R^r_n$ is equal to the maximum
sum of the lengths of $r$ disjoint, weakly increasing subsequences
of $(X_k)_{1\le k \le n}$, where by {\it disjoint} it is meant
that each element of $(X_k)_{1\le k \le n}$
occurs in at most one of the $r$ subsequences.
(See Lemma 1 of Section 3.2 in \cite{Fu}).
More general results of this sort, involving partial orderings of the
alphabet and associated antichains, are known
as Greene's Theorem \cite{Greene}.  However, such results are not
enough for our purpose.  Below we need a different
way of reconstructing disjoint subsequences.

We begin by examining
an arbitrary collection of $r$ disjoint, weakly increasing subsequences
of $(X_k)_{1\le k \le n}$, and show that we can always map
these $r$ subsequences onto another collection of 
$r$ disjoint, weakly increasing subsequences
whose properties will be amenable to our combinatorial analysis.

Specifically, with the number of rows $r$ fixed, suppose that,
for each $1 \le j \le r$, we have a weakly increasing subsequence 
$(X_{k_\ell^j}^j)_{1\le \ell \le n_j}$ of length $n_j \le n$,
and that the $r$ subsequences are disjoint.

We first construct the new subsequence 
$(\tilde{X}_{\tilde{k}_{\ell}^{1}}^1)_{1 \le \ell \le \tilde{n}_1}$
as follows.  First, place all $\alpha_1$s occurring among the $r$
original subsequences into 
$(\tilde{X}_{\tilde{k}_{\ell}^{1}}^1)_{1 \le \ell \le \tilde{n}_1}$,
if there are any.  If the last $\alpha_1$ occurs at the $n^{th}$
index, then 
$(\tilde{X}_{\tilde{k}_{\ell}^{1}}^1)_{1 \le \ell \le \tilde{n}_1}$
is complete.
Otherwise, place all $\alpha_2$s which occur after the
final $\alpha_1$ into 
$(\tilde{X}_{\tilde{k}_{\ell}^{1}}^1)_{1 \le \ell \le \tilde{n}_1}$,
if there are any.  If the last $\alpha_2$ occurs at the $n^{th}$
index, then 
$(\tilde{X}_{\tilde{k}_{\ell}^{1}}^1)_{1 \le \ell \le \tilde{n}_1}$
is complete.  Otherwise, continue adding, successively,
$\alpha_3, \dots, \alpha_{m-r+1}$ in the same manner.  Thus,
$(\tilde{X}_{\tilde{k}_{\ell}^{1}}^1)_{1 \le \ell \le \tilde{n}_1}$
consists of a weakly increasing sequence 
of length $\tilde{n}_1$
having values in
$\{ \alpha_1, \dots, \alpha_{m-r+1} \}$.

Next, we construct the new subsequence
$(\tilde{X}_{\tilde{k}_{\ell}^{2}}^2)_{1 \le \ell \le \tilde{n}_2}$
similarly.
By considering only those letters among the $r$ original
subsequences which have not already been moved to the first new subsequence,
start with the smallest available letter, $\alpha_2$,
and continue adding, successively, $\alpha_3, \dots, \alpha_{m-r+2}$.
Note that, crucially, all $\alpha_2$s added to 
$(\tilde{X}_{\tilde{k}_{\ell}^{2}}^2)_{1 \le \ell \le \tilde{n}_2}$
occur before the last index at which $\alpha_1$ was added to
the first subsequence.  More generally, each $\alpha_j$, 
$2 \le j \le m-r+2$, added to 
$(\tilde{X}_{\tilde{k}_{\ell}^{2}}^2)_{1 \le \ell \le \tilde{n}_2}$
occurs before the last $\alpha_{j-1}$ was added to the first subsequence.
Thus, $(\tilde{X}_{\tilde{k}_{\ell}^{2}}^2)_{1 \le \ell \le \tilde{n}_2}$
consists of a weakly increasing subsequence of length $\tilde{n}_2$
having values in $\{ \alpha_2, \dots, \alpha_{m-r+2} \}$.

The construction of 
$(\tilde{X}_{\tilde{k}_{\ell}^{j}}^j)_{1 \le \ell \le \tilde{n}_j}$,
for $3 \le j \le r$, continues in the same manner, 
with 
$(\tilde{X}_{\tilde{k}_{\ell}^{j}}^j)_{1 \le \ell \le \tilde{n}_j}$,
constructed from among the entries of the $r$ original subsequences
which were not moved into any of the first $j-1$ new subsequences,
so that
$(\tilde{X}_{\tilde{k}_{\ell}^{j}}^j)_{1 \le \ell \le \tilde{n}_j}$,
consists of a weakly increasing sequence 
of length $\tilde{n}_j$
having values in
$\{ \alpha_j, \dots, \alpha_{m-r+j} \}$.  It is possible that beyond some
$j \ge 2$ the new subsequences may be empty.

We claim that, indeed, the construction of the $r^{th}$ new subsequence
exhausts the set of available entries.  Indeed, without loss of generality,
assume that after we have created the $(r-1)^{th}$ new subsequence,
the set of available entries is non-empty,
and designate the location of the final $\alpha_{\ell}$ to be included
in the $j^{th}$ new subsequence by $k_{j,\ell}$, for
$1 \le j \le r$ and $1 \le \ell \le m$, with the convention that 
$k_{j, m-r+j} = n$.  (If no $\alpha_{\ell}$ was
available for inclusion, set $k_{j,\ell} = k_{j,\ell-1}$,
where $k_{j,0} = 0,$ for all $1 \le j \le r$.)
Clearly, all $\alpha_1,\alpha_2,\dots,\alpha_{r-1}$ have been
included in the first $r-1$ new subsequences.
If $r=m$, we are done: simply put the remaining $\alpha_r$s
into the $r^{th}$ new subsequence. 
If $r < m$, we may still ask whether there was,
for some $r+1 \le \ell \le m$, 
an $\alpha_{\ell}$
from among the available entries
which occurred before $k_{r,\ell-1}$.
Assume that there is such an $\alpha_{\ell}$.
Now by construction, $k_{j+1,\ell-r+j} \le k_{j,\ell-r+j-1},$
for $1 \le j \le r-1$.  Hence, there exist letters 
$\alpha_{j_1} < \alpha_{j_2} < \dots < \alpha_{j_r} \le \alpha_{\ell-1}$ 
among the original subsequences which occurred at or after $k_{r,\ell-1}$, and,
moreover, each letter must come from a different subsequence.
But since each original subsequence was weakly increasing, none
of them could have contained an $\alpha_{\ell}$ before
$k_{r,\ell-1}$, and we have a contradiction.

To better understand this construction, consider
the first row of Figure~\ref{figa},
which shows an initial sequence of
length $n=12$, with $m=4$ letters, 
broken into $r=3$ disjoint, weakly increasing subsequences 
of lengths $n_1=3, n_2=4$, and $n_3=3$, and so
with total length $10$.
The final three rows of the diagram show the results
of the operations described above, producing $3$ new
weakly increasing subsequences of length
$\tilde{n}_1=4, \tilde{n}_2=3$, and $\tilde{n}_3=3$.

\begin{figure}
  \begin{center}
   \includegraphics[width=0.8\textwidth]{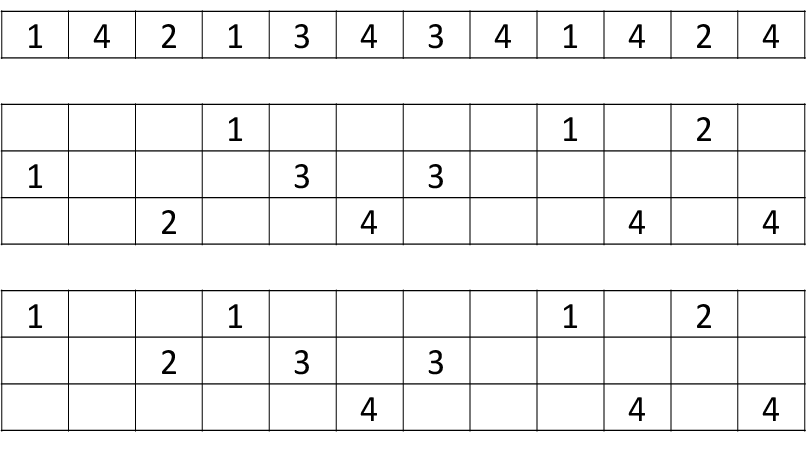}
  \end{center}
  \caption{Transformation of $r=3$ subsequences.}
  \label{figa}
\end{figure}

Hence, if we wish to find $r$ disjoint, weakly increasing subsequences
whose length sum is maximal, it suffices to consider only
those disjoint, weakly increasing subsequences for which
the final occurrence of the letter $\alpha_{\ell}$
in the subsequence $i$ happens after the final occurrence
in the subsequence $j$, whenever $i < j$.
Because such ranges do not overlap, if we wish to count
the number of $\alpha_{\ell}$s in the $j^{th}$ subsequence,
it suffices to simply count the number of $\alpha_s$s 
in $(X_k)_{1\le k \le n}$ over that range.

Indeed, returning to the fundamental combinatorial
objects of our development, the $a^j_k$, we see that
since $a^j_{\ell} - a^j_k$ counts the number
of $\alpha_j$s in the range $\ell + 1, \dots, k$,
we can describe the valid index ranges over which to
search for the maximal sum as
$J_{r,m} = \{(k_{j,\ell}, 1 \le j \le r, 0 \le \ell \le m): 
k_{j,\ell} = 0, 0 \le \ell \le j-1; 
k_{j,\ell} = n, m-r+j \le \ell \le m; 
k_{j,\ell-1} \le k_{j,\ell}, 1 \le \ell \le m; 
k_{j+1,\ell+1} \le k_{j,\ell}, 1 \le j \le r-1, 0 \le \ell \le m-1
\}.$
The constraints on the $k_{j,\ell}$ follow simply from the
fact that each subsequence is weakly increasing and that, moreover,
the intervals associated with a given letter do not overlap.
Figure~\ref{figb} indicates the relative positions of
each range, for $r=4$ and $m=7$.

\begin{figure}
  \begin{center}
    \includegraphics[width=0.8\textwidth]{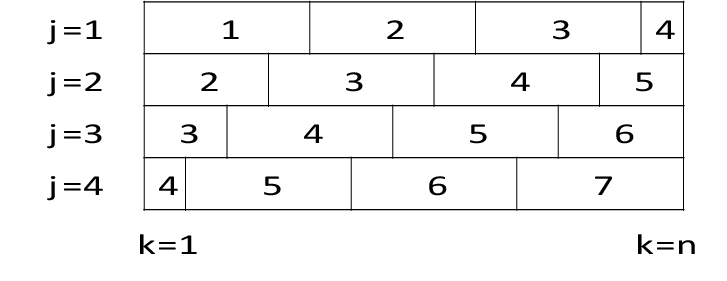}
  \end{center}
  \caption{Schematic diagram of $J_{r,m}$, for $r=4, m=7$.}
  \label{figb}
\end{figure}

Since the first possible letter of each subsequence grows
from $\alpha_1$ to $\alpha_r$, and the last possible letter grows
from $\alpha_{m-r+1}$ to $\alpha_m$, the result is proved.\CQFD
\end{Proof}

\begin{Rem}  The combinatorial representation obtained 
in the previous theorem has a 
known flavor, but is usually written in a different way.  Indeed, if $\lambda$ is the RSK 
shape image of the array of non-negative integers $\{b_{i,\ell}\}_{1\le i \le n, 1\le \ell \le m}$, then 
$\lambda_1 + \cdots + \lambda_r = \max_{\pi \in \Pi_r}\sum_{(i,\ell) \in \pi}b_{i,\ell}$, 
where $\Pi_r$ is the set of all collections of $r$ non-intersecting up/right lattice paths from 
$\{(1,1), \dots, (1,r)\}$ to $\{(n,m-r+1), \dots, (n,m)\}$.  
In our notation, $b_{i,\ell} = {\bf 1}_{X_i=\ell}$ and 
$k_{j, \ell}$ is the largest horizontal coordinate of the $j$-th path at level $\ell$.  
(See, e.g., \cite{HX} and the references therein, for details and a more complete bibliography.)

\end{Rem}

We are now ready to apply our asymptotic covariance results 
(Theorem \ref{thm2}), along 
with a Brownian sample-path approximation,
to the combinatorial expression \eqref{item7e}, 
and so obtain a Brownian functional
expression for the limiting shape of the Young diagrams for
all irreducible, aperiodic, homogeneous Markov chains.

Indeed, for each $1 \le r \le m$, let
the sum of the first $r$ rows of the Young diagrams
be given by

\begin{equation}\label{item7k}
V^r_n := \sum_{j=1}^r R^j_n =
  \max_{k_{j,\ell} \in J_{r,m}} \sum_{j=1}^r 
  \sum_{\ell=j}^{m-r+j} \left( a^{\ell}_{k_{j,\ell}} - a^{\ell}_{k_{j,\ell-1}} \right),
\end{equation}

\noindent where the index set $J_{r,m}$ is defined 
as in Theorem \ref{thm3}.  Define as before 
$T_k^r = \sum^k_{i=1}(W^r_i - \pi_r) = a_k^r - \pi_r k$, 
and so rewrite \eqref{item7k} as

\begin{align}\label{item7l}
V^r_n 
&=  \max_{k_{j,\ell} \in J_{r,m}} \sum_{j=1}^r \sum_{\ell=j}^{m-r+j} \left( \left(T^{\ell}_{k_{j,\ell}}   + \pi_{\ell}k_{j,\ell}  \right)
                                                                       -    \left(T^{\ell}_{k_{j,\ell-1}} + \pi_{\ell}k_{j,\ell-1}\right) \right)\nonumber\\
&=  \max_{k_{j,\ell} \in J_{r,m}}  \sum_{j=1}^r \sum_{\ell=j}^{m-r+j} \left(\left( T^{\ell}_{k_{j,\ell}}  - T^{\ell}_{k_{j,\ell-1}}\right)
                                       +  \pi_{\ell} \left(k_{j,\ell}   -           k_{j,\ell-1}\right)\right).
\end{align}

Next, let $\tau$ be a permutation of the indices $1,2,\dots,m$
such that $\pi_{\tau(1)} \ge \pi_{\tau(2)} \ge \cdots \ge \pi_{\tau(m)} > 0$.
Moreover, we demand that if $\pi_{\tau(i)} = \pi_{\tau(j)}$ for $i < j$,
then $\tau(i) < \tau(j)$.  (The permutation so defined is thus unique.)
Let $\nu_r = \sum_{j=1}^r \pi_{\tau(j)}$ be
the sum of the $r$ largest values among $\pi_1,\pi_2,\dots,\pi_m$.
We obtain, below, the limiting distribution of
$(V^r_n - \nu_r n)/\sqrt{n}$ 
as a Brownian functional.

To introduce Brownian sample-path approximations,
and for each $1 \le r \le m$, we first
define the asymptotic variance of 
$T^r_n$ as in \eqref{item6k}, by

\begin{equation}\label{item7n}
\sigma_r^2 := \lim_{n \rightarrow \infty} \frac{\text{Var } T^r_n}{n} = e_{r} \Sigma e_{r}^T ,
\end{equation}

\noindent and, for $r_1 \ne r_2$,
the asymptotic covariance of
$T^{r_1}_n$ and $T^{r_2}_n$ by

\begin{align}\label{item7o}
\sigma_{r_1,r_2} 
&:= \lim_{n \rightarrow \infty} \frac{\mbox{Cov}(T^{r_1}_n, T^{r_2}_n)}{n} = e_{r_1} \Sigma e_{r_2}^T ,
\end{align}

\noindent where $\Sigma$ is the covariance matrix of 
Theorem~\ref{thm2} associated with the transition matrix $P$.
For each $1 \le r \le m$,
we then let

\begin{equation}\label{item7p}
\hat B^r_n(t)=\frac{T^r_{[nt]} + (nt-[nt])(W^r_{[nt]+1} - \pi_r)}{\sigma_r\sqrt{n}},
\end{equation}

\noindent for $0 \le t \le 1$. This rescaling of
$[0,n]$ to $[0,1]$ calls for us to define a new
parameter set over which we will maximize a functional
arising from the expressions in \eqref{item7p}.
Indeed, for any positive integers $s$ and $d$, 
with $s \le d$, define
the set

\begin{align*}
I_{s,d}  = \Bigl\{(t_{j,\ell}, 1 \le j \le s, 0 \le \ell \le d): &t_{j,\ell} = 0, 0 \le \ell \le j-1;  
t_{j,\ell} = 1,  d-s+j \le \ell \le d;\nonumber\\
							    &t_{j,\ell-1} \le t_{j,\ell}, 1 \le \ell \le d; t_{j+1,\ell+1} \le t_{j,\ell}, 1 \le j \le s-1,     0 \le \ell \le d-1\Bigr\}.
\end{align*}

\noindent Note that the constraints 
$t_{j,j-1} = 0$ and $t_{j,d-s+j} = 1$, 
for $1 \le j \le s$,
force many of the $t_{j,\ell}$ to be zero or one.
We will denote the $s \times (d+1)$-tuple 
elements of $I_{s,d}$, by $(t_{.,.})$.
The structure of $I_{s,d}$ for $s=4$ and $d=7$ mirrors that of $J_{r,n}$ shown in Figure~\ref{figb}, but with the horizontal range now normalized to $[0,1]$.  The locations of $t_{j,\ell}$
would then be indicated by the horizontal lines within the diagram.


With this notation, we can write \eqref{item7l} as

\begin{align}\label{item7q}
\frac{V^r_n - n\nu_r }{\sqrt{n}}
&=  \max_{(t_{.,.}) \in I_{r,m}} 
    \Bigl\{      \sum_{j=1}^r \sum_{\ell=j}^{m-r+j} \sigma_{\ell}  \left( \hat B^{\ell}_n(t_{j,\ell})  - \hat B^{\ell}_n(t_{j,\ell-1}) \right) 
+ \sum_{j=1}^r \sum_{\ell=j}^{m-r+j}                \sqrt{n}(\pi_{\ell} - \pi_{\tau(j)})   \left( t_{j,\ell}  - t_{j,\ell-1} \right)      \Bigr\},  
\end{align}

\noindent since the piecewise linear nature of the expression being maximized requires the maximum to occur at parameter values of the form 
$t_{j,\ell} = k/n$, for $k = 0,1,...,n$.

Our analysis of \eqref{item7q} will yield the
following theorem, whose proof we defer to the
conclusion of the section.
This theorem gives, in particular, a full characterization 
of the limiting shape of the Young diagrams
in the non-uniform iid case.\\

\begin{theorem}\label{thm4}
Let $(X_n)_{n \ge 0}$ be an irreducible, aperiodic,
homogeneous Markov chain with finite state space
${\cal A}_m = \{\alpha_1 < \cdots < \alpha_m\}$,
transition  matrix $P$, and stationary distribution
$(\pi_1,\pi_2,\dots,\pi_m)$.
Let $\Sigma = (\sigma_{r,s})_{1\le r,s \le m}$ 
be the associated asymptotic covariance matrix,
as given in \eqref{item6la},
and let the law of $X_0$ be given by the
stationary distribution. 
Let $\tau$ be the permutation of $\{1,2,\dots,m\}$
such that $\pi_{\tau(i)} \ge \pi_{\tau(i+1)}$,
and $\tau(i) < \tau(j)$ whenever
$\pi_{\tau(i)} = \pi_{\tau(j)}$ and $i < j$.
For each $1 \le r \le m$,
let $V^r_n$ be the sum of the lengths 
of the first $r$ rows of the associated
Young diagrams,
and let $\nu_r = \sum_{j=1}^r \pi_{\tau(j)}$.  
Finally, let $d_r$ be
the multiplicity of $\pi_{\tau(r)}$, and let

\begin{equation*}
m_r = \begin{cases} 0,                                        &\text{if $\pi_{\tau(r)} = \pi_{\tau(1)}$,}\\
		    \max\{i: \pi_{\tau(i)} < \pi_{\tau(r)}\}, & \text{otherwise.}
      \end{cases}
\end{equation*}

\noindent Then, for each $1 \le r \le m$,

\begin{align}\label{item7r}
\frac{V^r_n -n\nu_r }{\sqrt{n}} \Longrightarrow  V^r_{\infty} &:= \sum_{i=1}^{m_r} \sigma_{\tau(i)} \tilde B^{\tau(i)}(1)  \nonumber\\
&\qquad  + \max_{I_{r-m_r,d_r}} 
   \sum_{j=1}^{r-m_r} \sum_{\ell=j}^{(d_r + m_r - r+j)}\!\!\sigma_{\tau(m_r + \ell)} \left(\tilde B^{\tau(m_r + \ell)}(t_{j,\ell})  
   - \tilde B^{\tau(m_r + \ell)}(t_{j,\ell-1})\right),
\end{align}

\noindent where the first sum on the right-hand side
of \eqref{item7r} is understood to be $0$, if $m_r = 0$.
Above, $\sigma_r^2 = \sigma_{r,r}$,
and $(\tilde B^1(t),\tilde B^2(t),\dots,\tilde B^m(t))$ is an $m$-dimensional Brownian
motion, with covariance matrix 
$\tilde{\Sigma} = ({\tilde\sigma}_{r,s})_{1\le r,s \le m}$ 
given by

\begin{equation}\label{item7s}
({\tilde\sigma}_{r,s})= t({\sigma}_{r,s})/\sigma_r \sigma_s,
\end{equation}

\noindent for $1\le r,s \le m$.
Moreover, for any $1 \le k \le m$,

\begin{equation}\label{item7saa}
\left(\frac{V^1_n -n\nu_1 }{\sqrt{n}},\frac{V^2_n -n\nu_2 }{\sqrt{n}},\dots,\frac{V^k_n -n\nu_k }{\sqrt{n}}\right) 
\Longrightarrow \left(V^1_{\infty},V^2_{\infty},\dots,V^k_{\infty} \right).
\end{equation}

\end{theorem}

\begin{Rem}
(i) The critical indices $d_r$ and $m_r$ in Theorem~\ref{thm4} are chosen so
that 
$$\pi_{\tau(m_r)} > \pi_{\tau(m_r + 1)} = \cdots = \pi_{\tau(r)} = \cdots =
	  \pi_{\tau(m_r + d_r)} > \pi_{\tau(m_r + d_r + 1)}.$$
	  
\noindent Thus, the functional in \eqref{item7r} consists of a sum of
$m_r$ Gaussian random variables and a maximal functional involving only $d_r$ of
the $m$ one-dimensional Brownian motions.  

(ii) Another, perhaps more natural, way of describing the covariance structure of the
$m$-dimensional Brownian motion in Theorem \ref{thm4} is to
note that $(\sigma_1 \tilde B^1(t),\sigma_2 \tilde B^2(t)$,
$\dots,\sigma_m \tilde B^m(t))$
has covariance matrix $t\Sigma$.
\end{Rem}

Let us now examine the case $r=1$ of the above theorem. Here, as
previously noted, $V^1_n = LI_n$.  
Since $m_1 = 0$, 
\eqref{item7r} becomes

\begin{align}\label{item7ra}
\frac{LI_n - n\pi_{max}}{\sqrt{n}} \Longrightarrow  
\max_{(t_{.,.}) \in I_{1,d_1}} 
    \sum_{\ell=1}^{d_1 } \sigma_{\tau(\ell)} 
    \left(\tilde B^{\tau(\ell)}(t_{1,\ell})  - \tilde B^{\tau(\ell)}(t_{1,\ell-1})\right),
\end{align}

\noindent where we have written 
$\pi_{max}$ for $\pi_{\tau(1)}.$ 
The functional in \eqref{item7ra} is similar to the one
obtained in the iid case in \cite{HL},  
the essential difference being, 
not in the form of the Brownian functional,
but rather in the covariance structure of the
Brownian motions.

To see precisely where this difference comes into play, 
note that if the transition matrix $P$ is cyclic, 
then the covariance matrix of the Brownian motion 
is also cyclic. Consider then the $3$-letter aperiodic, homogeneous, 
doubly-stochastic Markov case. Since the
Brownian covariance matrix is symmetric, and, moreover,
degenerate, an additional circularity constraint forces 
it to have the permutation-symmetric structure seen in the 
iid uniform case.  In particular, $LI_n$ will have as a limiting law, 
up to a scaling factor, the maximal eigenvalue of the traceless $3 \times 3$
GUE:

\begin{align}\label{item7rab}
\frac{LI_n - n/3}{\sqrt{n}} \Longrightarrow  
\sigma  \max_{(t_{.,.}) \in I_{1,3}} 
    \sum_{\ell=1}^{3}  \left( \tilde B^{\ell}(t_{1,\ell})  - \tilde B^{\ell}(t_{1,\ell-1})\right),
\end{align}

\noindent where $\sigma  = \sigma_{\ell},$
for all $1 \le \ell \le 3$, and with the
Brownian covariance matrix given by

\begin{equation*}
{\tilde \Sigma}= \begin{pmatrix}
     1      &-1/2   &-1/2\\
    -1/2   &1       &-1/2\\
    -1/2   &-1/2    &1
  \end{pmatrix}t,
\end{equation*}

\noindent and where we have used the fact
that $\tau(\ell) = \ell$, for all $1 \le \ell \le 3$.

However, when $m\ge 4$, the cyclicity constraint
does not force the Brownian covariance matrix
to have the permutation-symmetric structure, 
as the following example shows for $m=4$.

\begin{Example}\label{Ex1}
Consider the following doubly-stochastic, aperiodic,
cyclic transition matrix:

\begin{equation}\label{item7t}
P =   \begin{pmatrix}
    0.4   &0.3    &0.2   &0.1\\
    0.1   &0.4    &0.3   &0.2\\
    0.2   &0.1    &0.4   &0.3\\
    0.3   &0.2    &0.1   &0.4
  \end{pmatrix}.
\end{equation}

\noindent While the doubly-stochastic nature of $P$
ensures that the stationary distribution is
uniform, the covariance matrix of the limiting Brownian motion,
at three-decimal accuracy, is computed to be

\begin{equation}\label{item7ta}
{\tilde \Sigma} = \begin{pmatrix}
    1.000    &-0.357   &-0.287   &-0.357\\
    -0.357   &1.000    &-0.357   &-0.287\\
    -0.287   &-0.357   &1.000    &-0.357\\
    -0.357   &-0.287   &-0.357   &1.000
  \end{pmatrix}t,
\end{equation}

\noindent  and $\sigma_r^2 = \sigma^2 := 0.263$,
for each $1 \le r \le 4$.  Thus, 

\begin{align}\label{item7tb}
\frac{LI_n -n/4}{\sqrt{n}} \Longrightarrow  \sigma \max_{(t_{.,.}) \in I_{1,4}} 
      \sum_{\ell=1}^{4}  \left(\tilde B^{\ell}(t_{1,\ell})  - \tilde B^{\ell}(t_{1,\ell-1})\right),
\end{align}

\noindent for $1 \le r \le 4$.  However, while the form 
of the functional is the same as in
the iid uniform case (up to the constant),
the covariance structure of the Brownian motion in \eqref{item7ta}
differs from that of the uniform iid case, {\it i.e.}, from

\begin{equation}\label{item7tc}
\begin{pmatrix}
    1      &-1/3   &-1/3   &-1/3\\
    -1/3   &1      &-1/3   &-1/3\\
    -1/3   &-1/3   &1      &-1/3\\
    -1/3   &-1/3   &-1/3   &1
  \end{pmatrix}t,
\end{equation}

\noindent and so the limiting distribution in \eqref{item7tb}
is not that of the uniform iid case, see also \cite[Sec.~3]{BGH} (and the conjectured final  
result stated in the 
introduction there), as well as    
Theorem~\ref{thm7b} and Theorem~\ref{thm7c} below.

\end{Example}

We thus see that Kuperberg's conjecture
regarding the shape of the RSK Young diagrams 
for random sequences generated by
aperiodic, homogeneous, and cyclic matrices
\cite{Ku} is not true for general $m$-alphabets. 
By simply extending the first-row analysis above to the
second and third rows, we see that it is true for $m=3$.
However, it fails for $m \ge 4$.
In the next section we shall see that
for the cyclic case the structure of $\Sigma$ can
be described in a way which 
delineates precisely when the limiting law is the spectrum of the traceless GUE.

In the more general doubly stochastic case, 
we have the following corollary:

\begin{corollary}\label{cor1a}
Let the transition matrix $P$ of Theorem \ref{thm4} be doubly stochastic.
Then, for every $1 \le r \le m$, $m_r = 0, d_r = m$, and 

\begin{align}\label{item7sa}
\frac{V^r_n -rn/m}{\sqrt{n}} \Longrightarrow  \max_{(t_{.,.}) \in I_{r,m}} 
      \sum_{j=1}^{r} \sum_{\ell=j}^{m-r+j} \sigma_{\ell} \left(\tilde B^{\ell}(t_{j,\ell})  - \tilde B^{\ell}(t_{j,\ell-1})\right).
\end{align}

If, moreover, the matrix $P$ has all entries of $1/m$, then

\begin{align}\label{item7sb}
\frac{V^r_n -rn/m}{\sqrt{n}} \Longrightarrow    \frac{\sqrt{m-1}}{m}  \max_{(t_{.,.}) \in I_{r,m}} 
    \sum_{j=1}^{r} \sum_{\ell=j}^{m-r+j}  \left(\tilde B^{\ell}(t_{j,\ell})  - \tilde B^{\ell}(t_{j,\ell-1})\right),
\end{align}

\noindent and the covariance matrix in \eqref{item7s} 
has all its off-diagonal terms equal to $-1/(m-1)$.

\end{corollary}

\noindent \begin{Proof}
For each $1 \le r \le m$, $\pi_r = 1/m$, and so $\nu_r = r/m$, $m_r = 0$, 
and the multiplicity $d_r = m$.  
Moreover, the 
permutation $\tau$ is simply the identity permutation.
This proves \eqref{item7sa}.  If, moreover, all the transition
probabilities are $1/m$, then
the multinomial nature of the underlying combinatorial
quantities $a_k^r$ tells us that
$\sigma_r^2 = (1/m)(1- 1/m)$, for each $1 \le r \le m$, and that
$\rho_{r_1,r_2} = -1/(m-1)$, for each $r_1 \ne r_2$,
thus  proving \eqref{item7sb}.\CQFD
\end{Proof}

To see that the functional in \eqref{item7sa}
is generally different from the uniform iid case,
even for $m=3$, consider the following 
non-cyclic example:

\begin{Example}\label{Ex2}
Let a doubly-stochastic (but non-cyclic), 
aperiodic Markov chain have 
transition matrix

\begin{equation}\label{Ex2a}
P =   \begin{pmatrix}
    0.4   &0.6    &0.0\\
    0.6   &0.0    &0.4  \\
    0.0   &0.4    &0.6  \\
  \end{pmatrix}.
\end{equation}

\noindent Again,
the doubly-stochastic nature of $P$
ensures that the stationary distribution is
uniform, and in the present example, 
the asymptotic covariance matrix, 
at three-decimal accuracy, is computed to be

\begin{equation}\label{Ex2b}
\begin{pmatrix}  
    0.459    &0.049   &-0.506\\
    0.049    &0.086   &-0.136\\
   -0.506   &-0.136    &0.642
  \end{pmatrix}.
\end{equation}

\noindent Note that, even though we have a uniform
stationary distribution, the asymptotic
variances ({\it i.e.}, the diagonal terms of \eqref{Ex2b}) 
have dramatically different values.
Moreover, according to
Remark \ref{Remiidnonunif},
in the uniform iid case, the only possibility for
the Brownian covariance matrix is that
the off-diagonal terms have value $-1/2$.
However, the Brownian motion 
covariance matrix obtained from 
\eqref{Ex2b} is

\begin{equation}\label{Ex2c}
\begin{pmatrix}   
   1.000     &0.246   &-0.935\\
    0.246    &1.000   &-0.577\\
   -0.935   &-0.577    &1.000
  \end{pmatrix}t.
\end{equation}

\noindent Not only are the off-diagonal terms 
different from $-1/2$, but in some cases
are even positive. In
short, the functional in \eqref{item7sa}
has a distribution which differs from
{\it any} iid case (even non-uniform), see also \cite[Sec.~3]{BGH}.

\end{Example}

\begin{Rem}
Generalizing a result of Baryshnikov \cite{Ba} and
of Gravner, Tracy, and Widom \cite{GTW}
on the representation of the maximal eigenvalue
of an $m \times m$ element of the GUE, Doumerc \cite{Do} 
found a Brownian functional expression for all the
eigenvalues of an $m \times m$ element of the GUE.
Our expression in \eqref{item7sb}
is similar, with the exception that
our $m$-dimensional Brownian motion is constrained by a
zero-sum condition, and, moreover, has a 
different covariance structure.
(We note, moreover, that the parameters over which the 
Brownian functional of \cite{Do}  is maximized
might be intended to range over a slightly larger set
which corresponds to our $I_{r,m}$.)  
Using a path-transformation technique relating the
joint distribution of a certain transformation of
$n$ continuous processes to the joint distribution
of the processes conditioned never to leave the Weyl
chamber, O'Connell and Yor \cite{OY2} employed
queuing-theoretic arguments to obtain different types of
representations for the
entire spectrum of the $m \times m$ element of the GUE.
In a study of much more general transformations
of this type, Bougerol and Jeulin \cite{BJ} obtained this result as a special case 
and moreover these two representations are shown to be equivalent to each other 
in Biane, Bougerol and O'Connell \cite{BBO}.  
In both these works, the maximal and minimal eigenvalues
have the same representations and it is also the same as the one 
given by our Brownian functionals.  However, 
our representation for the rest of the spectrum is different from theirs and the  
equivalence with the works cited above is not that immediate.   
Indeed, in our case, a single supremum is taken over 
a rather involved set while in theirs, sequences of maxima and infima are taken over 
much simpler sets.   Another piece of work of interest is reported by Kirillov~\cite{Ki},  
where minima and maxima can be disentangled via classes of linear transformations (see Theorem~3.3 and 
Theorem~3.5 in \cite{Ki}).   Finally, \cite{BGH} presents equalities in law between the spectra of the principal 
minors of a GUE matrix 
and maximal functionals of independent Brownian motions obtained here. 
\end{Rem}

If $d_r=1$, {\it i.e.,} if the
$r^{th}$ most probable state is unique, then
the following result can be
viewed as lying at the other extreme 
from Corollary \ref{cor1a}:

\begin{corollary}\label{cor1b}
Let $1 \le r \le m$, and let $d_r = 1 $ in Theorem \ref{thm4}.
Then

\begin{align}\label{item7sc}
\frac{V^r_n -n\nu_r }{\sqrt{n}} \Longrightarrow  \sum_{i=1}^{r} \sigma_{\tau(i)} \tilde B^{\tau(i)}(1).
\end{align}

\end{corollary}

\noindent \begin{Proof}
If $d_r = 1$, then $m_r = r - 1$, and so the maximal term
of \eqref{item7r} contains only one summand, namely
$\sigma_{\tau(m_r + 1)} \tilde B^{\tau(m_r + 1)}(1) = \sigma_{\tau(r)} \tilde B^{\tau(r)}(1)$.
Including this term in the first summation term of \eqref{item7r}
proves \eqref{item7sc}. \CQFD
\end{Proof}

\begin{Rem}
The maximal term of the functional in \eqref{item7r}
is that of the doubly-stochastic,  $d_r$-letter case.
Indeed, the maximal term involves precisely $d_r$
Brownian motions over the $r-m_r$ rows.  Such a functional
would arise in a doubly-stochastic $d_r$-letter situation
with a covariance matrix consisting of the sub-matrix of
the original $\Sigma$ corresponding to the $d_r$ 
Brownian motions, as in Corollary \ref{cor1a}.
The Gaussian term corresponds to the functional 
of Corollary \ref{cor1b}. That is, in some sense,
the limiting law of \eqref{item7r} interpolates between
these two extreme cases.
\end{Rem}

\noindent \begin{Proof} {\bf (Theorem \ref{thm4})}
Since the $r=m$  case is trivial 
($V_n^m$ is then identically equal to $n$),
assume that $r < m$.
Recall the approximating functional \eqref{item7q}:

\begin{align}\label{item7u}
\frac{V^r_n - n\nu_r }{\sqrt{n}}
&=  \max_{I_{r,m}} 
    \biggl\{\sum_{j=1}^r\!\sum_{\ell=j}^{m-r+j}\!\!\!\sigma_{\ell} \left( \hat B^{\ell}_n(t_{j,\ell})  - 
    \hat B^{\ell}_n(t_{j,\ell-1})\right)
+ \sum_{j=1}^r\!\sum_{\ell=j}^{m-r+j}\!\!\sqrt{n}(\pi_{\ell} - \pi_{\tau(j)})   \left( t_{j,\ell}  - t_{j,\ell-1} \right)\biggr\}.
\end{align}

\noindent Introducing the notation $\Delta t_{j,\ell} := [t_{j,\ell-1}, t_{j,\ell}]$
and $M^{\ell}_n(\Delta t_{j,\ell}) := M^{\ell}_n(t_{j,\ell})  - M^{\ell}_n(t_{j,\ell-1})$,
for any $m$-dimensional process $M(t) = (M^1(t),M^2(t),\dots,M^m(t))$, $t \in [0,1]$,
we can rewrite \eqref{item7u} more compactly as

\begin{align}\label{item7ua}
\frac{V^r_n - n\nu_r }{\sqrt{n}}
=  \max_{I_{r,m}} 
    \biggl\{ \sum_{j=1}^r \sum_{\ell=j}^{m-r+j} \sigma_{\ell} \hat B^{\ell}_n(\Delta t_{j,\ell})
 -  \sqrt{n} \sum_{j=1}^r \sum_{\ell=j}^{m-r+j} (\pi_{\tau(j)} - \pi_{\ell})  |\Delta t_{j,\ell}| \biggr\}.  
\end{align}

The main idea of the proof to follow will be to show that the second summation of \eqref{item7ua} 
can, in effect, be eliminated by choosing the 
$\Delta t_{j,\ell}$ in an appropriate
manner.  We first claim that this sum is always non-negative.  Indeed, 

\begin{align*}
\sum_{j,\ell} (\pi_{\tau(j)} - \pi_{\ell})  |\Delta t_{j,\ell}| 
& = \sum_{j} \pi_{\tau(j)} \sum_{\ell} |\Delta t_{j,\ell}| - \sum_{j,\ell} \pi_{\ell} |\Delta t_{j,\ell}|\nonumber\\
& = \nu_r - \sum_{j,\ell} \pi_{\ell} |\Delta t_{j,\ell}|.\nonumber
\end{align*}

\noindent Now the sum $\sum_{j,\ell} \pi_{\ell} |\Delta t_{j,\ell}|$ is maximized by choosing $\Delta t_{j,\ell}$ in such a way that 
$\sum_{\ell} |\Delta t_{j,\ell}| = 1$, for $r$ indices $j$ having the $r$ largest probabilities $\pi_{j}$.  (One may do so by taking, {\it e.g.}, 
$\Delta t_{j,j} = [0,1]$ for $j = \tau(1), ... ,\tau(r)$).  This then gives us
$\nu_r - \sum_{j,\ell} \pi_{\ell} |\Delta t_{j,\ell}| \ge \nu_r - \nu_r = 0$ and the claim is proved. 

To eliminate the second summation of \eqref{item7ua}, we thus intuitively wish to set $|\Delta t_{j,\ell}| = 0$ for any $j$ for which $\pi_{j} < \pi_{\tau(r)}$.  To this end, we define a restricted set of parameters
$I_{r,m}(\varepsilon) = \{(t_{j,\ell}) \in I_{r,m}: \sum_{j,\ell} |\Delta t_{j,\ell}|{\bf 1}_{\{\pi_{\ell} < \pi_{\tau(r)}\}} \le \varepsilon r, \;
\sum_{j,\ell} |\Delta t_{j,\ell}|{\bf 1}_{\{\pi_{\ell} = \pi_{\tau(m_r)}\}} 
= m_r\}$ if $m_r > 0$ and 
$I_{r,m}(\varepsilon) = \{(t_{j,\ell}) \in I_{r,m}: \sum_{j,\ell} |\Delta t_{j,\ell}|{\bf 1}_{\{\pi_{\ell} < \pi_{\tau(r)}\}} \le \varepsilon r \}$ if $m_r = 0$.  The idea here is for $I_{r,m}(\varepsilon)$ to reflect parameter choices
that control intervals associated with small probabilities while at the same time maximizing the intervals associated with the $m_r$ largest probabilities.
With this parameter set in hand, we see that, provided $I_{r,m}(0) \ne \emptyset$, 

\begin{align}\label{this1}
\max_{I_{r,m}} \sum_{j=1}^{r} \sum_{\ell=j}^{m-r+j} \left(\sigma_{\ell} \hat{B}_n^{\ell}(\Delta t_{j,\ell}) 
- \sqrt{n} \left( \pi_{\tau(j)} - \pi_{\ell}\right) |\Delta t_{j,\ell}|\right) 
\ge \max_{I_{r,m}(0)} \sum_{j=1}^{r} \sum_{\ell=j}^{m-r+j} \sigma_{\ell} \hat{B}_n^{\ell}(\Delta t_{j,\ell}).
\end{align}

Moreover, by the
Invariance Principle and the Continuous Mapping Theorem,

\begin{align}\label{this2}
\max_{I_{r,m}(0)} \sum_{j=1}^{r} \sum_{\ell=j}^{m-r+j} \sigma_{\ell} \hat{B}_n^{\ell}(\Delta t_{j,\ell})
    \Longrightarrow \max_{I_{r,m}(0)} \sum_{j=1}^{r} \sum_{\ell=j}^{m-r+j} \sigma_{\ell} \tilde B^{\ell}(\Delta t_{j,\ell}).
\end{align}

We claim that, indeed, $I_{r,m}(0) \ne \emptyset$, 
and that, moreover,

\begin{align}\label{this3}
\max_{I_{r,m}} \sum_{j=1}^{r} \sum_{\ell=j}^{m-r+j} \left(\sigma_{\ell} \hat{B}_n^{\ell}(\Delta t_{j,\ell}) 
- \sqrt{n} \left( \pi_{\tau(j)} - \pi_{\ell}\right) |\Delta t_{j,\ell}| \right)
\Longrightarrow \max_{I_{r,m}(0)} \sum_{j=1}^{r} \sum_{\ell=j}^{m-r+j} \sigma_{\ell} \tilde B^{\ell}(\Delta t_{j,\ell}).
\end{align}

\noindent We will prove that $I_{r,m}(0) \ne \emptyset$
by creating a bijection between 
$I_{r,m}(0)$ and $I_{r-m_r,d_r}$.
To this end, for $1 \le i \le m_r$, let 
$\tilde{I}_{\tau(i),i} = [u_{\tau(i),i-1},u_{\tau(i),i}] = [0,1]$.
Next, choose any $(u_{.,.}) \in I_{r-m_r,d_r}$, 
and define further intervals $\tilde{I}_{\tau(m_r+j),\ell}$
for $1 \le j \le r-m_r$ and $1 \le \ell \le d_r$ as follows. 
If $u_{j-1,\ell} > 0$, let $\tilde{I}_{\tau(m_r+j),\ell}$
be the half-open interval
$\tilde{I}_{\tau(m_r+j),\ell} = \Delta u_{j,\ell} = (u_{j-1,\ell}, u_{j,\ell}]$,
and if $u_{j-1,\ell} = 0$, instead let $\tilde{I}_{\tau(m_r+j),\ell}$
be the closed interval 
$\tilde{I}_{\tau(m_r+j),\ell} = \Delta u_{j,\ell} = [0, u_{j,\ell}]$.

We now create a partition of these intervals in a
manner which relies on the ideas used in the 
proof of Theorem \ref{thm3}.
Consider the {\it set} of points 
$\{ u_{j,\ell}\}_{(1 \le j \le r-m_r, 1 \le \ell \le d_r)}$,
and order them as
$s_0 := 0 < s_1 < \cdots < s_{\kappa-1} < s_{\kappa}:= 1$,
for some integer $\kappa$,
and let 
$\Delta s_1 = [0,s_1]$, and
$\Delta s_q = (s_{q-1},s_q]$, 
for $2 \le q \le \kappa$.

Trivially, for each $1 \le q \le \kappa$, 
and for each $1 \le i \le m_r$, 
$\Delta s_q \subset \tilde{I}_{\tau(i),i}$.
Moreover, for each $1 \le j \le r-m_r$, there exists a unique
$\ell(j,q)$ such that 
$\Delta s_q \subset \tilde{I}_{\tau(m_r+j),\ell(j,q)}$.
For each $q$, consider the set of indices
$A_q := \{\tau(1), \dots, \tau(m_r)\} \cup \{ \tau(m_r + \ell(1,q)), \dots, \tau(m_r + \ell(r-m_r,q)) \}$,
and order these $r$ elements of $A_q$ as
$1 \le \tilde{\ell}(1,q) < \cdots < \tilde{\ell}(r,q) \le m$.

Using these partitions, we examine,
with foresight, the following functional
of a general $m$-dimensional process $(M(t))_{t\ge 0}$:

{\allowdisplaybreaks
\begin{align}
& \sum_{i=1}^{m_r} M^{\tau(i)}(1) 
  + \sum_{j=1}^{(r-m_r)} \sum_{\ell=j}^{(r-m_r+d_r-1)} M^{\tau(m_r + \ell)}(\Delta u_{j,\ell})\label{this4a}\\
&\qquad= \sum_{i=1}^{m_r} \left( \sum_{q=1}^{\kappa}  M^{\tau(i)}(\Delta s_q)  \right)
+ \sum_{j=1}^{(r-m_r)} \sum_{\ell=j}^{(r-m_r+d_r-1)} \left( \sum_{q:\Delta s_q \subset \tilde{I}_{\tau(m_r+j),\ell} } 
          M^{\tau(m_r + \ell)}(\Delta s_q)\right)\nonumber\\
&\qquad= \sum_{q=1}^{\kappa} \left( \sum_{i=1}^{m_r}  M^{\tau(i)}(\Delta s_q)  
             +  \sum_{j=1}^{(r-m_r)}  M^{\tau(m_r + \ell(j,q))}(\Delta s_q)   \right)\nonumber\\
&\qquad= \sum_{q=1}^{\kappa}   \sum_{j=1}^{r}        M^{\tilde{\ell}(j,q)}(\Delta s_q) 
 = \sum_{j=1}^{r}        \sum_{q=1}^{\kappa}   M^{\tilde{\ell}(j,q)}(\Delta s_q)\nonumber\\
&\qquad= \sum_{j=1}^{r}        \sum_{\ell=1}^{r}     M^{\tilde{\ell}(j,q)}(\Delta t_{j,\ell})\label{this4},
\end{align}
}

\noindent where, for each $1 \le j \le r$,
and for each $1 \le \ell \le m$,
$t_{j,\ell} := \max\{s_q: \ell \ge \tilde{\ell}(j,q) \}$.
(That is, for each $j$, we collapse together intervals
$\Delta s_q$ corresponding to the same component
$M^{\ell}$.)  Now, since our functional in 
\eqref{this4} has non-trivial summands only for
$\ell$ such that $\pi_{\tau(\ell)} \ge \pi_{\tau(r)}$, 
we have shown that
$(t_{.,.}) \in I_{r,m}(0)$.

The following example illustrates this argument.  Suppose
we have an alphabet of size
$m=10$, with 
$$(\pi_1, \pi_2, \dots, \pi_{10}) = (1/15,2/15,1/10,2/15,2/15,1/30,1/10,2/15,1/15,1/10).$$
Then,
$$\pi_{\tau(1)} = \pi_{\tau(2)} =\pi_{\tau(3)}  =\pi_{\tau(4)} = 2/15, \quad m_1 = m_2 = m_3 = m_4 = 0, \quad d_1=d_2=d_3=d_4=4,$$
$$\pi_{\tau(5)} = \pi_{\tau(6)} = \pi_{\tau(7)} = 1/10, \quad m_5 = m_6 =m_7 = 4, \quad d_5=d_6=d_7=3,$$
$$\pi_{\tau(8)} = \pi_{\tau(9)} = 1/15, \quad m_8=m_9 = 7, \quad d_8=d_9=2,$$
$$\pi_{\tau(10)} = 1/30, \quad m_{10} = 9, \quad d_{10}=1.$$

In particular, note that the two largest, distinct
probability values are $2/15$ and $1/10$,
of multiplicities $4$ and $3$, respectively.
We consider the case $r=6$ and show how 
$I_{r-m_r,d_r} = I_{6-4,3} = I_{2,3}$ corresponds
to an element of $I_{r,m}(0) = I_{6,10}(0)$.
Figure~\ref{figd} shows a typical element of 
the unconstrained index set $I_{6,10}$.

\begin{figure}
  \begin{center}
    \includegraphics[width=0.8\textwidth]{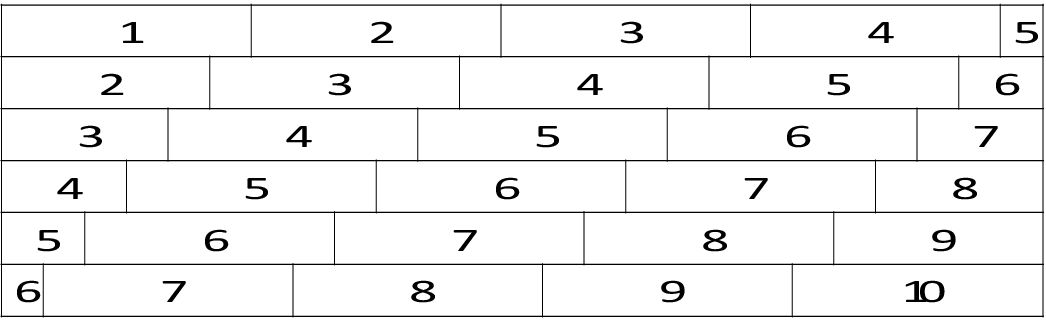}
  \end{center}
  \caption{A typical element of $I_{6,10}$.}
  \label{figd}
\end{figure}

Now $\tau(1) = 2, \tau(2)=4, \tau(3)=5, \tau(4) = 8, \tau(5) = 3,$
and $\tau(6) = 7$.  Our construction begins with
the amalgamation of the four rows corresponding
to the indices for which $\pi_i$ is strictly
greater than $\pi_{\tau(r)} = \pi_{\tau(6)} = 1/10$, with
$I_{2,3}$.  This amalgamation is shown in Figure~\ref{fige}.

\begin{figure}
  \begin{center}
    \includegraphics[width=0.7\textwidth]{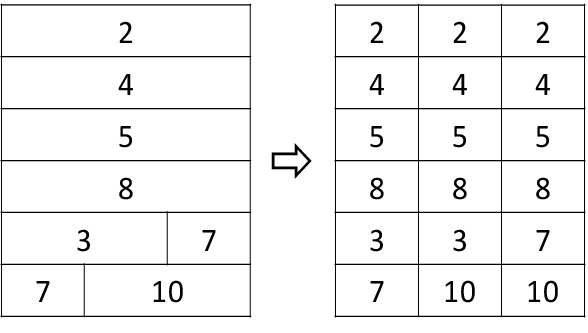}
  \end{center}
  \caption{Amalgamating top $4$ rows with $I_{2,3}$ (bottom $2$ rows).}
  \label{fige}
\end{figure}

Finally, we simply reorder each vertical column in the
original order of the indices,
as shown in Figure~\ref{figf}.
We see that, first of all,
we have constructed an element of $I_{6,10}$.  Moreover,
since we have four rows whose indices are associated
with the maximum value, and a remaining two rows whose indices
are associated with $\pi_{\tau(6)}$, we indeed have
an element of $I_{6,10}(0)$.  Note that the $6 \times 4 = 24$
free indices in $I_{6,10}$ (corresponding to the locations of
the $24$ vertical bars in Figure~\ref{figd}) have been
reduced to only two indices in $I_{6,10}(0)$.

\begin{figure}
  \begin{center}
    \includegraphics[width=0.7\textwidth]{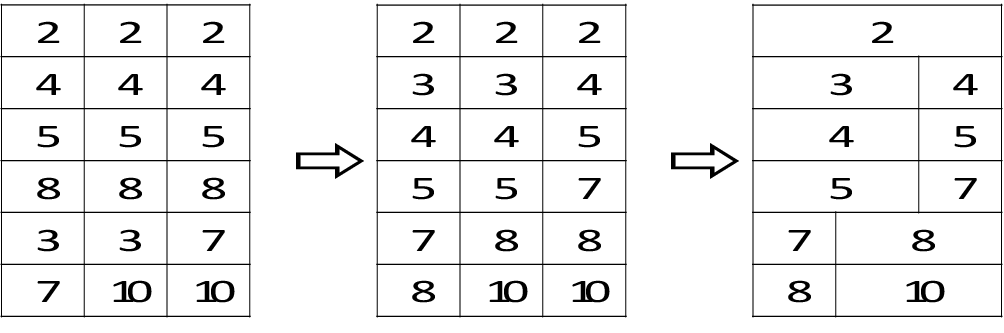}
  \end{center}
  \caption{Reordering vertically to obtain an element in $I_{6,10}(0).$}
  \label{figf}
\end{figure}

In addition, we may essentially reverse this construction,
starting with an element of $I_{r,m}(0)$ ($\ne \emptyset$), 
and so obtain an
element of $I_{r-m_r,d_r}$.    Indeed,
from the definitions of
$I_{r,m}(0)$ and $\nu_r$ we know that

$$\nu_r = \sum_{j=1}^{r} \pi_{\tau(j)} =  \sum_{j=1}^{r} \sum_{\ell=j}^{m-r+j} \pi_{\ell} |\Delta t_{j,\ell}|,$$

\noindent for any $(t_{.,.}) \in I_{r,m}(0)$.
However, we also have

{\allowdisplaybreaks
\begin{align*}
\sum_{j=1}^{r} \sum_{\ell=j}^{m-r+j} \pi_{\ell} |\Delta t_{j,\ell}| 
&= {\bf 1}_{\{m_r > 0\}} \Bigl( \sum_{j=1}^{r} \sum_{\ell=j}^{m-r+j} {\bf 1}_{\{\pi_{\tau(\ell)} \ge \pi_{\tau(m_r)}\}}  \pi_{\ell} |\Delta t_{j,\ell}|\nonumber\\
&  \qquad + \sum_{j=1}^{r} \sum_{\ell=j}^{m-r+j} {\bf 1}_{\{\pi_{\tau(\ell)} < \pi_{\tau(m_r)}\}}  \pi_{\ell} |\Delta t_{j,\ell}| \Bigr) \nonumber\\
& \qquad + {\bf 1}_{\{m_r = 0\}} \pi_{\tau(1)} \sum_{j=1}^{r} \sum_{\ell=j}^{m-r+j} |\Delta t_{j,\ell}|\nonumber\\
&\le {\bf 1}_{\{m_r > 0\}}( (\pi_{\tau(1)} + \cdots + \pi_{\tau(m_r)}) + (r-m_r)\pi_{\tau(r)}) 
+  {\bf 1}_{\{m_r = 0\}} r \pi_{\tau(1)} \nonumber\\
&= \nu_r,
\end{align*}
}

\noindent with equality holding throughout if and only if
$m_r = 0$ or $m_r > 0$ and 
$\sum_{j=1}^r |\Delta t_{j,\ell}| = 1,$
for all $\ell$ such that $\pi_{\tau(\ell)} \ge \pi_{\tau(m_r)}$,
and that, moreover,
$\sum_{j=1}^{r} \sum_{\ell=j}^{m-r+j} {\bf 1}_{\{\pi_{\tau(\ell)} = \pi_{\tau(r)}\}} |\Delta t_{j,\ell}| = r - m_r$.
If $m_r > 0$, then, for any $(t_{.,.}) \in I_{r,m}(0)$,
we may start with \eqref{this4},
and use again the permutation of the indices employed there.
We thus obtain the first term of \eqref{this4a},
which corresponds to the condition 
$\sum_{j=1}^r |\Delta t_{j,\ell}| = 1,$
for all $\ell$ such that $\pi_{\tau(\ell)} \ge \pi_{\tau(m_r)}$,
and also the second term of \eqref{this4a},
which corresponds to the other condition
$\sum_{j=1}^{r} \sum_{\ell=j}^{m-r+j} {\bf 1}_{\{\pi_{\tau(\ell)} = \pi_{\tau(r)}\}} |\Delta t_{j,\ell}| = r - m_r$.
If $m_r = 0$ the same reasoning holds, except that the first
term in \eqref{this4a} is taken to be zero.

Having thus established
a bijection between
$I_{r,m}(0)$ and $I_{r-m_r,d_r}$,
we may thus maximize over these two parameter sets, and so,
for any process $(M(t))_{t\ge 0}$, obtain the general result

\begin{align}\label{this5}
\sum_{i=1}^{m_r} M^{\tau(i)}(1) 
  + \max_{I_{r-m_r,d_r}} \sum_{j=1}^{(r-m_r)} \sum_{\ell=j}^{(r-m_r+d_r-1)} M^{\tau(m_r + \ell)}(\Delta u_{j,\ell}) 
  =  \max_{I_{r,m}(0)} \sum_{j=1}^{r}  \sum_{\ell=1}^{m-r+j}         M^{\tilde{\ell}(j,q)}(\Delta t_{j,\ell}).
\end{align}

We now proceed to show that \eqref{this3} holds.
First, fix $c > 0$, and, for each $1 \le \ell \le m$, set

\begin{equation}\label{this6}
c_{\ell} =   \begin{cases}
  c, \qquad \text{if $\pi_{\ell} < \pi_{\tau(r)}$,}\\
  0, \qquad \text{otherwise.}
  \end{cases}
\end{equation}

Next, let 
$\widehat{M}_n^{\ell}(t) = \sigma_{\ell} \hat{B}_n^{\ell}(t) - c_{\ell}t$,
and let 
$M^{\ell}(t) = \sigma_{\ell} \tilde B^{\ell}(t) - c_{\ell}t$.
Then, for $n$ large enough, namely, 
for $\sqrt{n} > c/(\pi_{\tau(r)} - \pi_{\tau(m_r + d_r +1)})$, 
we have that, almost surely,
for any $(t_{.,.}) \in I_{r,m}$,

\begin{align}\label{this7}
\sum_{j=1}^{r}  \sum_{\ell=1}^{m-r+j} \widehat{M}_n^{\ell}(\Delta t_{j,\ell}) 
\ge \sum_{j=1}^{r} \sum_{\ell=j}^{m-r+j} 
  \left(\sigma_{\ell}  \hat{B}_n^{\ell}(\Delta t_{j,\ell}) 
- \sqrt{n} \left( \pi_{\tau(j)} - \pi_{\ell}\right) |\Delta t_{j,\ell}|\right).
\end{align}

\noindent Hence, almost surely, both

\begin{align}\label{this8}
\max_{I_{r,m}} \sum_{j=1}^{r}  \sum_{\ell=1}^{m-r+j} \widehat{M}_n^{\ell}(\Delta t_{j,\ell}) 
\ge \max_{I_{r,m}}  \sum_{j=1}^{r} \sum_{\ell=j}^{m-r+j} \left( \sigma_{\ell} \hat{B}_n^{\ell}(\Delta t_{j,\ell}) 
- \sqrt{n} \left( \pi_{\tau(j)} - \pi_{\ell}\right) |\Delta t_{j,\ell}| \right),
\end{align}

\noindent and

\begin{align}\label{this9}
\max_{I_{r,m}(0)} \sum_{j=1}^{r}  \sum_{\ell=1}^{m-r+j} \widehat{M}_n^{\ell}(\Delta t_{j,\ell})
  = \max_{I_{r,m}(0)} \sum_{j=1}^{r} \sum_{\ell=j}^{m-r+j} \sigma_{\ell} \hat{B}_n^{\ell}(\Delta t_{j,\ell}).
\end{align}

Now choose any $z > 0$.  Then

\begin{align}\label{this10}
&\bbp \biggl(  \max_{I_{r,m}}  \sum_{j=1}^{r} \sum_{\ell=j}^{m-r+j} \left(\sigma_{\ell} \hat{B}_n^{\ell}(\Delta t_{j,\ell}) 
      - \sqrt{n} \left( \pi_{\tau(j)} - \pi_{\ell}\right) |\Delta t_{j,\ell}|\right) - \max_{I_{r,m}(0)} \sum_{j=1}^{r} \sum_{\ell=j}^{m-r+j} \sigma_{\ell} \hat{B}_n^{\ell}(\Delta s_{j,\ell}) > z \biggr) \nonumber\\
& \quad \quad \quad \le \bbp\left(\max_{I_{r,m}} \sum_{j=1}^{r}  \sum_{\ell=1}^{m-r+j} \widehat{M}_n^{\ell}(\Delta t_{j,\ell})
  - \max_{I_{r,m}(0)} \sum_{j=1}^{r}  \sum_{\ell=1}^{m-r+j} \widehat{M}_n^{\ell}(\Delta t_{j,\ell}) > z \right),
\end{align}

\noindent so that

{\allowdisplaybreaks
\begin{align}\label{this11}
&\limsup_{n\rightarrow \infty} 
  \bbp \biggl( \max_{I_{r,m}}  \sum_{j=1}^{r} \sum_{\ell=j}^{m-r+j} \left(\sigma_{\ell} \hat{B}_n^{\ell}(\Delta t_{j,\ell}) 
       - \sqrt{n} \left( \pi_{\tau(j)} - \pi_{\ell}\right) |\Delta t_{j,\ell}|\right)  - \max_{I_{r,m}(0)} 
       \sum_{j=1}^{r} \sum_{\ell=j}^{m-r+j} \sigma_{\ell} \hat{B}_n^{\ell}(\Delta s_{j,\ell}) > z \biggr)\nonumber\\
&\ \quad \quad \le \limsup_{n\rightarrow \infty} 
   \bbp\biggl( \max_{I_{r,m}} \sum_{j=1}^{r}  \sum_{\ell=1}^{m-r+j} \widehat{M}_n^{\ell}(\Delta t_{j,\ell})
    - \max_{I_{r,m}(0)} \sum_{j=1}^{r}  \sum_{\ell=1}^{m-r+j} \widehat{M}_n^{\ell}(\Delta t_{j,\ell}) > z  \biggr)\nonumber\\
& \quad \quad = \bbp\left(\max_{I_{r,m}} \sum_{j=1}^{r}  \sum_{\ell=1}^{m-r+j} M^{\ell}(\Delta t_{j,\ell})
     - \max_{I_{r,m}(0)} \sum_{j=1}^{r}  \sum_{\ell=1}^{m-r+j} M^{\ell}(\Delta t_{j,\ell}) > z  \right),
\end{align}
}

\noindent by the Invariance Principle and the Continuous Mapping Theorem.  
Now, for any $0 \le \varepsilon \le 1$, $I_{r,m}(0) \subset I_{r,m}(\varepsilon) \subset I_{r,m}(1) \subset I_{r,m}$.
We bound \eqref{this11} using this family of subsets as follows:

{\allowdisplaybreaks
\begin{align}\label{this12}
& \bbp\left(\max_{I_{r,m}} \sum_{j=1}^{r}  \sum_{\ell=1}^{m-r+j} M^{\ell}(\Delta t_{j,\ell})
    - \max_{I_{r,m}(0)} \sum_{j=1}^{r}  \sum_{\ell=1}^{m-r+j} M^{\ell}(\Delta t_{j,\ell}) > z  \right)\nonumber\\
&\le \bbp\left(\max_{I_{r,m}(\varepsilon)} \sum_{j=1}^{r}  \sum_{\ell=1}^{m-r+j} M^{\ell}(\Delta t_{j,\ell})
    - \max_{I_{r,m}(0)} \sum_{j=1}^{r}  \sum_{\ell=1}^{m-r+j} M^{\ell}(\Delta t_{j,\ell}) > z  \right)\nonumber\\
&\quad +  \bbp\left(\max_{I_{r,m} \backslash I_{r,m}(\varepsilon)} \sum_{j=1}^{r}  \sum_{\ell=1}^{m-r+j} M^{\ell}(\Delta t_{j,\ell})
    - \max_{I_{r,m}(0)} \sum_{j=1}^{r}  \sum_{\ell=1}^{m-r+j} M^{\ell}(\Delta t_{j,\ell}) > z  \right)\nonumber\\
&\le \bbp\left(\max_{I_{r,m}(\varepsilon)} \sum_{j=1}^{r}  \sum_{\ell=1}^{m-r+j} \tilde B^{\ell}(\Delta t_{j,\ell})
    - \max_{I_{r,m}(0)} \sum_{j=1}^{r}  \sum_{\ell=1}^{m-r+j} \tilde B^{\ell}(\Delta s_{j,\ell}) > z  \right)\nonumber\\
&\quad +  \bbp\left(\max_{I_{r,m} \setminus I_{r,m}(\varepsilon)} \sum_{j=1}^{r}  \sum_{\ell=1}^{m-r+j} \tilde B^{\ell}(\Delta t_{j,\ell})
    - \max_{I_{r,m}(0)} \sum_{j=1}^{r}  \sum_{\ell=1}^{m-r+j} \tilde B^{\ell}(\Delta s_{j,\ell}) > z + \varepsilon r c  \right)\nonumber\\
&\le \bbp\left(\max_{I_{r,m}(\varepsilon)} \sum_{j=1}^{r}  \sum_{\ell=1}^{m-r+j} \tilde B^{\ell}(\Delta t_{j,\ell})
    - \max_{I_{r,m}(0)} \sum_{j=1}^{r}  \sum_{\ell=1}^{m-r+j} \tilde B^{\ell}(\Delta s_{j,\ell}) > z  \right)\nonumber\\
&\quad +  \bbp \biggl( \max_{I_{r,m}} \sum_{j=1}^{r}  \sum_{\ell=1}^{m-r+j} \tilde B^{\ell}(\Delta t_{j,\ell})
- \max_{I_{r,m}(0)} \sum_{j=1}^{r}  \sum_{\ell=1}^{m-r+j} \tilde B^{\ell}(\Delta s_{j,\ell}) > z + \varepsilon r c  \biggr).
\end{align}
}

We can now take the limsup in \eqref{this12}, as
$c \rightarrow \infty$, and then, as
$\varepsilon \rightarrow 0$, and
so establish convergence to zero in probability,{\it i.e.,}

\begin{align}\label{this13}
\max_{I_{r,m}}  \sum_{j=1}^{r} \sum_{\ell=j}^{m-r+j} \left(\sigma_{\ell} \hat{B}_n^{\ell}(\Delta t_{j,\ell}) 
       - \sqrt{n} \left( \pi_{\tau(j)} - \pi_{\ell}\right) |\Delta t_{j,\ell}|  \right)  - \max_{I_{r,m}(0)} \sum_{j=1}^{r} 
       \sum_{\ell=j}^{m-r+j} \sigma_{\ell} \hat{B}_n^{\ell}(\Delta s_{j,\ell}) 
\stackrel{\bbp}{\longrightarrow} 0.
\end{align}

\noindent Since

\begin{equation}\label{this13b}
\max_{I_{r,m}(0)} \sum_{j=1}^{r} \sum_{\ell=j}^{m-r+j} \sigma_{\ell} \hat{B}_n^{\ell}(\Delta s_{j,\ell})
\Longrightarrow \max_{I_{r,m}(0)} \sum_{j=1}^{r} \sum_{\ell=j}^{m-r+j} \sigma_{\ell} \tilde B^{\ell}(\Delta s_{j,\ell}),
\end{equation}

\noindent by the Converging Together Lemma, we have 
proved \eqref{this3}.  Equation \eqref{item7r}
of the theorem follows from the bijection between
$I_{r,m}(0)$ and $I_{r-m_r,d_r}$ described in
the general result \eqref{this5}.

Finally, we can obtain the convergence of the joint distribution
in \eqref{item7saa} in the following
manner.  Given any $(\theta_1,\theta_2,\dots,\theta_r) \in \bbr^r$,
we have

{\allowdisplaybreaks
\begin{align}\label{this13c}
&\sum_{k=1}^{r} \theta_k \left( \frac{V^k_n -n\nu_k }{\sqrt{n}} \right)\nonumber\\
&\qquad=  \sum_{k=1}^{r} \theta_k \biggl( \max_{I_{k,m}}  \sum_{j=1}^{k} \sum_{\ell=j}^{m-k+j} \left(\sigma_{\ell} \hat{B}_n^{\ell}(\Delta t_{j,\ell}) 
       - \sqrt{n} \left( \pi_{\tau(j)} - \pi_{\ell}\right) |\Delta t_{j,\ell}|  \right)   \biggr) \nonumber\\
&\qquad=  \sum_{k=1}^{r} \theta_k \biggl( \max_{I_{k,m}}  \sum_{j=1}^{k} \sum_{\ell=j}^{m-k+j} \left(\sigma_{\ell} \hat{B}_n^{\ell}(\Delta t_{j,\ell}) 
       - \sqrt{n} \left( \pi_{\tau(j)} - \pi_{\ell}\right) |\Delta t_{j,\ell}|  \right)   \nonumber\\
&\qquad - \max_{I_{k,m}(0)} \sum_{j=1}^{k} \sum_{\ell=j}^{m-k+j} \sigma_{\ell} \hat{B}_n^{\ell}(\Delta s_{j,\ell}) \biggr)
+ \sum_{k=1}^{r} \theta_k \biggl( \max_{I_{k,m}(0)} \sum_{j=1}^{k} \sum_{\ell=j}^{m-k+j} \sigma_{\ell} \hat{B}_n^{\ell}(\Delta s_{j,\ell})   \biggr).
\end{align}
}

Now from \eqref{this13}, the first summation on the right-hand side of \eqref{this13c}
converges to zero in probability, as $n \rightarrow \infty$.  Moreover,
the second summation 
is a continuous functional of $(\hat{B}_n^{1},\hat{B}_n^{2},\dots,\hat{B}_n^{m})$,
and so, by the Invariance Principle and Continuous Mapping Theorem,
converges.  Then the Converging Together Lemma, along with
the bijection result \eqref{this5}, gives

\begin{align}\label{this13d}
\sum_{k=1}^{r} \theta_k \left( \frac{V^k_n -n\nu_k }{\sqrt{n}} \right) \Longrightarrow \sum_{k=1}^{r} 
\theta_k \biggl( \max_{I_{k,m}(0)} \sum_{j=1}^{k} \sum_{\ell=j}^{m-k+j} \sigma_{\ell} \tilde B^{\ell}(\Delta s_{j,\ell})  \biggr)
= \sum_{k=1}^{r} \theta_k V^k_{\infty}.
\end{align}

Since \eqref{this13d} holds for arbitrary $(\theta_1,\theta_2,\dots,\theta_r) \in \bbr^r$,
by the Cram\'er-Wold Theorem, we have the joint convergence result
\eqref{item7saa}.\CQFD
\end{Proof}

Since the shape of the Young diagrams is
more naturally expressed in terms of the
$R_n^k$, rather than of the $V_n^k$,
we may restate the results of the
previous theorem as follows:

\begin{theorem}\label{thm5}
Let $(X_n)_{n \ge 0}$ be an irreducible, aperiodic,
homogeneous Markov chain with finite state space
${\cal A}_m = \{\alpha_1 < \cdots < \alpha_m\}$,
and with stationary distribution
$(\pi_1,\pi_2,\dots,\pi_m)$.
Then, in the notations of Theorem~\ref{thm4},

\begin{align}\label{this14}
\left( \frac{R_n^1 -n\pi_{\tau(1)}}{\sqrt{n}},\frac{R_n^2 -n\pi_{\tau(2)}}{\sqrt{n}}, \dots, \frac{R_n^m -n\pi_{\tau(m)}}{\sqrt{n}} \right) 
    \Longrightarrow (R_{\infty}^1, R_{\infty}^2, \dots, R_{\infty}^m),
\end{align}

\noindent where

\begin{align}\label{this15a}
R_{\infty}^1
&=  \max_{I_{1,d_1}} 
    \sum_{\ell=1}^{d_1} \sigma_{\tau(\ell)} \left(\tilde B^{\tau(\ell)}(t_{1,\ell})  - \tilde B^{\tau(\ell)}(t_{1,\ell-1}) \right),
\end{align}

\noindent and, for each $2 \le k \le m$,

\begin{align}\label{this15b}
&R_{\infty}^k
=  \sum_{i=m_{k-1}+1}^{m_k} \sigma_{\tau(i)} \tilde B^{\tau(i)}(1)  \nonumber\\
&\quad  + \max_{I_{k-m_k,d_k}} 
   \sum_{j=1}^{k-m_k} 
     \sum_{\ell=j}^{(d_k + m_k- k+j)}
       \sigma_{\tau(m_k + \ell)} \tilde B^{\tau(m_k + \ell)}(\Delta t_{j,\ell})\nonumber\\
&\quad  - \max_{I_{k-1-m_{k-1},d_{k-1}}} 
   \sum_{j=1}^{k-1-m_{k-1}} 
     \sum_{\ell=j}^{(d_{k-1} + m_{k-1}- k+1+j)}
       \sigma_{\tau(m_{k-1} + \ell)} \tilde B^{\tau(m_{k-1} + \ell)}(\Delta t_{j,\ell}),
\end{align}

\noindent where we use the notation
$\tilde B^{s}(\Delta t_{j,\ell}) = \tilde B^{s}(t_{j,\ell}) - \tilde B^{s}(t_{j,\ell-1})$,
for any $1 \le s \le m$, $1 \le j \le k$, and $1 \le \ell \le m$,
and where the first sum on the right-hand side
of \eqref{this15b} is understood to be $0$, if $m_k = m_{k-1}$.

\end{theorem}

\noindent \begin{Proof}
First, $R_n^1 = V_n^1$, and,
for each $2 \le k \le m$, 
$R_n^k = V_n^k - V_n^{k-1}$.
Expressing these equalities
at the multivariate level,
we have

{\allowdisplaybreaks
\begin{align}\label{this15c}
&\left( \frac{R_n^1 -n\pi_{\tau(1)}}{\sqrt{n}},\frac{R_n^2 -n\pi_{\tau(2)}}{\sqrt{n}}, \dots, \frac{R_n^m -n\pi_{\tau(m)}}{\sqrt{n}} \right)  \nonumber\\
&\qquad\qquad= \left( \frac{V_n^1 -n\pi_{\tau(1)}}{\sqrt{n}},\frac{V_n^2 - V_n^1 -n\pi_{\tau(2)}}{\sqrt{n}}, \dots, \frac{V_n^m -V_n^{m-1} -n\pi_{\tau(m)}}{\sqrt{n}} \right) \nonumber\\
&\qquad\qquad= \left( \frac{V_n^1 -n\nu_1 }{\sqrt{n}},\frac{V_n^2 -n\nu_2 }{\sqrt{n}}, \dots, \frac{V_n^m - n\nu_m}{\sqrt{n}} \right) - \left(0,\frac{V_n^1 -n\nu_1 }{\sqrt{n}}, \dots, \frac{V_n^{m-1} - n\nu_{m-1} }{\sqrt{n}} \right) \nonumber\\
&\qquad\qquad\Longrightarrow (V_{\infty}^1, V_{\infty}^2,\dots, V_{\infty}^m) - (0, V_{\infty}^1, \dots, V_{\infty}^{m-1}) := (R_{\infty}^1, R_{\infty}^2, \dots, R_{\infty}^m),
\end{align}
}

\noindent where the weak convergence
follows immediately 
from the Continuous Mapping Theorem,
since the transformation is linear.  Equations \eqref{this15a} and \eqref{this15b}
follow simply from the Brownian expressions
for $(V_{\infty}^1, V_{\infty}^2,\dots, V_{\infty}^m)$
in Theorem \ref{thm4}.\CQFD
\end{Proof}

If all $m$ letters have unique stationary
probabilities, 
then we have the following corollary to
Theorem~\ref{thm5}:

\begin{corollary}\label{cor1c}
If the stationary distribution of Theorem~\ref{thm5}
is such that each $\pi_r$ is unique, then

\begin{align}\label{this15d}
\left( \frac{R_n^1 -n\pi_{\tau(1)}}{\sqrt{n}},\frac{R_n^2 -n\pi_{\tau(2)}}{\sqrt{n}}, \dots, \frac{R_n^m -n\pi_{\tau(m)}}{\sqrt{n}} \right) 
    \Longrightarrow N((0,0,\dots,0),\Sigma).
\end{align}

In other words, the limiting distribution is identical in law to
the spectrum of the diagonal matrix
$D = \text{diag}\{Z_1,Z_2,\dots,Z_m\}$, where
$(Z_1,Z_2,\dots,Z_m)$ is a centered normal
random vector with covariance matrix $\Sigma$.
\end{corollary}

\noindent \begin{Proof}
Now, for all $1 \le k \le m$, 
$d_k = 1 $, and
$m_k = k-1$, so that

\begin{align*}
R_{\infty}^1
&=  \max_{I_{1,d_1}} 
    \sum_{\ell=1}^{d_1} \sigma_{\tau(\ell)} \left(\tilde B^{\tau(\ell)}(t_{1,\ell})  - \tilde B^{\tau(\ell)}(t_{1,\ell-1}) \right) = \sigma_{\tau(1)}\tilde B^{\tau(1)}(1),
\end{align*}

\noindent and, for each $2 \le k \le m$,

{\allowdisplaybreaks
\begin{align*}
&R_{\infty}^k
=  \sum_{i=m_{k-1}+1}^{m_k} \sigma_{\tau(i)} \tilde B^{\tau(i)}(1)   + \max_{I_{k-m_k,d_k}} 
   \sum_{j=1}^{k-m_k} 
     \sum_{\ell=j}^{(d_k + m_k- k+j)}
       \sigma_{\tau(m_k + \ell)} \tilde B^{\tau(m_k + \ell)}(\Delta t_{j,\ell})\nonumber\\
&\qquad \qquad - \max_{I_{k-1-m_{k-1},d_{k-1}}} 
   \sum_{j=1}^{k-1-m_{k-1}} 
     \sum_{\ell=j}^{(d_{k-1} + m_{k-1}- k+1+j)}
       \sigma_{\tau(m_{k-1} + \ell)} \tilde B^{\tau(m_{k-1} + \ell)}(\Delta t_{j,\ell})\nonumber\\
&\quad=  \sum_{i=k-1}^{k-1} \sigma_{\tau(i)} \tilde B^{\tau(i)}(1)  + \max_{I_{1,1}} 
   \sum_{j=1}^{1} 
     \sum_{\ell=j}^{j}
       \sigma_{\tau(k - 1 + \ell)} \tilde B^{\tau(k - 1 + \ell)}(\Delta t_{j,\ell})\nonumber\\
&\qquad \qquad - \max_{I_{1,1}} 
   \sum_{j=1}^{1} 
     \sum_{\ell=j}^{j}
       \sigma_{\tau(k-2 + \ell)} \tilde B^{\tau(k-2 + \ell)}(\Delta t_{j,\ell})\nonumber\\
&\quad=   \sigma_{\tau(k-1)} \tilde B^{\tau(k-1)}(1) 
         +\sigma_{\tau(k)} \tilde B^{\tau(k)}(1)
         -\sigma_{\tau(k-1)} \tilde B^{\tau(k-1)}(1)\nonumber\\
&\quad=   \sigma_{\tau(k)} \tilde B^{\tau(k)}(1).
\end{align*}
}

\noindent Moreover, the joint law result
 for $(R_{\infty}^1,R_{\infty}^2,\dots,R_{\infty}^m)$
holds as well, and this is clearly a multivariate normal distribution, with mean
$(0,0,\dots,0)$ and covariance matrix $\Sigma$.  
Since the spectrum of a diagonal matrix consists of its diagonal elements,
the final claim of the corollary holds.\CQFD
\end{Proof}

\begin{Rem}\label{RMTconnect}
The joint law of 
$(R_{\infty}^1, R_{\infty}^2, \dots, R_{\infty}^m)$
in the iid uniform alphabet case is identical to the
joint law of the eigenvalues of an
$m \times m$ traceless GUE matrix.
Corollary \ref{cor1c} also gives a spectral characterization
for the unique probability case, in particular,
for a non-uniform iid alphabet with unique stationary
probabilities.  
This is consistent with
the characterization of the limiting law of
$LI_n$ in the non-uniform iid case, due to 
Its, Tracy, and Widom \cite{ITW1,ITW2}, as that of the
largest eigenvalue of the block associated with
the most probable letters among a
direct sum of independent GUE matrices
whose dimensions correspond
to the multiplicities $d_r$ of
Theorem \ref{thm4} and \ref{thm5}, 
subject to the condition
that $\sum_{r=1}^m \sqrt{\pi_{\tau(r)}}X_r = 0$,
where $X_1,X_2,\dots,X_m$ are the diagonal
elements of the random matrix.
\end{Rem}

\begin{Rem}\label{GenTrace}
The difference between the zero-trace condition
$\sum_{r=1}^m X_r = 0$ and
the generalized traceless condition
$\sum_{r=1}^m \sqrt{\pi_{\tau(r)}} X_r = 0$
amounts to nothing more than a difference
in the choice of scaling for each row
$R^{r}_n$.
We will find it more natural to express our
results using the normalization associated
with the zero-trace condition
$\sum_{r=1}^m X_r = 0$
\end{Rem}

\section{Fine Structure of the Brownian Functional}

So far, we have seen that the limiting shape of the RSK random Young
diagrams generated by an aperiodic, irreducible, homogeneous
Markov chain can be expressed as a Brownian functional.
The form of this functional
is similar to the iid case;
the essential difference is in the covariance
structure of the Brownian motion.
We begin our study of the consequences of this difference.

In the iid uniform $m$-alphabet case, Johansson \cite{Jo} proved
that the limiting shape of the Young diagrams had a
joint law which is that of the spectrum of an $m \times m$
traceless GUE matrix. An immediate consequence of this
result is that the limiting shape of the Young diagrams 
contains simple symmetries, {\it e.g.,}  for each $1 \le r \le m$,
$R^r_{\infty} \stackrel{\cal{L}}{=} - R^{m+1-r}_{\infty}.$
Now, as was seen in Corollary \ref{cor1a} 
of Theorem \ref{thm4}, the form of the Brownian  functional 
in the doubly stochastic case involved only the maximal term.
We will see that that there
is also a pleasing symmetry to the limiting shape of 
Young diagrams in the doubly stochastic case 
by examining a natural bijection
between the parameter set $I_{r,m}$ and
$I_{m-r,m}$, for any $1 \le r \le m-1.$
Indeed, this result will follow as a corollary
to the following, more general, theorem:

\begin{theorem}\label{thm6}
The limiting functionals of Theorem~\ref{thm4} enjoy
the following symmetry property:
for every $1 \le r \le m-1$,

{\allowdisplaybreaks
\begin{align}\label{this16}
&V^r_{\infty} := \sum_{i=1}^{m_r} \sigma_{\tau(i)} \tilde B^{\tau(i)}(1)  + \max_{t(\cdot,\cdot) \in I_{r-m_r,d_r}} 
   \sum_{j=1}^{r-m_r} 
   \sum_{\ell=j}^{(m_r+d_r-r+j)} \sigma_{\tau(m_r + \ell)} \tilde B^{\tau(m_r + \ell)}(\Delta t_{j,\ell})\nonumber\\
&\qquad \stackrel{\cal{L}}{=}
\sum_{i=m_r+d_r+1}^{m} \sigma_{\tau(i)} \tilde B^{\tau(i)}(1) + \max_{u(\cdot,\cdot) \in I_{m_r+d_r-r,d_r}} 
     \sum_{j=1}^{m_r+d_r-r} 
     \quad \sum_{\ell=j}^{r-m_r+j} \sigma_{\tau(m_r + \ell)} \tilde B^{\tau(m_r + \ell)}(\Delta u_{j,\ell}),
\end{align}
}

\noindent where $\tilde B^{\ell}(\Delta) :=
\tilde B^{\ell}(t) - \tilde B^{\ell}(s)$, for $\Delta = [s,t]$,
and where the non-maximal terms on the left and right-hand
sides of \eqref{this16} are identically zero if $m_r = 0$,
or $m_r + d_r = m$, respectively. 

\end{theorem}

\begin{Rem}\label{Remsym}
Recall that, from the definitions of $m_r$ and $d_r$,
the non-maximal summation terms on the left and 
right-hand sides of \eqref{this16}
reflect the letters which have, respectively, greater and smaller
stationary probabilities than $\pi_{\tau(r)}$.
Recall, moreover, that the maximal terms are associated with the indices
having the same stationary probability as $\pi_{\tau(r)}$.
The maximal term on the left-hand side of \eqref{this16}
involves a summation over $r-m_r$ rows,
while the one on the right-hand side involves
$m_{r+1}-r$ rows.  Thus, in a sense, the two maximal
terms in \eqref{this16} split $d_r = m_{r+1}-m_r$ rows between
themselves.  In summary, the functional on the
right-hand side of \eqref{this16} corresponds to
the sum of the lengths of the $m-r$ {\it bottom} rows of the 
Young diagrams.
\end{Rem}

\noindent \begin{Proof} 
Without loss of generality, we may assume that
$\tau(j) = j$, for all $1 \le j \le m$.
Fix $1 \le r \le m-1$, and
for any point $\tilde t$ in the index set $I_{r-m_r,d_r}$,
define the interval $\Delta t_{j+m_r,\ell}$, 
for $1 \le j \le r-m_r$ and $1 \le \ell \le d_r$
as follows.
If $\tilde t_{j,\ell-1} = 0$, let
$\Delta t_{j+m_r,\ell} = [0,\tilde t_{j,\ell}]$, and if
$\tilde t_{j,\ell-1} > 0$, let
$\Delta t_{j+m_r,\ell} = (\tilde t_{j,\ell-1}, \tilde t_{j,\ell}]$.
Furthermore, for each 
$1 \le j \le m_r$ or $m_{r+1} < j \le m$, 
set $\Delta t_{j,\ell} = [0,1]$, for $j = \ell$,
$\Delta t_{j,\ell} = \{0\}$, for $0 \le \ell < j$, and
$\Delta t_{j,\ell} = \{1\}$, for $j < \ell \le m.$
Next, as in the proof of Theorem \ref{thm4},
consider the {\it set} of points 
$\{\tilde t_{j,\ell}\}_{(1 \le j \le r-m_r, 1 \le \ell \le d_r)}$,
and order them as
$s_0 := 0 < s_1 < \cdots < s_{\kappa-1} < s_{\kappa}:= 1$,
for some integer $\kappa$,
and let 
$\Delta s_1 = [0,s_1]$, and
$\Delta s_q = (s_{q-1},s_q]$, 
for $2 \le q \le \kappa$.

Now, for each $1 \le q \le \kappa$, 
let $A_q$ consist of the indices $\ell$ for which
$\Delta s_q \cap \Delta t_{j,\ell} \ne \emptyset$.
Then, almost surely, 

{\allowdisplaybreaks
\begin{align}\label{this18}
\sum_{i=1}^{m_r} \sigma_{i} \tilde B^{i}(1)  
+ \sum_{j=1}^{r-m_r}    \sum_{\ell=j}^{(m_r+d_r-r+j)} \sigma_{m_r + \ell} \tilde B^{m_r + \ell}(\Delta t_{j,\ell}) 
& = \sum_{j=1}^{r} \sum_{\ell=1}^{m} \sigma_{\ell} \tilde B^{\ell}(\Delta t_{j,\ell})\nonumber\\
& = \sum_{j=1}^{r} \sum_{q=1}^{\kappa}  \sum_{\ell=1}^{m} \sigma_{\ell} \tilde B^{\ell}(\Delta t_{j,\ell} \cap \Delta s_q)\nonumber\\
& = \sum_{j=1}^{r} \sum_{q=1}^{\kappa}  \sum_{\ell \in A_q} \sigma_{\ell} \tilde B^{\ell}(\Delta s_q).
\end{align}
}

Now by the ``stairstep'' properties of $I_{r,m}$
there are precisely
$r$ elements in each $A_q$.  
Letting  $\tilde{A}_q = \{1,\dots,m\} \setminus A_q$,
for each $1 \le q \le \kappa$, we thus see that each 
$\tilde{A}_q$ contains exactly $m-r$ elements.
Let  $\tilde{\ell}_{j,q}$ be the $j^{th}$ smallest
element of $\tilde{A}_q$.  We claim that for each $1 \le j \le m-r$,
the sequence 
$\tilde{\ell}_{j,1}, \tilde{\ell}_{j,2}, \dots, \tilde{\ell}_{j,\kappa}$.
is weakly decreasing. 

Indeed, fix $1 \le j \le m-r$ and $1 \le q \le \kappa-1$,
and suppose that $\tilde{\ell}_{j,q}$ is less than all
the elements of $A_q$.  Then, by the properties of $I_{r,m}$,
the least element of $A_{q+1}$ is no smaller, so that
the $j^{th}$ smallest element of $\tilde{A}_q$,
$\tilde{\ell}_{j,q+1}$ is also $\tilde{\ell}_{j,q}$.
Next, suppose that $\tilde{\ell}_{j,q}$ is
greater than $k \ge 1$ elements of $A_q$.
Thus, $\tilde{\ell}_{j,q} = j + k$.  Then 
there are at most $k$ elements of $A_{q+1}$ which
are less than or equal to 
$\tilde{\ell}_{j,q}$, by the properties
of $I_{r,m}$.  But this implies that there
are at least $j$ elements of $\tilde{A}_{q+1}$
which are less than or equal to 
$\tilde{\ell}_{j,q}$.  Thus,
$\tilde{\ell}_{j,q+1} \le \tilde{\ell}_{j,q}$,
and the claim is proved.

Moreover, since each $A_q$ contains 
$\{1,2,\dots,m_r\}$, we see that necessarily
each $\tilde{A}_q$ contains 
$\{m_r+d_r+1, m_r+d_r+2,\dots, m\}.$

For each $1 \le j \le m-r$,
we may now amalgamate the intervals $\Delta s_q$
to obtain a partition of the unit interval.
Specifically, for each $1 \le j \le m-r$,
and each $1 \le \ell \le m$, let $\tilde{u}_{j,\ell}$ be the
smallest $s_q$ such that
$\tilde{\ell}_{j,q+1} \le \ell$. 
(We define $\tilde{u}_{j,0} = 1$, for all $1 \le j \le m-r$.)

Finally, and most crucially, recall that
$\sum_{\ell=1}^m \sigma_{\ell} \tilde B^{\ell}(t) = 0$, for all $t$.  
Then since 
$(\tilde B^1,\tilde B^2,\dots,\tilde B^m)  
\stackrel{{\cal L}}{=} (-\tilde B^1,-\tilde B^2,\dots,-\tilde B^m)$,

{\allowdisplaybreaks
\begin{align}\label{this19}
\sum_{j=1}^{r}   \sum_{q=1}^{\kappa}  \sum_{\ell \in A_q}         \sigma_{\ell} \tilde B^{\ell}(\Delta s_q) 
&= \sum_{j=1}^{m-r} \sum_{q=1}^{\kappa}  \sum_{\ell \in \tilde{A}_q} \left(-\sigma_{\ell} \tilde B^{\ell}(\Delta s_q)\right)\nonumber\\
& = -\sum_{i=m_r+d_r+1}^{m} \sigma_{i} \tilde B^{i}(1) 
           - \sum_{j=1}^{m_r+d_r-r} 
               \quad \sum_{\ell=1}^{m} \sigma_{m_r+\ell} \tilde B^{m_r+\ell}(\Delta u_{j,\ell})\nonumber\\
& \stackrel{\cal{L}}{=} 
\sum_{i=m_r+d_r+1}^{m} \sigma_{i} \tilde B^{i}(1) 
+ \sum_{j=1}^{m_r+d_r-r} 
     \quad \sum_{\ell=1}^{m}  \sigma_{m_r+\ell} \tilde B^{m_r+\ell}(\Delta u_{j,\ell}),
\end{align}
}

\noindent where 
$\Delta u_{j,\ell} = [u_{j,\ell-1},u_{j,\ell}]$. 
But, by the way we ordered each $A_q$, we must have 
$\Delta u_{j_1,\ell} \cap \Delta u_{j_2,\ell} = \emptyset$, for any $j_1 \ne j_2$.  Thus,
$u \in I_{m_r+d_r-r,d_r}$, and so we may restrict the summation over 
$\ell$ in \eqref{this19} to $\ell=j,\dots, r-m_r+j$, since the remaining terms are zero.
Equation \eqref{this16} follows immediately
by taking the maxima over $I_{r-m_r,d_r}$ 
and $I_{m_r+d_r-r,d_r}$ over the left-hand
and right-hand sides, respectively, of \eqref{this19}.\CQFD
\end{Proof}

For doubly stochastic transition matrices, the
symmetry is even more apparent:

\begin{corollary}\label{cor2}
Let the transition matrix $P$ of Theorem \ref{thm4} be doubly stochastic.
Then, for every $1 \le r \le m-1$,

\begin{align}\label{this19a}
&V_{\infty}^r : = \max_{t(\cdot,\cdot) \in I_{r,m}} \sum_{j=1}^{r} \sum_{\ell=j}^{m-r+j} 
\sigma_{\ell} \left(\tilde B^{\ell}(t_{j,\ell})  - \tilde B^{\ell}(t_{j,\ell-1})\right)\nonumber\\
&\qquad \stackrel{\cal{L}}{=}
\max_{u(\cdot,\cdot) \in I_{m-r,m}} \sum_{j=1}^{m-r} \sum_{\ell=j}^{r+j} 
\sigma_{\ell} \left(\tilde B^{\ell}(u_{j,\ell})  - \tilde B^{\ell}(u_{j,\ell-1})\right) := V_{\infty}^{m-r},
\end{align}

\noindent and so

\begin{align}\label{this19b}
\lim_{n \rightarrow \infty}   \frac{\sum_{j=1}^r      R_n^j - rn/m}{\sqrt{n}} 
\stackrel{\cal{L}}{=} \lim_{n \rightarrow \infty} \frac{rn/m - \sum_{j=m-r+1}^m  R_n^j}{\sqrt{n}}.
\end{align}

\noindent Moreover,

\begin{equation}\label{this19c}
(V_{\infty}^1,\dots,V_{\infty}^r) \stackrel{\cal{L}}{=}
(V_{\infty}^{m-1},\dots,V_{\infty}^{m-r}) .
\end{equation}

\end{corollary}

\noindent \begin{Proof} 
Since $m_r = 0$ and $d_r = m$ for all $1 \le r \le m$, 
the non-maximal terms on both sides of \eqref{this16}
disappear, and we have \eqref{this19a}.

To prove \eqref{this19b}, recall that $V_n^m = \sum_{j=1}^m R_n^j = n$.
Then, from the result just proved,

{\allowdisplaybreaks
\begin{align}\label{this20}
\frac{V_n^{m-r} - (m-r)n/m}{\sqrt{n}} 
&= \frac{\sum_{j=1}^{m-r} R_n^j- (m-r)n/m}{\sqrt{n}} \nonumber\\
&= \frac{\left(n - \sum_{j=m-r+1}^m R_n^j\right)- (m-r)n/m}{\sqrt{n}} \nonumber\\
&= \frac{rn/m - \sum_{j=m-r+1}^m R_n^j}{\sqrt{n}}\nonumber\\
&\Longrightarrow V_{\infty}^{m-r} \stackrel{\cal{L}}{=} V_{\infty}^{r},
\end{align}
}

\noindent and we have established the claimed symmetry.

Finally, the extension of \eqref{this19a} to \eqref{this19c} 
follows from a standard Cram\'er-Wold argument.\CQFD
\end{Proof}

\begin{Rem}
Since $R^{m}_{\infty} = -V^{m-1}_{\infty}$,
almost surely, Corollary \ref{cor2} states that
$R^{m}_{\infty} \stackrel{\cal{L}}{=} -R^{1}_{\infty}$.
From the symmetry of the Brownian motion, we thus
see that $R^{m}_{\infty}$ may be represented as
a {\it minimal} Brownian functional:

\begin{equation*}
R_{\infty}^m
=  \min_{I_{1,m}} 
    \sum_{\ell=1}^{m} \sigma_{\ell} \left(\tilde B^{\ell}(t_{1,\ell})  
                - \tilde B^{\ell}(t_{1,\ell-1}) \right),
\end{equation*}

\noindent a result also noted in \cite{OY2}.

\end{Rem}

Turning again to the cyclic case, recall that,
for $m \ge 4$, the limiting shape of the Young diagrams
in general differs from that of the iid uniform case.  The
following proposition characterizes the asymptotic covariance
matrices of such Markov chains.

\begin{proposition}\label{prop2}
Let $P$ be the $m \times m$ transition matrix
of an aperiodic, irreducible, cyclic Markov chain
on an $m$-letter, totally ordered alphabet, 
${\cal A}_m = \{\alpha_1 < \alpha_2 < \cdots < \alpha_m\}$, 
with 

\begin{equation}\label{this21}
P =
  \begin{pmatrix} 
  a_1 & a_m &\cdots & a_3 & a_2\\
  a_2 & a_1 &\ddots &&     a_3\\
  \vdots & \ddots & \ddots & \ddots &\vdots\\
  a_{m-1} && \ddots  & a_1  &a_m\\
  a_m &a_{m-1} &\cdots &a_2 &a_1
  \end{pmatrix}.
\end{equation}

\noindent For $1 \le j \le m$, let 
$\lambda_j = \sum_{k=1}^m a_k \omega^{(k-1)(j-1)}$ be an eigenvalue of $P$, 
where $\omega = exp (2\pi i/m)$ is the $m^{th}$ 
principal root of unity.  Let also  
$\gamma_j =  \lambda_j/(1-\lambda_j)$,  
for $2 \le j \le m$, and
$\beta_j = \cos (2\pi j/m)$,
for $0 \le j \le m$.  Then, the asymptotic covariance matrix $\Sigma$ 
is given by:\\

\noindent For $m = 2m_0 + 1$,

\begin{equation}\label{this22}
\Sigma = \frac{m-1}{m^2}M^{(1)} + \frac{4}{m^2} \sum_{j=2}^{m_0+1} Re(\gamma_j) M^{(j)},
\end{equation}

\noindent and for $m = 2m_0$,

\begin{equation}\label{this23}
\Sigma = \frac{m-1}{m^2}M^{(1)} + \frac{4}{m^2} \sum_{j=2}^{m_0} Re(\gamma_j) M^{(j)} + \frac{2}{m^2} \gamma_{m_0+1} M^{(m_0+1)},
\end{equation}

\noindent where $M^{(j)}$ 
is an $m \times m$  Toeplitz matrix with entries
$(M^{(j)})_{k,\ell} = \beta_{(j-1)|k-\ell|}$, for
$2 \le j \le m$, 
and $(M^{(1)})_{k,\ell} = \delta_{k,\ell} - (1-\delta_{k,\ell})/(m-1)$,
for $j=1$.
\end{proposition}

\noindent \begin{Proof} It is straightforward, and classical,
to verify that, for each $1 \le j \le m$,
$(1,\omega^{j-1}, \omega^{2(j-1)},\dots, \omega^{(m-1)(j-1)})$
is a left eigenvector of $P$, with eigenvalue 
$\lambda_j = \sum_{k=1}^m a_k \omega^{(k-1)(j-1)}$.  We can thus
write our standard Jordan decomposition of $P$, 
which in this case is a true diagonalization,
as $P = S^{-1} \Lambda S$,
where $\Lambda = \text{diag}(1,\lambda_2,\dots,\lambda_m)$,

\begin{equation}\label{this24}
S = 
  \begin{pmatrix} 
  1 & 1 &\cdots & 1 & 1\\
  1 & \omega   & \omega^2 &\cdots & \omega^{m-1}\\
  1 & \omega^2 & \omega^4 &\cdots & \omega^{2(m-1)}\\
  \vdots & \vdots & \ddots & \ddots &\vdots\\
  1 & \omega^{m-1} & \omega^{2(m-1)} &\cdots & \omega^{(m-1)^2}
  \end{pmatrix},
\end{equation}

\noindent and 

{\allowdisplaybreaks
\begin{equation}\label{this25}
S^{-1} = 
  \frac{1}{m}
  \begin{pmatrix} 
  1 & 1 &\cdots & 1 & 1\\
  1 & \omega^{-1}   & \omega^{-2} &\cdots & \omega^{-(m-1)}\\
  1 & \omega^{-2} & \omega^{-4} &\cdots & \omega^{-2(m-1)}\\
  \vdots & \vdots & \ddots & \ddots &\vdots\\
  1 & \omega^{-(m-1)} & \omega^{-2(m-1)} &\cdots & \omega^{-(m-1)^2}
  \end{pmatrix}.
\end{equation}
}

In the present cyclic, and hence, doubly stochastic case, we know
that $\Sigma = (1/m)(I + S^{-1}DS + (S^{-1}DS)^T)$,
where, as usual, 
$D =\text{diag}(\gamma_1,\gamma_2,\dots, \gamma_m)$
$= \text{diag}(-1/2,\lambda_2/(1-\lambda_2),\dots, \lambda_m/(1-\lambda_m))$.
We can then compute the entries of
$S^{-1}DS$ as follows:

{\allowdisplaybreaks
\begin{align}\label{this26}
(S^{-1}DS )_{j_1,j_2}
  &=  \sum_{k,\ell} (S^{-1})_{j_1,k} (D)_{k,\ell} (S)_{\ell,j_2}\nonumber\\
  &=  \sum_{k,\ell}  \frac{1}{m} (\omega^{-(j_1-1)(k-1)}) (\delta_{k,\ell} \gamma_k) (\omega^{(j_2-1)(\ell-1)})\nonumber\\
  &=  \sum_{k=1}^m   \frac{\gamma_k}{m} \omega^{(j_2-j_1)(k-1)}\nonumber\\
  &=  \frac{1}{m} \left(-\frac{1}{2} + \sum_{k=2}^m   \gamma_k \omega^{(j_2-j_1)(k-1)}\right),
\end{align}
}

\noindent for all $1 \le j_1, j_2, \le m$.
The entries of the asymptotic covariance matrix $\Sigma$
can thus be written as

\begin{align}\label{this27}
\sigma_{j_1,j_2}
  &= \frac{1}{m}\left(\delta_{j_1,j_2} + (S^{-1}DS )_{j_1,j_2} + (S^{-1}DS )_{j_2,j_1}\right)\nonumber\\
  &= \frac{1}{m}\left(\delta_{j_1,j_2} + \frac{1}{m} \left(-1 + \sum_{k=2}^m  \gamma_k ( \omega^{(j_2-j_1)(k-1)} + \omega^{(j_1-j_2)(k-1)} ) \right) \right)\nonumber\\
  &= \frac{m-1}{m^2} M^{(1)}_{j_1,j_2} + \frac{2}{m^2} \sum_{k=2}^m  \gamma_k  \beta_{|j_2-j_1|(k-1)} ,
\end{align}

\noindent for all $1 \le j_1, j_2, \le m$.

Next, note that since $\lambda_{m+2-k} = \bar {\lambda_k}$,
{\it i.e.}, the complex conjugate of $\lambda_{m+2-k}$,
we have $\gamma_{m+2-k} = \bar {\gamma_k}$,
for all $2 \le k \le m$.
Moreover, since 
$\beta_{|j_2-j_1|(k-1)} = \beta_{|j_2-j_1|((m+2-k)-1)}$,
we can write \eqref{this27} more symmetrically as
\eqref{this22} or \eqref{this23}, depending on whether
$m$ is odd or even, respectively, and in the latter case,
we also use that $\gamma_{m_0+1}$ is real, since 
$\omega^{m_0} = -1$.\CQFD
\end{Proof}

Let us again examine the cases $m=3$ and $m=4$.
In the former case, we have  

\begin{equation*}
M^{(1)} = 
\begin{pmatrix} 
  1 & -1/2 & -1/2\\
  -1/2 & 1 & -1/2\\
  -1/2 & -1/2 &1
  \end{pmatrix}.
\end{equation*}

\noindent But for $m=3$, $\beta_1 = -1/2 = \beta_2$, 
and so $M^{(1)} = M^{(2)}$.  Hence

\begin{equation}\label{Sig3}
\Sigma = \frac{2}{9} M^{(1)} + \frac{4}{9}Re(\gamma_2)M^{(2)} = \frac{2}{9}(1 + 2Re(\gamma_2))M^{(1)}.
\end{equation}

Hence, for $m=3$, cyclicity {\it always} produces a rescaled
version of the uniform iid case, with the rescaling
factor given by $1 + 2Re(\gamma_2)$.

For $m=4$, however,

\begin{equation*}
M^{(1)} = 
\begin{pmatrix} 
  1 & -1/3 & -1/3 &-1/3\\
  -1/3 & 1 & -1/3 &-1/3\\
  -1/3 & -1/3 &1  &-1/3\\
  -1/3 & -1/3 & -1/3 &1
  \end{pmatrix},
\end{equation*}

\noindent and $\beta_1 = 0$,
$\beta_2 = -1$, and $\beta_3 = 0$.  Thus,

\begin{equation*}
M^{(2)} = 
  \begin{pmatrix} 
    1 & 0 & -1 & 0\\
    0 & 1 & 0 & -1\\
    -1 & 0 & 1 & 0\\
    0 & -1 & 0 & 1
  \end{pmatrix},
\end{equation*}

\noindent and 
\begin{equation*}
M^{(3)} = 
  \begin{pmatrix} 
    1 & -1 & 1 & -1\\
    -1 & 1 & -1 & 1\\
    1 & -1 & 1 & -1\\
    -1 & 1 & -1 & 1
 \end{pmatrix}.
\end{equation*}

\noindent In this case, we have

\begin{align*}
\Sigma = \frac{3}{16} M^{(1)} + \frac{4}{16}Re(\gamma_2)M^{(2)} + \frac{2}{16}\gamma_3M^{(3)}.
\end{align*}

\noindent Next, note that
$2M^{(2)} + M^{(3)} = 3M^{(1)}$.
Then, if $Re(\gamma_2) = \gamma_3$,

\begin{align}\label{Sig4}
\Sigma &= \frac{3}{16} M^{(1)} + \frac{4}{16}Re(\gamma_2)M^{(2)} + \frac{2}{16}\gamma_3M^{(3)}\nonumber\\
&= \frac{3}{16} M^{(1)} + \frac{3}{16}(2Re(\gamma_2)M^{(1)})\nonumber\\
&= \frac{3}{16}(1 + 2Re(\gamma_2))M^{(1)},
\end{align}

\noindent so that there is still a
rescaled version of the iid case in a non-iid 
cyclic setting.  Indeed, since we know that
$\lambda_2 = a_1 + ia_2 - a_3 - ia_4 = (a_1-a_3) + i(a_2-a_4)$
and $\lambda_3 = a_1 - a_2 + a_3 - a_4$, 
we find that 

$$Re(\gamma_2) = \frac{a_2+2a_3+a_4}{(a_2+2a_3+a_4)^2 + (a_2-a_4)^2} - 1,$$

\noindent and $\gamma_3 = 1/(2(a_2+a_4)) - 1$.
A short calculation then shows that
$Re(\gamma_2) = \gamma_3$ if and only if
$a_3^2 = a_2 a_4$.  We thus have a
complete characterization of all $4$-letter,
cyclic Markov chains whose Young diagrams have
the same limiting shape as the uniform iid case.
In particular, choosing $a_2=a_4=a$,
for some $0 < a < 1/3$, leads to $a_3 = a$
and $a_1 = 1 - 3a$.  If, moreover,
$a = 1/4$, we have again the iid uniform case.
For $a \ne 1/4$, however, we may view the
Markov chain as a ``lazy'' version
of the uniform iid case.

Note that the scaling factor in both \eqref{Sig3}
and \eqref{Sig4} is $1 + 2Re(\gamma_2)$.
The following theorem shows that, in fact,
such a scaling factor occurs for general $m$,
and gives a spectral characterization of 
cyclic transition matrices giving a permutation-symmetric $\Sigma$ and thus to an iid
limiting shape.

\begin{theorem}\label{thm7b}
Let $P$ be the $m \times m$ transition matrix
of an aperiodic, irreducible, cyclic Markov chain
on an $m$-letter, totally ordered alphabet given in
Proposition~\ref{prop2}.  Then the asymptotic
covariance matrix $\Sigma$ is a rescaled version
of the iid uniform covariance matrix $\Sigma_{iidu}:= ((m-1)/m^2)M^{(1)}$ if
and only if 

\begin{align}\label{this21b}
Re\left(\frac{\lambda_j}{1-\lambda_j}\right) = \gamma,  \qquad \text{for all $2 \le j \le m$},
\end{align}

\noindent for some real constant $\gamma$.
Moreover, the scaling is then given by

\begin{equation}\label{this21ba}
\Sigma = (1 + 2\gamma)\Sigma_{iidu}.
\end{equation}

\end{theorem}

\noindent \begin{Proof} We first claim
that the system of matrix equations

\begin{equation}\label{this23c}
\sum_{j=2}^{m} b_j M^{(j)} = M^{(1)}
\end{equation}

\noindent has a solution
$b_j = 1/(m-1)$, for all $2 \le j \le m$.
Indeed, revisiting \eqref{this27}, we can
express each $M^{(j)}$ as

\begin{align}\label{this23d}
M^{(j)} =  \tilde M^{(j)} + \tilde M^{(m-j+2)},
\end{align}

\noindent where 
$(\tilde M^{(j)})_{k,\ell} = \omega^{(j-1)(\ell-k)}/2$,
for all $1 \le k,\ell \le m$, so that
\eqref{this23c} becomes

\begin{align}\label{this23e}
M^{(1)} = \sum_{j=2}^{m} b_j \left(\tilde M^{(j)} + \tilde M^{(m-j+2)}\right) 
= \sum_{j=2}^{m} (b_j+b_{m-j+2}) \tilde M^{(j)} 
= \sum_{j=2}^{m}  \tilde b_j \tilde M^{(j)},
\end{align}

\noindent where $\tilde b_j := (b_j+b_{m-j+2})/2$, 
for $2 \le j \le m$.

Now, clearly,
each $\tilde M^{(j)}$ is cyclic, so that in solving
\eqref{this23e} we need only examine
the $m$ entries in the first rows of the 
matrices.  We can thus reduce \eqref{this23e}
to the  
$m \times (m-1)$ system of equations

\begin{equation}\label{this23f}
  \begin{pmatrix}
  1 & 1 &1 &\cdots & 1\\
  \omega &\omega^2 &\omega^3 &\cdots &\omega^{m-1}\\
  \omega^2 &\omega^4 &\omega^6 &\cdots &\omega^{2(m-1)}\\
  \vdots & \vdots  &\vdots  &\ddots &\vdots\\
  \omega^{m-1} &\omega^{2(m-1)} &\omega^{3(m-1)} &\cdots &\omega^{(m-1)^2}
  \end{pmatrix}
  \begin{pmatrix}
  \tilde b_2\\
  \tilde b_3\\
  \vdots\\
  \tilde b_{m}
  \end{pmatrix}
  =
  \begin{pmatrix}
  1\\
  \frac{-1}{m-1}\\
  \frac{-1}{m-1}\\
  \vdots\\
  \frac{-1}{m-1}
  \end{pmatrix}.
\end{equation}

\noindent Since each of the last $m-1$ rows of the matrix
in \eqref{this23f} sums to $-1$, it is clear that
$\tilde b_j = 1/(m-1)$ is a solution to
the system.  To see that this solution is, 
in fact, unique, consider the 
$(m-1) \times (m-1)$ sub-matrix consisting of
the last $m-1$ rows of the matrix in \eqref{this23f},
namely,

\begin{equation}\label{this23g}
  \begin{pmatrix}
  \omega &\omega^2 &\omega^3 &\cdots &\omega^{m-1}\\
  \omega^2 &\omega^4 &\omega^6 &\cdots &\omega^{2(m-1)}\\
  \vdots & \vdots  &\vdots  &\ddots &\vdots\\
  \omega^{m-1} &\omega^{2(m-1)} &\omega^{3(m-1)} &\cdots &\omega^{(m-1)^2}
  \end{pmatrix}.
\end{equation}

\noindent Now this matrix, which is very closely
related to the Fourier matrix which arises in discrete
Fourier transform problems, is in fact invertible,
and can be shown to have one eigenvalue of $-1$,
and $m-2$ eigenvalues of the form
$\pm \sqrt{m}$ and $\pm i\sqrt{m}$, so that the
modulus of the determinant is $m^{(m-2)/2} \ne 0$.
Thus, the solution $\tilde b_j = 1/(m-1)$ is
unique, and since $\tilde b_j = (b_j+b_{m-j+2})/2$,
for all $2 \le j \le m$, we conclude that
$b_j = 1/(m-1)$ is a solution to \eqref{this23c} as well, for all $2 \le j \le m$,
and the claim is proved.

We can now use Proposition~\ref{prop2} to simplify
the asymptotic covariance matrix decomposition
as follows:

{\allowdisplaybreaks
\begin{align}\label{this23h}
\Sigma 
&= \frac{m-1}{m^2} M^{(1)} + \frac{2}{m^2}\sum_{k=2}^{m}\gamma_k M^{(k)}\nonumber\\
&= \frac{m-1}{m^2} M^{(1)} + 2\gamma \frac{1}{m^2}\sum_{k=2}^{m}M^{(k)}\nonumber\\
&= \frac{m-1}{m^2} M^{(1)} + 2\gamma \frac{m-1}{m^2}M^{(1)}\nonumber\\
&= (1+2\gamma)\frac{m-1}{m^2} M^{(1)}\nonumber\\
&= (1+2\gamma)\Sigma_{iidu},
\end{align}
}

\noindent where $\gamma = Re(\gamma_j)$, 
for all $2 \le j \le m$.  If the real parts of
$\gamma_j$ are not all identical, then the
uniqueness of the solution of \eqref{this23c}
implies that no such simplification is possible,
and the theorem is proved.\CQFD
\end{Proof}

\begin{Rem}\label{nonvac}
To see that the condition in \eqref{this21b}
is not vacuous for any $m$, recall that
for $m=4$, the ``lazy'' chain
has the iid limiting shape.  This is true
for general $m$: if $a_2 = a_3 = \cdots = a_m = a$,
for some $0 < a < 1/(m-1)$, then 
$\lambda_j = 1-m)a$, for all $2 \le j \le m$.
Trivially, then, $\gamma_j = 1/ma - 1 := \gamma$,
for all $2 \le j \le m$, so that the conditions
of Theorem~\ref{thm7b} are satisfied, and the
scaling factor is given by
$1 + 2\gamma = (2 -ma)/ma$.
Even in the $m=4$ case, however, we saw that
there were other, more general, cyclic
transition matrices which gave rise to the
iid limiting distribution.
\end{Rem}

The previous proposition indicates precisely when we
may expect the limiting shape of a cyclic Markov
chain to be spectrum of the traceless GUE. 
Now the first-order behavior of all rows
of the Young diagrams is $n/m + O(\sqrt{n})$ for cyclic 
Markov chains.  Although this differs from the first-order
behavior in the non-uniform iid case, one may still 
ask whether the limiting shape for some cyclic Markov chains might
still be that of some non-uniform iid case.  In fact,
this can never occur: cyclicity ensures that the 
asymptotic covariance matrix is also cyclic, and thus
cannot be equal to the asymptotic covariance matrix
of any non-uniform iid case.\\

Another class of Markov chains whose 
asymptotic covariance matrices can be easily studied
is the class of {\it reversible} Markov chains, {\it i.e.},
those with transition matrices such that 
$\pi_i P_{i,j} = \pi_j P_{j,i}$, for all $1 \le i,j \le m$.
The following theorem describes the asymptotic covariance matrix of
such Markov chains and, in the doubly stochastic case,
gives necessary and sufficient conditions for recovering
a rescaled uniform iid asymptotic covariance matrix.

\begin{theorem}\label{thm7c}
Let $P$ be the transition matrix
of an aperiodic, irreducible, reversible Markov chain
on an $m$-letter, totally ordered alphabet 
${\cal A}_m = \{\alpha_1 < \alpha_2 < \cdots < \alpha_m\}$. 
Then the asymptotic covariance matrix $\Sigma$
is given by

\begin{equation}\label{this21c}
\Sigma = \Pi^{1/2}S^T(I + 2 D) S\Pi^{1/2},
\end{equation}

\noindent where $P = (S\Pi^{1/2})^{-1} \Lambda (S\Pi^{1/2})$ is
the diagonalization of $P$.  
If, moreover, $P$ is doubly stochastic, {\it i.e.,}
if $P$ is symmetric, then $\Sigma$ is a rescaled version of the 
iid uniform case with transition matrix $P_{uiid}$ if and only if 
$P = \alpha P_{uiid} + (1 - \alpha) I$,
for $0 < \alpha \le m/(m-1)$.
\end{theorem}

\noindent \begin{Proof} It is elementary that $P$ is
similar to a symmetric matrix $Q$.  Specifically,
$P = \Pi^{-1/2}Q\Pi^{1/2}$, where $\Pi$ is the diagonal
matrix containing the stationary distribution.  
Since $Q$ is symmetric, it is diagonalizable 
with orthogonal eigenvectors.  Writing
$Q = S^T \Lambda S$, we have
$P = (S\Pi^{1/2})^{-1} \Lambda (S\Pi^{1/2})$. 
But then

\begin{align}\label{this21ca}
\Sigma &= \Pi + \Pi(S\Pi^{1/2})^{-1} D(S\Pi^{1/2}) + \left(S\Pi^{1/2})^{-1} D(S\Pi^{1/2}\right)^T \Pi\nonumber\\
       &= \Pi + 2\Pi^{1/2}S^T D S\Pi^{1/2} \nonumber\\
       &= \Pi^{1/2}(I + 2S^T D S)\Pi^{1/2} \nonumber\\
       &= \Pi^{1/2}S^T(I + 2 D) S\Pi^{1/2}.
\end{align}

\noindent where $D = \text{diag}(-1/2,\lambda_2/(1-\lambda_2),\dots,\lambda_m/(1-\lambda_m))$.
Writing $\Lambda_{\Sigma} := (I + 2D)$, we find that

\begin{equation*}
\Lambda_{\Sigma} = 
\text{diag}\left(0,\frac{1+\lambda_2}{1-\lambda_2},\dots,\frac{1+\lambda_m}{1-\lambda_m}\right).
\end{equation*}

In the doubly stochastic case, \eqref{this21ca} 
constitutes the diagonalization of $\Sigma$, since
$\Pi = (1/m)I$. If $\Sigma = a\Sigma_{uiid}$, for some $a > 0$,
then the spectrum of $\Sigma$ must consist of $0$ and
$m-1$ identical positive eigenvalues.  Therefore,
$\lambda_2 = \lambda_3 = \cdots = \lambda_m = \lambda$, so that
$\Lambda = \text{diag}(1,\lambda,\lambda,\dots,\lambda)$.
Thus $\Lambda_{\Sigma} = (1 + \lambda)(I - \Lambda)/(1-\lambda)^2$, 
so that \eqref{this21ca} now becomes

\begin{align}\label{this21cb}
\Sigma  = \frac{1}{m}S^T \Lambda_{\Sigma} S = \frac{1+\lambda}{m(1-\lambda)^2} S^T (I - \Lambda) S = \frac{1+\lambda}{m(1-\lambda)^2}(I - P).
\end{align}

\noindent In particular, $\Sigma_{uiid} = (I - P_{uiid})/m$.
Since $\Sigma = a\Sigma_{uiid}$, for some $a > 0$, it 
follows that

\begin{equation}\label{this21cc}
P = a\frac{(1-\lambda)^2}{1+\lambda}P_{uiid} 
   + \left(1 - a\frac{(1-\lambda)^2}{1+\lambda}\right)I.
\end{equation}

\noindent After setting $\alpha = a(1-\lambda)^2/(1+\lambda)$,
and checking that $P$ is irreducible
for $0 < \alpha \le m/(m-1)$, 
the theorem is proved.\CQFD
\end{Proof}

\begin{Rem}\label{reversible}
Recall (see the discussion preceding Theorem~\ref{thm7b}) that in the $4$-letter cyclic case, where $P$
had initial column entries $a_1, a_2,a_3,$ and $a_4$, 
the criterion for obtaining a rescaling of the covariance matrix
in the iid case was that $a_3^2 = a_2 a_4$.  If we further demand
that $P$ be symmetric, then $a_2 = a_4$, 
and the criterion is refined to
$a_2 = a_3 = a_4 = a$, which is consistent with Theorem~\ref{thm7c},
as $P = 4aP_{uiid} + (1-4a)I$,
for $0 < 4a \le 4/3$.
\end{Rem}

Thus far we have expressed our limiting laws in terms
of Brownian functionals whose Brownian motions have a
non-trivial covariance structure arising directly from
the specific nature of the transition matrix.  It is of
interest to instead express the limiting laws in
terms of {\it standard} Brownian motions.

Since the asymptotic covariance matrix $\Sigma$ 
is non-negative definite, we can find an $m \times m$
matrix $C$ such that $\Sigma = CC^T$. 
Clearly, we then have

\begin{equation}\label{this40}
(\sigma_1 \tilde B^1(t),\sigma_2 \tilde B^2(t),\dots,\sigma_m \tilde B^m(t))^T
= C(B^1(t),B^2(t),\dots,B^m(t))^T,
\end{equation}

\noindent where $(B^1(t),B^2(t),\dots,B^m(t))^T$ is a standard,
$m$-dimensional Brownian motion.

In order to simplify notation, we will assume that $\tau(\ell) = \ell$,
for all $\ell$, and so write our main result \eqref{item7r} 
in Theorem \ref{thm4} as

\begin{align}\label{this41}
\frac{V^r_n -\nu_r n}{\sqrt{n}} &\Longrightarrow  \sum_{k=1}^{m_r} \sigma_{k} \tilde B^{k}(1) 
+ \max_{I_{r-m_r,d_r}} 
    \sum_{j=1}^{r-m_r} \sum_{\ell=j}^{(d_r + m_r - r+j)} \sigma_{m_r + \ell} \tilde B^{m_r + \ell}(\Delta t_{j,\ell}) := V^r_{\infty}.
\end{align}

\noindent Simply substituting \eqref{this40} into
\eqref{this41} immediately yields

\begin{align}\label{this42}
V^r_{\infty}
&=\sum_{k=1}^{m_r} \left(\sum_{i=1}^m C_{k,i} B^{i}(1) \right)  + \max_{I_{r-m_r,d_r}} 
             \sum_{j=1}^{r-m_r} \sum_{\ell=j}^{(d_r + m_r - r+j)} 
               \left(\sum_{i=1}^m C_{m_r + \ell,i} B^{i}(\Delta t_{j,\ell})\right)\nonumber\\
&=\sum_{i=1}^m \sum_{k=1}^{m_r} C_{k,i} B^{i}(1)  + \max_{I_{r-m_r,d_r}} 
             \sum_{i=1}^m  \sum_{j=1}^{r-m_r} \sum_{\ell=j}^{(d_r + m_r - r+j)} 
                C_{m_r + \ell,i} B^{i}(\Delta t_{j,\ell}).
\end{align}

Now the first term in \eqref{this42} is simply a Gaussian term
whose variance can be computed explicitly.  Unfortunately,
the maximal term does not in general succumb to any
significant simplifications.  However, in the 
iid case, we can further simplify \eqref{this42} in a
very satisfying way.  Indeed, since, in the iid case, we have
$\sigma_k^2 = \pi_k(1-\pi_k)$ and, for $k \ne \ell$,
$\sigma_{k,\ell} = - \pi_k \pi_{\ell}$, one can quickly
check that $C$ can be chosen so that
$C_{k,k} = \sqrt{\pi_k} -\sqrt{\pi_k}\pi_k$, and, for $k \ne \ell$,
$C_{k,\ell} = -\sqrt{\pi_{\ell}}\pi_k$.  
Moreover, for all $m_r + 1 \le k \le m_r + d_r$,
$\pi_k$ = $\pi_{m_r + 1} = \pi_r$.
Then, within the
maximal term, 
$C_{m_r + \ell,i} = \sqrt{\pi_r} - \pi_r\sqrt{\pi_r}$,
for $i = m_r + \ell$, and
$C_{m_r + \ell,i} = - \pi_{r}\sqrt{\pi_i}$,
for $i \ne m_r + \ell$.
With the convention that $\nu_0 = 0$, 
we can then express \eqref{this42} as 

{\allowdisplaybreaks
\begin{align}\label{this43}
V^r_{\infty}
&=\sum_{i=1}^{m_r}  \sqrt{\pi_i} B^{i}(1) 
  + \sum_{i=1}^m \sum_{k=1}^{m_r} (-\sqrt{\pi_{i}}\pi_k) B^{i}(1)\nonumber \\ 
& \qquad  +  \max_{I_{r-m_r,d_r}} \biggl\{
                 \sum_{j=1}^{r-m_r} \sum_{\ell=j}^{(d_r + m_r - r+j)} 
                        \sqrt{\pi_r} B^{m_r+\ell}(\Delta t_{j,\ell}) +  \sum_{i=1}^m  \sum_{j=1}^{r-m_r} \sum_{\ell=j}^{(d_r + m_r - r+j)} 
                (- \pi_{r}\sqrt{\pi_i}) B^{i}(\Delta t_{j,\ell})\biggr\}\nonumber\\
&=\sum_{i=1}^{m_r}  \sqrt{\pi_i} B^{i}(1) 
  - \sum_{i=1}^m \sqrt{\pi_{i}} B^{i}(1)  \sum_{k=1}^{m_r}  \pi_k \nonumber \\ 
& \qquad  +  \sqrt{\pi_r} \max_{I_{r-m_r,d_r}} \biggl\{
                 \sum_{j=1}^{r-m_r} \sum_{\ell=j}^{(d_r + m_r - r+j)} 
                        B^{m_r+\ell}(\Delta t_{j,\ell}) -  \sqrt{\pi_{r}}\sum_{i=1}^m 
                        \sqrt{\pi_{i}} \sum_{j=1}^{r-m_r} \sum_{\ell=j}^{(d_r + m_r - r+j)} 
                 B^{i}(\Delta t_{j,\ell})\biggr\}\nonumber\\
&=\biggl\{\sum_{i=1}^{m_r}  \sqrt{\pi_i} B^{i}(1) 
  - \nu_{m_r} \sum_{i=1}^m \sqrt{\pi_{i}} B^{i}(1) 
  - \pi_r (r-m_r)\sum_{i=1}^m \sqrt{\pi_{i}} B^{i}(1)\biggr\}\nonumber \\
& \qquad  + \sqrt{\pi_r} \max_{I_{r-m_r,d_r}} 
                 \sum_{j=1}^{r-m_r} \sum_{\ell=j}^{(d_r + m_r - r+j)} 
                        B^{m_r+\ell}(\Delta t_{j,\ell})\nonumber\\
&=\biggl\{\sum_{i=1}^{m_r}  \sqrt{\pi_i} B^{i}(1) 
    - \nu_{r} \sum_{i=1}^m \sqrt{\pi_{i}}  B^{i}(1)\biggr\}  + \sqrt{\pi_r} \max_{I_{r-m_r,d_r}} 
                 \sum_{j=1}^{r-m_r} \sum_{\ell=j}^{(d_r + m_r - r+j)} 
                        B^{m_r+\ell}(\Delta t_{j,\ell})\nonumber\\
&=\Bigl\{(1-\nu_{r})   \sum_{i=1}^{m_r}  \sqrt{\pi_i}B^{i}(1)
  - \nu_{r} \sum_{i=m_r+d_r+1}^m \sqrt{\pi_{i}} B^{i}(1)\Bigr\}\nonumber\\
&\qquad  + \sqrt{\pi_r}\Bigl\{- \nu_{r}\sum_{i=m_r+1}^{m_r+d_r} B^{i}(1) + \max_{I_{r-m_r,d_r}} 
                 \sum_{j=1}^{r-m_r} \sum_{\ell=j}^{(d_r + m_r - r+j)} 
                        B^{m_r+\ell}(\Delta t_{j,\ell})\Bigr\}.
\end{align}
}

Note that the first two Gaussian term of \eqref{this43}
are independent of the remaining two Gaussian-maximal expression
terms.

Following Glynn and Whitt\cite{GW} and Baryshnikov\cite{Ba},
who studied the Brownian functional

$$D_m = \max_{I_{1,m}} \sum_{\ell=1}^{m} B^{\ell}(\Delta t_{\ell}),$$

\noindent we define the following, more general, Brownian functional:

\begin{equation}\label{this44}
D_{r,m} := \max_{I_{r,m}} 
                 \sum_{j=1}^{r} \sum_{\ell=j}^{(m-r+j)} 
                        B^{\ell}(\Delta t_{j,\ell}),
\end{equation}

\noindent where $1 \le r \le m$.  Clearly, the maximal term in
\eqref{this43} has just such a form.  We also remark that, see \cite{BGH}, 
$D_{r,m}$ corresponds to the sum of the $r$ largest eigenvalues
of an $m \times m$ GUE matrix.

To better understand \eqref{this43}, we may, without much
loss in generality, focus on the first block, that is,
values of $r$ such that $m_r = 0$.  The first  Gaussian
term of \eqref{this43} thus vanishes, and, writing
$\pi_{max}$ for $\pi_r$, we have

{\allowdisplaybreaks
\begin{align}\label{this45}
V^r_{\infty} = -r\pi_{max}  \sum_{i=d_1+1}^m \sqrt{\pi_{i}} B^{i}(1) + \sqrt{\pi_{max}} 
     \left(-r \pi_{max} \sum_{i=1}^{d_1} B^{i}(1) + D_{r,d_r}\right).
\end{align}
}

\noindent  In the uniform iid case, the first Gaussian term of \eqref{this45}
itself vanishes, since $d_r = d_1 = m$, and we have

\begin{align}\label{this46}
V^r_{\infty} = \frac{1}{\sqrt{m}}
     \left(-\frac{r}{m}\sum_{i=1}^{m}  B^{i}(1) + D_{r,m}\right) := \frac{H_{r,m}}{\sqrt{m}}.
\end{align}

\noindent Furthermore, and still specializing
\eqref{this45} to $r=1$, 

{\allowdisplaybreaks
\begin{align}\label{this47}
\frac{LI_n - n\pi_{max}}{\sqrt{n}} 
&\Longrightarrow -\pi_{max}  \sum_{i=d_1+1}^m \sqrt{\pi_{i}} B^{i}(1) + \sqrt{\pi_{max}} 
     \left(-\pi_{max} \sum_{i=1}^{d_1} B^{i}(1) + D_{1,d_1}\right)\nonumber\\
&\quad = -\pi_{max}  \sum_{i=d_1+1}^m \sqrt{\pi_{i}} B^{i}(1)  + \sqrt{\pi_{max}} 
       \left(\frac{1}{d_1}-\pi_{max}\right) \sum_{i=1}^{d_1}B^{i}(1) + \sqrt{\pi_{max}}H_{1,d_1}.
\end{align}
}

One can easily compute the variance of the Gaussian terms
in \eqref{this47} to be $\pi_{max}(1/d_1 - \pi_{max})$.
For $d_1 = 1$, this becomes $\pi_{max}(1 - \pi_{max})$,
and since $H_{1,1} = 0$ a.s., the limiting distribution is
$N(0,\pi_{max}(1 - \pi_{max}))$.  For $d_1 = m$, the variance
of the Gaussian term vanishes, and we again recover 
\eqref{this46}, namely, $H_{1,m}/\sqrt{m}$.
Both of these results are consistent with those of the
authors' previous paper \cite{HL}.

\section{Connections to Random Matrix Theory}

For iid uniform $m$-letter alphabets, 
the limiting law of the shape of the Young diagrams 
corresponds to the joint distribution 
of the eigenvalues of an $m \times m$ 
matrix from the traceless GUE \cite{Jo}. 
In the non-uniform iid case, we further noted that
Its, Tracy, and Widom \cite{ITW1,ITW2} have 
described the limiting law of the length of the longest weakly 
increasing subsequences as that of the 
maximal eigenvalue  
of a random matrix consisting of independent diagonal
blocks, each of which is a matrix
from the GUE.  The size
of each block depends
upon the multiplicity
of the corresponding stationary probability.
In addition, there is a zero-trace
condition involving the stationary probabilities
on the composite matrix (see \cite{HX} for the RSK diagrams result).

As a first step in extending these
connections between Brownian functionals and
spectra of random matrices, recall the general case when 
the stationary probabilities are all distinct
(see Remark \ref{RMTconnect}).
Our Brownian functionals then have
no true maximal terms, so that the limiting
shape, $(R_{\infty}^1,R_{\infty}^2,\dots,R_{\infty}^m)$
is simply multivariate normal, with
covariance matrix $\Sigma$ (or, more precisely,
the matrix obtained by permuting the rows and columns 
of $\Sigma$ using $\tau$, the 
permutation of $\{1,2,\dots,m\}$
previously defined).  Trivially,
this limiting law corresponds to the spectrum of a diagonal
matrix whose elements are multivariate normal
with the same covariance matrix $\Sigma$.

We can see that this general result is consistent
with the iid non-uniform case having distinct
probabilities.  Indeed, each block is of size $1$,
and is rescaled so that the variance is 
$\pi_{\tau(i)}(1-\pi_{\tau(i)})$, for
$1 \le i \le m$.  
Because of this rescaling, instead of having
a generalized zero-trace condition, as in the
non-rescaled matrices used in \cite{ITW1,ITW2},
our condition is rather a true zero-trace
condition.  This zero-trace condition is clear,
since the covariance matrix for {\it any} iid
case (uniform and non-uniform alike) is that of a 
multinomial distribution with parameters
$(n=1;\pi_{\tau(1)},\pi_{\tau(2)},\dots,\pi_{\tau(m)})$,
and any $(Y_1,Y_2,\dots,Y_m)$ having such a 
distribution of course satisfies
$\sum_{i=1}^m Y_i = 1$, so that
$\textvar(\sum_{i=1}^m Y_i) = 0$, which implies 
the zero-trace condition for
 $(R_{\infty}^1,R_{\infty}^2,\dots,R_{\infty}^m)$.

Next, consider the case when each stationary 
probability has multiplicity no greater than $2$.  
One may conjecture that
$(R_{\infty}^1,R_{\infty}^2,\dots,R_{\infty}^m)$
is equal in law to the spectrum of
a direct sum of certain $2 \times 2$
and/or $1 \times 1$ random matrices. 
Specifically,
let $\kappa \le m$ be the number of distinct probabilities
among the stationary distributions. Then the
composite matrix consists of a direct sum
of $\kappa$ GUE matrices which are as follows.
First, the overall diagonal 
$(X_1, X_2, \dots, X_m)$ of the matrix has
a $N(0,\Sigma)$ distribution.
Next, if $d_r = 1$, then the GUE matrix is
simply the $1 \times 1$ matrix $(X_r)$.
Finally, if $d_r=2$, then the GUE matrix
is the $2 \times 2$ matrix

\begin{equation*}
  \begin{pmatrix}
     X_{m_r+1} &Y_{m_r+1} + iZ_{m_r+1}\\
     Y_{m_r+1} - iZ_{m_r+1}   &X_{m_r+2}
  \end{pmatrix},
\end{equation*}

\noindent whose off-diagonal random variables 
$Y_{m_r+1}$ and $Z_{m_r+1}$
are iid, centered, normal random variables,
independent of all other random variables in the
overall matrix, 
with variance

\begin{equation*}
(\sigma^2_{m_r+1} 
   - 2\rho_{m_r+1,m_r+2}\sigma_{m_r+1}\sigma_{m_r+2}
   +\sigma^2_{m_r+2})/4.
\end{equation*}

\noindent Such a conjecture would imply the following
marginal result regarding a single block of such a matrix,
which without loss of generality we take to
be the first block.  In fact, this result is a genuine extension of the 
connection between maximal Brownian functionals and random matrices
beyond the standard GUE, traceless GUE, and purely diagonal cases
already discussed.

\begin{theorem}\label{thmgentwobytwo}
Let $d_1 = 2$ and $\tau(r) = r$, for all $1 \le r \le m$.
Then $(R_{\infty}^1,R_{\infty}^2) = (V_{\infty}^1,V_{\infty}^2-V_{\infty}^1)$
is distributed as the spectrum 
$\{(\lambda_1,\lambda_2): \lambda_1 \ge \lambda_2\}$ of the
$2 \times 2$ Gaussian Hermitian matrix

\begin{equation}\label{rmt1}
M :=
  \begin{pmatrix}
     X_1 &Y_1 + iZ_1 \\
     Y_1 - iZ_1   &X_2
  \end{pmatrix},
\end{equation}

\noindent where $(X_1,X_2)$ is a pair of centered,
bivariate normal random variables with covariance matrix

\begin{equation}\label{rmt2}
\Sigma_2 =
  \begin{pmatrix}
     \tilde{\sigma}^2_1              &\tilde{\rho}\tilde{\sigma}_1\tilde{\sigma}_2\\
     \tilde{\rho}\tilde{\sigma}_1\tilde{\sigma}_2      &\tilde{\sigma}^2_2
  \end{pmatrix},
\end{equation}
\noindent and $Y_1$ and $Z_1$ are iid, centered normal
random variables, independent of $(X_1,X_2)$, with
variance  

\begin{equation}\label{offdiagvar}
(\tilde{\sigma}^2_{1}  - 2\tilde{\rho}\tilde{\sigma}_{1}\tilde{\sigma}_{2} +\tilde{\sigma}^2_{2})/4.
\end{equation}

\end{theorem}

\noindent \begin{Proof} 
We prove the equivalent result that
$(V_{\infty}^1,V_{\infty}^2)$
is distributed as the pair
$(\lambda_1,\lambda_1+\lambda_2)$. Now, by definition,


\begin{align}\label{rmt3}
(V_{\infty}^1,V_{\infty}^2)
&= \Bigl( \max_{0\le t \le 1} 
     \bigl(\tilde{\sigma}_1\tilde{B}^1(t) 
       + \tilde{\sigma}_2\tilde{B}^2(1) - \tilde{\sigma}_2\tilde{B}^2(t) \bigr),  \tilde{\sigma}_1\tilde{B}^1(1) + \tilde{\sigma}_2\tilde{B}^2(1)  \Bigr)\nonumber\\
&= \Bigl(\tilde{\sigma}_2\tilde{B}^2(1)
     + \max_{0\le t \le 1} 
     \bigl(\tilde{\sigma}_1\tilde{B}^1(t) 
       - \tilde{\sigma}_2\tilde{B}^2(t) \bigr), \tilde{\sigma}_1\tilde{B}^1(1) + \tilde{\sigma}_2\tilde{B}^2(1)  \Bigr).
\end{align}

We now simplify \eqref{rmt3}, by
introducing new Brownian motions and then
decomposing the resulting expression into
two independent parts.  To do so, begin by
defining the new variances and correlation coefficients
$\sigma_1^2 :=  \tilde{\sigma}^2_2$,
$\sigma_2^2 :=  \tilde{\sigma}^2_1 - 
2\tilde{\rho}\tilde{\sigma}_1\tilde{\sigma}_2 + \tilde{\sigma}^2_2$,
and
$\rho := (\tilde{\rho}\tilde{\sigma}_1 - \tilde{\sigma}_2)/
\sqrt{\tilde{\sigma}^2_1 - 
2\tilde{\rho}\tilde{\sigma}_1\tilde{\sigma}_2 + \tilde{\sigma}^2_2}$.
Then it is easily verified that 
$B^1(t) := \tilde B^2(t)$, and
$B^2(t) := (\tilde \sigma_1 \tilde B^1(t) - \tilde \sigma_2 \tilde B^2(t))/\sigma_2 $,
where the equalities are pointwise,
are (dependent) standard Brownian motions, and
\eqref{rmt3} becomes

\begin{align}\label{rmt4}
(V_{\infty}^1,V_{\infty}^2)
&= \bigl(\sigma_1 B^1(1) + \sigma_2 \max_{0\le t \le 1} B^2(t),
          2\sigma_1 B^1(1) + \sigma_2 B^2(1)\bigr)\nonumber\\
&= \Bigl( (\sigma_1 B^1(1) - \rho \sigma_1 B^2(1)) 
         + \sigma_2 \Bigl(\rho \frac{\sigma_1}{\sigma_2}B^2(1) +  \max_{0\le t \le 1} B^2(t) \Bigr),\nonumber\\
&\qquad \qquad 2(\sigma_1 B^1(1) - \rho \sigma_1 B^2(1))  
              + (\sigma_2 + 2\rho \sigma_1) B^2(1)\Bigr).
\end{align}

\noindent Note that 
$B^1(t) - \rho B^2(t)$ is independent of
$B^2(t)$ and has variance 
$\sigma_1^2(1-\rho^2)t$.
Introducing the Brownian functional

\begin{equation}\label{rmt4b}
U(\beta) = \left(\beta - \frac12\right)B^2(1) + \max_{0\le t \le 1} B^2(t),
\end{equation}

\noindent $\beta \in \bbr$,
and using 
$\sigma_1^2, \sigma_2^2$, and $\rho$ above,
\eqref{rmt4} becomes

\begin{align}\label{rmt5}
&(V_{\infty}^1,V_{\infty}^2)
\stackrel{{\cal L}}{=} \sigma_1 \sqrt{1-\rho^2}Z (1,2)
 + \biggl(\sigma_2 U\left(\frac{1}{2} - \rho \frac{\sigma_1}{\sigma_2} \right),(\sigma_2 + 2\rho \sigma_1) B^2(1) \biggr)\nonumber\\
&\quad= \frac{\tilde{\sigma}_1\tilde{\sigma}_2\sqrt{1-\tilde \rho^2}}
        {\sqrt{\tilde{\sigma}^2_1 - 2\tilde{\rho}\tilde{\sigma}_1\tilde{\sigma}_2 + \tilde{\sigma}^2_2}}
        Z(1,2)\nonumber\\
&\qquad + \biggl(\sqrt{\tilde{\sigma}^2_1 - 2\tilde{\rho}\tilde{\sigma}_1\tilde{\sigma}_2 + \tilde{\sigma}^2_2}\quad
             U\biggl(\frac{\tilde{\sigma}^2_1 - \tilde{\sigma}^2_2}
                         {2\sqrt{\tilde{\sigma}^2_1 - 2\tilde{\rho}\tilde{\sigma}_1\tilde{\sigma}_2 
                          + \tilde{\sigma}^2_2} }  \biggr),  2(\tilde{\sigma}^2_1 - \tilde{\sigma}^2_2) B^2(1) \biggr),
\end{align}

\noindent where $Z$ is a standard normal
random variable independent of the sigma-field
generated by $B^2$.

Turning now to the eigenvalues' distributions,
we first consider the centered, multivariate normal
random variables $(W_1,W_2)$, having covariance matrix

\begin{equation*}
  \begin{pmatrix}
     \sigma^2_1              &\rho\sigma_1\sigma_2 \\
     \rho\sigma_1\sigma_2      &\sigma^2_2
  \end{pmatrix},
\end{equation*}

\noindent and let $W_3$ and $W_4$ be two iid, 
centered, normal random variables, independent of
$(W_1,W_2)$, with variance $\sigma_2^2$.
Then it is classical that

\begin{equation*}
\Bigl(W_2, \sqrt{W_2^2 + W_3^2 + W_4^2}\Bigr) \stackrel{{\cal L}}{=}
\sigma_2\bigl(B(1),2\max_{0\le t \le 1}B(t)-B(1)\bigr),
\end{equation*}

\noindent or, equivalently,

\begin{equation}\label{rmt6}
\biggl(W_2, \beta W_2 + \frac12 \sqrt{W_2^2 + W_3^2 + W_4^2}\biggr) \stackrel{{\cal L}}{=}
\sigma_2(B(1),U(\beta)),
\end{equation}

\noindent where $B$ is a standard Brownian
motion, and $U(\beta)$, $\beta \in \bbr$,
is defined in terms of $B$, rather than
in terms of $B^2$, as in \eqref{rmt4b}.
Then consider the random variable

\begin{align}\label{rmt7}
\tilde \lambda 
&:= W_1 + \sqrt{W_2^2 + W_3^2 + W_4^2}\nonumber\\
&= \biggl(W_1 - \rho\frac{\sigma_1}{\sigma_2}\biggr) 
       + \biggl(\rho\frac{\sigma_1}{\sigma_2} + \sqrt{W_2^2 + W_3^2 + W_4^2}\biggr).
\end{align}

\noindent Using \eqref{rmt6},
and noting that the variance of the first term
in \eqref{rmt7} is $\sigma_1^2(1-\rho^2)$,
it is easy to see that

\begin{equation}\label{rmt8}
\tilde \lambda 
\stackrel{{\cal L}}{=} \sigma_1\sqrt{1-\rho^2}Z 
                         + 2\sigma_2 U\Bigl(\frac{\rho\sigma_1}{2\sigma_2}\Bigr),
\end{equation}

\noindent where $Z$ is a standard normal random variable
independent of $B$.

We now apply this result to the eigenvalues of the 
matrix $M$ in \eqref{rmt1}, namely, to

\begin{equation}\label{rmt9}
\lambda_1 
= \biggl(\frac{X_1 + X_2}{2}\biggr)
  + \sqrt{\biggl(\frac{X_1 - X_2}{2}\biggr) + Y_1^2 + Z_1^2}, 
\end{equation}

\noindent and

\begin{equation}\label{rmt10}
\lambda_2 
= \biggl(\frac{X_1 + X_2}{2}\biggr)
  - \sqrt{\biggl(\frac{X_1 - X_2}{2}\biggr) + Y_1^2 + Z_1^2}.
\end{equation}

\noindent Letting $W_1 = (X_1 + X_2)/2$, $W_2 = (X_1 - X_2)/2$,
$W_3 = Y_1$, and $W_4 = Z_1$, we have

\begin{align}\label{rmt11}
(\lambda_1, \lambda_1+\lambda_2) 
&= \biggl( W_1 + \sqrt{W_2^2 + W_3^2 + W_4^2}, 2W_1 \biggr)\nonumber\\
&=  \biggl( \biggl(W_1 -\hat \rho \frac{\hat \sigma_1}{\hat \sigma_2}W_2\biggr)
       + 2\biggl( \hat \rho \frac{\hat \sigma_1}{2\hat \sigma_2}W_2
                 + \frac12 \sqrt{W_2^2 + W_3^2 + W_4^2}\biggr),\nonumber\\ 
&\qquad \qquad 2\biggl(W_1 -\hat \rho \frac{\hat \sigma_1}{\hat \sigma_2}W_2\biggr)
       + 2\hat \rho \frac{\hat \sigma_1}{\hat \sigma_2}W_2\biggr)\nonumber\\
&= \biggl(W_1 -\hat \rho \frac{\hat \sigma_1}{\hat \sigma_2}W_2\biggr)(1,2)\nonumber\\
&\qquad \qquad + \biggl(\hat \rho \frac{\hat \sigma_1}{2\hat \sigma_2}W_2
                           + \frac12 \sqrt{W_2^2 + W_3^2 + W_4^2},
               2\hat \rho \frac{\hat \sigma_1}{\hat \sigma_2}W_2\biggr),
\end{align}

\noindent where
$\hat \sigma_1^2 = (\tilde{\sigma}^2_1 
+2\tilde{\rho}\tilde{\sigma}_1\tilde{\sigma}_2 + \tilde{\sigma}^2_2)/4$,
$\hat \sigma_2^2 = (\tilde{\sigma}^2_1 
-2\tilde{\rho}\tilde{\sigma}_1\tilde{\sigma}_2 + \tilde{\sigma}^2_2)/4$,
and
$\hat \rho \hat \sigma_1^2 \hat \sigma_2^2 = (\tilde{\sigma}^2_1 - \tilde{\sigma}^2_2)/4$.
Noting that the variance of 
$W_1 - (\hat \rho \hat \sigma_1/ \hat \sigma_2) W_2$
is $\hat \sigma_1^2(1-\hat \rho^2) = \sigma_1^2(1-\rho^2)$,
and that, moreover,
$\beta := \hat \rho \hat \sigma_1 / 2\hat \sigma_2 = 
(\tilde{\sigma}^2_1 - \tilde{\sigma}^2_2) /
(2\sqrt{\tilde{\sigma}^2_1 - 2\tilde{\rho}\tilde{\sigma}_1\tilde{\sigma}_2 
                          + \tilde{\sigma}^2_2})$,
we find that

\begin{align}\label{rmt12}
(\lambda_1, \lambda_1+\lambda_2) 
&= \hat \sigma_1 \sqrt{1 - \hat \rho^2}Z(1,2)
   + \Bigl(2\hat \sigma_2 U\Bigl(\frac{\hat \rho \hat \sigma_1}{2 \sigma_2}\Bigr),
      2 \hat \rho \frac{\hat \sigma_1}{\hat \sigma_2}B^2(1) \Bigr)\nonumber\\
&= \sigma_1\sqrt{1- \rho^2}Z(1,2) + \sigma_2\bigl(U(\beta), 4\beta B^2(1) \bigr)\nonumber\\
&\stackrel{{\cal L}}{=} (V_{\infty}^1,V_{\infty}^2),
\end{align}

\noindent and we have our identity in law.\CQFD
\end{Proof}

The preceding theorem marks a first step in describing 
the increasingly robust connection between 
Gaussian random matrices and maximal functionals
of several Brownian motions.  We next turn to another very striking
relationship between a different class of $2 \times 2$ Gaussian random matrices and a maximal 
functional of two independent {\it Brownian bridges}.  As we will see, this result
is, in fact, more fundamental than the connections noted thus far in that it furnishes
a {\it conditional} equivalence in law from which many by-now classical results may be derived.

By a {\it standard Brownian bridge terminating at $b$} 
we mean a process $(\dot{B}(t;b))_{0 \le t \le 1}$ 
equal in law to $(B(t) - tB(1) + bt)_{0 \le t \le 1}$, where 
$B$ is a standard Brownian motion.

\begin{theorem}\label{brownbridgethm}
Let $(\dot{B}^1(t;b_1))_{0\le t \le 1}$ and $(\dot{B}^2(t;b_2))_{0\le t \le 1}$ be two independent 
standard Brownian bridges
terminating at $b_1$ and $b_2$, respectively, and let

\begin{equation}\label{brownbridge1}
M_{bb} :=
  \begin{pmatrix}
     b_1 &Y_1 + iZ_1 \\
     Y_1 - iZ_1   &b_2
  \end{pmatrix},
\end{equation}

\noindent where $b_1$ and $b_2$ are {\it real constants}, 
while $Y_1$ and $Z_1$ are independent centered normal random variables
with variance $1/2$.  Then the 
largest eigenvalue of $M_{bb}$ has the same law as

\begin{equation}\label{vbb}
V_{bb}:= \max_{0\le t \le 1}\bigl((\dot{B}^1(t;b_1) - \dot{B}^1(0;b_1)) 
+ (\dot{B}^2(1;b_2) - \dot{B}^2(t;b_2))\bigr).
\end{equation}

\noindent

\end{theorem}

\noindent \begin{Proof} We first simplify the familiar-looking expression
\eqref{vbb} and obtain

\begin{align}\label{brownbridge2}
V_{bb}
&:=       \max_{0\le t \le 1}\bigl((\dot{B}^1(t;b_1) - \dot{B}^1(0;b_1)) + (\dot{B}^2(1;b_2) - \dot{B}^2(t;b_2))\bigr)\nonumber\\
&=        \max_{0\le t \le 1}\bigl( \dot{B}^1(t;b_1) + (b_2 - \dot{B}^2(t;b_2))\bigr)\nonumber\\
&=  b_2 + \max_{0\le t \le 1}\bigl( \dot{B}^1(t;b_1) - \dot{B}^2(t;b_2)\bigr).
\end{align}

Since a linear combination of two independent Brownian bridges is 
again a Brownian bridge, by elementary rescaling we obtain the equality in law

\begin{equation}\label{brownbridge3}
V_{bb} \stackrel{{\cal L}}{=} b_2 + \sqrt{2} \max_{0\le t \le 1} \dot{B} \biggl(t; \frac{b_1 - b_2}{\sqrt{2}} \biggr),
\end{equation}

\noindent where $\dot{B}$ is a standard Brownian bridge
terminating at $(b_1-b_2)/\sqrt{2}$.

Next, recall the elementary result (which may be easily obtained from the Reflection
Principle) stating that the maximum $M(b)$ of a standard Brownian bridge terminating
at $b$ has a distribution given by

$$\bbp(M(b) \ge a) = \exp{(-2a(a-b))},$$

\noindent for any $a \ge \max(b,0).$

Applying this result to \eqref{brownbridge3}, we find that

{\allowdisplaybreaks
\begin{align}\label{brownbridge4}
\bbp(V_{bb} \ge a) 
&= \bbp\biggl( \max_{0\le t \le 1} \dot{B} \biggl(t; \frac{b_1 - b_2}{\sqrt{2}} \biggr) \ge \frac{a-b_2}{\sqrt{2}}    \biggr) \nonumber\\
&= \exp{\biggl( -2 \biggl(\frac{a-b_2}{\sqrt{2}}\biggr)  \biggl(\frac{a-b_2}{\sqrt{2}}  - \frac{b_1-b_2}{\sqrt{2}}  \biggr) \biggr) }\nonumber\\
&= \exp{(-(a-b_1)(a-b_2))},
\end{align}
}

\noindent for all $a \ge \max(b_1, b_2)$.

To complete the proof of the claim, we simply note that the 
eigenvalues of $M_{bb}$ satisfy the characteristic equation

$$(\lambda - b_1)(\lambda - b_2) - Y_1^2 - Z_1^2 = 0.$$

\noindent Writing $\lambda_1$ for the largest eigenvalue
of $M_{bb}$, we see that the event $\{ \lambda_1 \ge a \}$,  
for $a \ge \max(b_1, b_2)$, may be rewritten as the event
$\{  Y_1^2 + Z_1^2 \ge (a-b_1)(a-b_2) \}$.  But 
$Y_1^2 + Z_1^2$ is distributed as a (rescaled) 
Chi-squared random variable 
with two degrees of freedom and is in fact
exponential with parameter $1$.  Thus,

\begin{align}\label{brownbridge5}
\bbp(\lambda_1 \ge a) 
&= \bbp( Y_1^2 + Z_1^2 \ge (a-b_1)(a-b_2)) \nonumber\\
&= \exp{(-(a-b_1)(a-b_2))},
\end{align}

\noindent and our claim is proved.\CQFD
\end{Proof}

\begin{Rem}\label{bridgetosomewhere}
The significance of Theorem \ref{brownbridgethm} 
lies in the fact that we can de-condition $b_1$ and $b_2$ 
to obtain again the connection between 
maximal Brownian functionals and both the GUE or the 
traceless GUE, as shown next.  Indeed, replace $b_1$ and $b_2$ 
by iid standard normal random variables $X_1$ and $X_2$ 
which are independent of $(Y_1,Z_1)$ in the random matrix setting.
Then, the joint density of $(V_{bb},X_1,X_2)$ is given by

\begin{equation}\label{decon1}
g_3(v,x_1,x_2) = (2v - x_1 - x_2)\varphi(x_1)\varphi(x_2)\exp{(-(v-x_1)(v-x_2))},
\end{equation}

\noindent for $x_1, x_2 \in \bbr$ and $v \ge \max(x_1, x_2)$, where $\varphi$ is the standard 
normal density.  

Let us first examine the constant-trace case by letting 
$w_1 = x_1 + x_2$ and $w_2 = x_1 - x_2$
and then conditioning the joint law in \eqref{decon1} on $w_1$.  An elementary
calculation shows that the conditional density is simply

\begin{equation}\label{decon2}
h(v|w_1) = \frac{4}{\sqrt{\pi}} \Bigl(v - \frac{w_1}{2}\Bigr)^2 \exp{\Bigl(v - \frac{w_1}{2}\Bigr)^2},
\end{equation}

\noindent for $v \ge w_1/2$.  In particular, for $w_1 = 0$, we do indeed recover
the traceless case:

\begin{equation}\label{decon3}
h(v|w_1=0) = \frac{4}{\sqrt{\pi}} v^2 e^{-v^2},
\end{equation}

\noindent for $v \ge 0$, which is the density associated with
the square root of a (rescaled) Chi-squared random varible with
three degrees of freedom, as expected.

Next, let us return to the joint density in \eqref{decon1} and successively integrate out
$x_1$ and $x_2$.  Writing $\Phi$ for the standard normal
cumulative distribution function, we find that the joint density of $(V_{bb},X_2)$
is given by

\begin{equation}\label{decon4}
g_2(v,x_2) = \varphi(v)\left(\varphi(x_2) +  v\Phi(x_2)\right),
\end{equation}

\noindent for $-\infty < x_2 \le v < +\infty$.  Integrating once more, we obtain
the density of the largest eigenvalue of the $2 \times 2$ GUE,

\begin{equation}\label{decon5}
g_1(v) = \varphi(v)\left( (1 + v^2)\Phi(v) + v\varphi(v)\right),
\end{equation}

\noindent for all $v \in \bbr$.  We have thus recovered the traceless and 
standard GUE distributions from our conditional result.

\end{Rem}

The results above are all ultimately based
on the trivial fact that we can solve a quadratic equation.
While larger random matrices do not in general yield
to such straightforward analyses, there are 
nonetheless certain classes that can be linked to distributions of maximal Brownian functionals.
The following development consolidates several of these and leads
to a more general conjecture.

Let $B$ be  
$m$-dimensional standard Brownian motion and 
introduce the notation

\begin{equation*}
 \text{diag}(A)= \begin{pmatrix} A_{1,1}\\ \vdots \\ A_{m,m}\end{pmatrix},
\end{equation*}

\noindent for any $m\times m$ matrix $A$ and

\begin{equation*}
\text{Diag}\begin{pmatrix} v_1\\ \vdots\\ v_m\end{pmatrix}=
\begin{pmatrix} v_1 &&  0\\ &\ddots &\\ 0 && v_m\end{pmatrix},
\end{equation*}

\noindent for any vector $v = (v_1,\dots,v_m)^ T$.
Now recall the Brownian functional introduced by Glynn and Whitt \cite {GW},
namely,

\begin{equation*}
D_m=\max\sum^m_{k=1} B^k(\udel t_k),   
\end{equation*}
\noindent where the maximum is taken over the usual Weyl chamber, and from results of 
Baryshnikov  \cite{Ba} and Gravner, Tracy, and Widom \cite {GTW}, 
$D_m\stackrel \call = \lam_{\max}(M)$,  
for $M$ an $m \times m$ element of the GUE.  
To extend this identity to a broader class of random matrices, first let $A = aI + {1}_m b^T$, 
where $a\in\bbr$, $b\in\bbr^m$, 
and where ${1}_m = (1,1,\dots,1)^T\in\bbr^m$.  Next, letting $\tilde B=AB$, the corresponding 
maximal Brownian functional associated with $\tilde B$ is such that 

\begin{align*}
\tilde D_m &:= \max \sum_{k=1}^m \tilde B^k(\udel t_k)\\
&=\max \sum_{k=1}^m\left( aB^k(\udel t_k)+\sum_{j=1}^m b_jB^j(\udel t_k)\right)\\
&= aD_m+\sum_{j=1}^m b_j\sum_{k=1}^m B^j(\udel t_k)\\
&= aD_m+b^TB(1)\\
&\stackrel\call = a\lam_{\max}(M) + b^T\diag (M)\\
&= \lam_{\max}(aM+b^T\diag (M)I)\\
&= \lam_{\max}(\tilde M), 
\end{align*}
\noindent where the equality in law follows from the results of \cite{BGH}, and where $\tilde M=aM+b^T\diag (M)I$.\\

This rank-one perturbation result, 
when combined with Theorem~\ref{thm5} or Corollary~\ref{cor1a} and taking 
$a=1$ and $b_i = -1/m$, $i=1, \dots, m$, recovers \cite{TW} only     
using combinatorial and probabilistic techniques and therefore bypassing the analytical ones 
present there.  It also shares, as shown next, its amenability to analysis 
with two further cases: the cyclic case for $m = 2$ or $3$, and also the case 
$A=aP$, where $P$ is an arbitrary permutation matrix and $a$ an arbitrary real.  
Write the general $m\times m$ cyclic matrix
$A$ as $A=\sum^{m-1}_{k=0}a_kP^k$, 
where $P$ is the permutation matrix
$P=\begin{pmatrix} 
0 & &  &1\\
1 & 0  &\bigcirc\\
& 1 & \ddots &\\
\bigcirc && 1 & 0\end{pmatrix}$.  
Such a cyclic $A$ also yields a cyclic covariance matrix
$\Sigma := AA^T$, and moreover, 
$$\tilde D_m := \max\sum_{k=1}^m\sum_{j=1}^m A_{kj}B^j(\udel t_k)=\max
\sum_{k=1}^m\sum_{j=1}^m a_{k-j} B^j(\udel t_k)=\max\sum^{m-1}_{r=0} a_n\sum^m_{k=1}B^{k-r}
(\udel t_k).$$  We claim that for all three cases above,
$$\tilde D_m \stackrel\call=\lam_{\max}\left(f(\sigma(A))M_0+
\begin{pmatrix}X_1 && 0\\
& \ddots\\
0 & & X_m\end{pmatrix}\right),$$ 
where $f$ is an as-yet unspecified function 
of, $\sigma(A)$, the spectrum of
$A$, $(X_1,\dots, X_m)^T\sim N(0,\Sigma)$, $\Sigma=AA^T$, 
and $M_0$ is a `` GUE matrix" with all diagonal entries of $0$.  Writing 
$\tilde M =f(\sigma(A))M_0+\begin{pmatrix} X_1 && 0\\
&\ddots\\ 0 & & X_m\end{pmatrix}$, we then have $\tilde D_m \stackrel\call=\lam_{\max}(\tilde M)$.  
Let us first look at the simple case of $A=aP$, 
where $P$ is an arbitrary permutation matrix with 
$P_{i,j}=\del_{i,\sigma(j)}$, $\sigma$ a permutation of $\{1,2,\dots,m\}.$ 
Then $\tilde B=aPB\stackrel\call= aB$, and so
$\tilde D_m=\max\sum^m_{k=1}aB^{\sigma (k)} (\udel t_k)\stackrel \call=
aD_m$, while $\tilde M=aM_0+\begin{pmatrix} X_1 && 0\\ &\ddots \\
0 && X_m\end{pmatrix}$, where $(X_1,\dots, X_m)^T\sim N(0,\Sigma)$, 
$\Sigma=AA^T=a^2I$.  But since $\tilde M=aM_0 + a\udiag (\diag (M)) = aM$, 
$M\in GUE$, we simply recover the base case $\tilde D_m\stackrel \call=
a\lam_{\max}(M)$ and see that $f(\sigma (A))=a$. \\

Next, consider the degenerate cyclic case where $a_0=a_1=\cdots = a_{m-1}=a$,
so that $A=a1_m 1_m^T$. But this is just a special case of the rank-one perturbation 
 $A = aI + {1}_m b^T$, so that
$\tilde D_m\stackrel\call= a1_m^T B(1)\sim N(0,ma^2)$, and 
$\tilde M=0M+a1_m^T(X_1,\dots, X_m)^TI=a\sum_{k=1}^m X_kI$, so that
$\lam_{\max} (\tilde M)=a\sum_{k=1}^m X_k\sim N(0,ma^2)$. In this pure Gaussian
case, $f(\sigma (A))=0$.\\

More generally, as previously shown, for $A=aI+{1}_m b^T$, 
\begin{align*}
\tilde D_m&\stackrel \call=\lam_{\max}(a_M+(b^T1_m)\diag (M)I)\\
&= \lam_{\max}(aM_0+(a+b^T1_m)\diag(M)I)\\
&=\lam_{\max}(aM_0+A\diag (M)I), 
\end{align*}
\noindent and since $A\diag (M)\sim N(0,\Sigma)$, $f(\sigma (A))=a$. \\

Let us now look at the low-dimensional cyclic cases in turn.
For the $2\times 2$ cyclic case,
let $A=a_0I+a_1\begin{pmatrix}0 & 1\\ 1 & 0\end{pmatrix}=(a_0-a_1)I+
a_1 1_2 1_2^T$. Then 

\begin{align*}
\tilde M&=(a_0-a_1)M+a_1 1_2^T\diag (M)I\\
& =(a_0-a_1)M_0+(a_0-a_1)\udiag(\diag (M))+a_1 1_2^T\diag (M)I\\
&= (a_0-a_1)M_0 + \udiag (((a_0-a_1)I+a_1 1_2^T1_2)\diag (M))\\
&= (a_0-a_1)M_0 + \udiag (A\diag (M)), 
\end{align*}
and here $A\diag (M)\sim N(0,\Sigma)$ and so $f(\sigma (A))=a_0-a_1$.  
Finally, consider the $3$-dimensional cyclic case.  Let 
$A=\begin{pmatrix} a_0 & a_2 & a_1\\ a_1 & a_0 & a_2\\
a_2 & a_1 & a_0\end{pmatrix}$, then  
$\Sigma = AA^T= \sum_{k=1}^3 a^2_k I+ \sum_{k\neq\ell} a_ka_{\ell}(1_3 1_3^T-I) 
= (\sum_{k=1}^3 a^2_k-\sum_{k\neq\ell} a_ka_{\ell})I+(\sum_{k\neq\ell}a_ka_{\ell}){1}_3{1}_3^T$. But we can also write 
$\Sigma = (\alpha I+\beta 1_31_3^T) (\alpha I + \beta 1_3 1_3^T)^T=
\alpha^2I+(2\alpha\beta + 3\beta^2)1_3 1_3^T$. Thus $\alpha^2= \sum_{k=1}^3 a_k^2-
\sum_{k\neq\ell}a_ka_\ell$, therefore $3\beta^2+2\alpha\beta=\sum_{k\neq\ell}a_ka_{\ell}$ and so 
$\beta=-\alpha/3 \pm 1/3\sqrt{\sum_{k=1}^3 a^2_k-3\sum_{k\neq\ell} a_ka_{\ell}}$.  Now, take
$\alpha = \sqrt{\sum_{k=1}^3 a^2_k-\sum_{k\neq\ell} a_ka_{\ell}}$ and 
$\beta=-\alpha/3 + \sqrt{\sum_{k=1}^3 a^2_k-3\sum_{k\neq\ell} a_ka_{\ell}}/3$.  But, 
$\tilde D_3\stackrel\call=\lam_{\max}(\alpha M_0+\udiag(
\alpha I+\beta 1_3 1_3^T)\diag (M))$, so that
$\alpha I+\beta 1_3 1_3^T\diag (M)\sim N(0,\Sigma)$, where $f(\sigma (A))=\alpha$.  
Noting that $M_0$ and the matrix $\begin{pmatrix} X_1 && 0\\ &\ddots\\
0 && X_m\end{pmatrix}$ are independent, we can summarize all the above 
results in the following theorem:

\begin{theorem}\label{cyclicRMTlins}
Let $B=(B(t))_{0\le t \le 1}$ be a standard $m$-dimensional Brownian motion, $A$ 
an $m \times m$ real matrix, and $\tilde B = AB$, so that $\tilde B$ has
covariance matrix $\Sigma t= AA^T t$.  Then, in the cases noted in the table below,
the maximal functional $\tilde D_m= \max\sum^m_{k=1} \tilde B^k(\udel t_k)$ 
is equal in law to the largest eigenvalue of an $m \times m$ 
Hermitian Gaussian random matrix having diagonal elements 
distributed as $N(0,\Sigma)$ and off-diagonal ones with total variance given
by $f(\sigma(A))$, independent of each other and of the diagonal entries.

\bigskip
\begin{tabular}{|c|c|c|}
\hline
Case & $\sigma (A)$ & $f(\sigma (A))$\\
\hline
$m\times m$, $A=aP, a\in \bbr,$ $P$ perm. & $\{a,a,\dots, a\}$ & $a$\\
\hline
$m\times m$, $A=aI+1_mb^T$ & $\{a+b^T 1_m, a,\dots, a\}$ & $a$\\
\hline
$2\times 2$ cyclic $A=\sum^1_{k=0}a_k P^k$ & $\{ a_0 + a_1, a_0-a_1\}$
& $a_0-a_1$\\
\hline
$3\times 3$ cyclic $A=\sum^2_{k=0}a_kP^k$ & $\{a_0+a_1+a_2,$
& $\alpha=\sqrt{\sum_{k=1}^3 a^2_k-\sum_{k\neq\ell} a_ka_{\ell}}$\\
& $a_0+a_1u+a_2u^2$ & $= |a_0+a_1u+a_2u^2|$\\
& $a_0+a_1u^2+a_2u\},$  & $=|a_0+a_1u^2+a_2u|$\\
&  $u=e^{2i\pi/3}$ &\\
\hline
\end{tabular}
\bigskip

\end{theorem}

\section{Concluding Remarks}

In this paper, we have obtained the limiting shape of
Young diagrams generated by an aperiodic, irreducible,
homogeneous Markov chain on a finite state alphabet.
The following remarks indicate natural directions
in which our results in some cases can, 
and in other cases, may hope to, be extended.\\

\noindent \textbullet \quad  Our limiting 
theorems have all been proved assuming 
that the initial distribution is 
the stationary one. However, 
such results as Theorem 2 of Derriennic and Lin
\cite{DL} allow to extend our framework to initial
distributions started at a specified state.
Indeed, in this case, {\it i.e.,} if for
some $k=1,\dots,m$, $\bbp(X_0 = \alpha_k) = 1$,
the asymptotic covariance matrix is still given
by \eqref{item6la}, and, for example, Theorem \ref{thm4}
remains valid.  For an arbitrary initial distribution,
what is needed in this non-stationary context is
an invariance principle.  More generally, our
results continue to hold for $k^{th}$-order Markov chains,
and in fact, they extend to {\it any} sequence for which
both an asymptotic covariance matrix and an 
invariance principle exist.\\

\noindent \textbullet \quad Our limiting theorems 
have only been proved for finite alphabets.
However, from the authors' previous work \cite{HL},
it is known that for countably infinite iid alphabets, 
$LI_n$ has a limiting law corresponding to that 
of a non-uniform, finite-alphabet.  Hence,
for a countably infinite-alphabet Markov chain (subject
to additional constraints), 
we might still be able to obtain limiting laws 
of the form developed in this paper.\\

\noindent \textbullet \quad By using
appropriate existing concentration inequalities,
one can expect to establish
the convergence of the moments of
the rows of the diagrams.\\

\noindent \textbullet \quad 
One field in which the connection between Brownian
functionals and random matrix theory has been
exploited is in Queuing Theory.  The development
below, following O'Connell and Yor \cite{OY1},
shows how Brownian functionals of the sort we have
studied arise as generalizations of standard queuing
models.

Let $A(s,t]$ and $S(s,t]$,
$-\infty < s < t < \infty$,
be two independent Poisson
point process on $\bbr$, with intensity measures
$\lambda$ and $\mu$, respectively, with $0 < \lambda < \mu$.
Here $A$ represents the {\it arrivals} process, and $S$ 
the {\it service time} process, at a queue consisting
of a single server. The condition $\lambda < \mu$ 
ensures that the {\it queue length}

\begin{equation}\label{qitem1}
Q(t) = \sup_{-\infty < s \le t} \left\{ A(s,t] - S(s,t] \right\},
\end{equation}

\noindent is a.s.~finite, for any $t \in \bbr$.  
Then, defining the {\it departure} process

\begin{equation}\label{qitem2}
D(s,t] = A(s,t] - (Q(t)-Q(s)),
\end{equation}

\noindent which is simply the number of arrivals
during $(s,t]$ less the change in the queue length
during $(s,t]$,
the classical problem is to determine the distribution
of $D(s,t]$.  The answer to this problem is given
by {\it Burke's Theorem} \cite{Burke} (see Theorem $1$ of \cite{OY1}):

\begin{theorem}\label{Bkthm}
$D$ is a Poisson process with intensity $\lambda$, and
$\{D(s,t],s\le t\}$ is independent of $\{Q(s),s\ge t\}$.
\end{theorem}

That is, $D$ has the same law as the arrivals process $A$.
Moreover, since,the queue length after time $t$ is independent
of the process $D$ up to time $t$, 
one may take the departures from the first queue
and use them as inputs to a second queue, and observe that
the departure process from the second queue also has the law
of $A$.  Proceeding in this way, one generalizes to a {\it tandem queue} 
of $n$ servers, each taking the departures from the previous queue
as its arrivals process.

One can further generalize this model to a {\it Brownian queue
in tandem} in the following manner.  Let
$B,B^1,B^2,\dots, B^n$ be independent, standard Brownian motions
on $\bbr$,
and write $B^k(s,t) = B^k(t) - B^k(s)$, for each $k$ and $s < t$, 
and similarly for $B$.  Let $c > 0$ be a constant, and define,
in complete analogy to \eqref{qitem1} and \eqref{qitem2},

\begin{equation}\label{qitem3}
q_1(t) = \sup_{-\infty < s \le t} \left\{ B(s,t) + B^1(s,t) - c(t-s) \right\},
\end{equation}

\noindent and, for $s < t$,

\begin{equation}\label{qitem4}
d_1(s,t) = B(s,t) - (q_1(t)-q_1(s)).
\end{equation}

\noindent  For $k = 2,3,\dots,n$, let

\begin{equation}\label{qitem5}
q_k(t) = \sup_{-\infty < s \le t} \left\{ d_{k-1}(s,t) + B^k(s,t) - c(t-s) \right\},
\end{equation}

\noindent and, for $s < t$,

\begin{equation}\label{qitem6}
d_k(s,t) = d_{k-1}(s,t) - (q_k(t)-q_k(s)).
\end{equation}

Here $B$ is the arrivals process for the first queue, 
$d_{k-1}$ is the arrivals process for the $k^{th}$ queue ($k \ge 2)$,
and  $ct - B^k(t)$ is the service process for the $k^{th}$ queue,
for all $k$.
Using the ideas employed in Burke's Theorem, 
it can be shown that the generalized queue lengths
$q_1(0), q_2(0), \dots, q_n(0)$
are iid random variables.  Moreover, they are 
exponentially distributed with mean $1/c$.

Using the definitions in \eqref{qitem3}-\eqref{qitem6},
and a simple inductive argument, one finds that

\begin{equation}\label{qitem7}
\sum_{k=0}^n q_k(0) 
=\sup_{t > 0} 
        \Bigl\{ B(-t,0) - ct + L_n(t)\Bigr\},
\end{equation}

\noindent where

\begin{equation}\label{qitem8}
L_n(t) = \sup_{{\scriptstyle 0\le s_1\le\cdots}{\le s_{n-1}\le t} }
                 \{B^1(-t,-s_{n-1}) + \cdots + B^n(-s_1,0)\}.
\end{equation}

By Brownian rescaling, we observe that 

\begin{align}\label{qitem9}
L_n(t) 
&\stackrel{ {\cal L} }{=} \sqrt{t}\sup_{{\scriptstyle 0\le s_1\le\cdots}{\le s_{n-1}\le 1} }
                 \{B^1(-1,-s_{n-1}) + \cdots + B^n(-s_1,0)\} \stackrel{ {\cal L} }{=} \sqrt{t}V^1_{\infty},
\end{align}

\noindent where the functional $V^1_{\infty}$ is
as in Theorem \ref{thm4}, with associated $n \times n$
covariance matrix $\Sigma = tI_n$ and parameter set
$I_{1,n}$.  Thus, $L_n(t)$ may
be thought of as a process version of this $V^1_{\infty}$.

The generalized Brownian queues 
in \eqref{qitem3}-\eqref{qitem6}
involved
{\it independent} Brownian motions.  These can be
replaced by Brownian motions 
$B^1$,
$\dots $,
$B^n$ for which
$(\sigma_1 B^1(t),\dots ,\sigma_n B^n(t))$
has (nontrivial) covariance matrix $t\Sigma$.  
Whether or not we keep the initial 
arrival process $B$
independent of $(B^1$,$\dots$,$B^n)$, we 
no longer have that
$q_1(0), q_2(0), \dots, q_n(0)$
are iid random variables, 
due to the dependence among the service times
$ct-B^k(t)$, but we do still have the identity
\eqref{qitem7} and \eqref{qitem9} 
relating the total occupancy
of the queue at time zero to 
$V_{\infty}^1$.
More importantly, our generalizations of the Brownian
functionals $L_n(t)$ above can be used to describe
the joint law of the input/output of each queue.\\

\noindent \textbullet \quad An important topic connecting
much of random matrix theory to other problems, such
as the shape of random Young diagrams, is the field of
orthogonal polynomials. (See, {\it e.g.}, \cite{Jo}.)
It would be of great interest to see what, if any, classes of
orthogonal polynomials are associated with the present paper.\\

\end{document}